\documentclass[9pt]{elsarticle}

\usepackage{hyperref} % for autoref
\usepackage{amsmath} % for eqref
\usepackage{subcaption} % for subfigure
\usepackage{amsfonts} % for mathbb
\usepackage{bm} % for bm
\usepackage{bbm} % for indicator function
\usepackage{stmaryrd} % for \llbracket and \rrbracket 
\usepackage{makecell} % for multiple lines within a cell of table
\usepackage{multirow} % for multirow within a table
\usepackage[misc]{ifsym} % for symbol Letter
\usepackage{appendix} % for appendix
\usepackage{color} % for color
\usepackage{fullpage}
\usepackage{cancel}

\makeatletter
\newcommand{\pushright}[1]{\ifmeasuring@#1\else\omit\hfill$\displaystyle#1$\fi\ignorespaces}
\newcommand{\pushleft}[1]{\ifmeasuring@#1\else\omit$\displaystyle#1$\hfill\fi\ignorespaces}
\makeatother

\usepackage{mathtools}

\usepackage{algorithm}
\usepackage{algorithmic}
\usepackage{eqparbox}

\usepackage{etoolbox}  % patch def of algorithmic environment
\makeatletter
\patchcmd{\algorithmic}{\addtolength{\ALC@tlm}{\leftmargin} }{\addtolength{\ALC@tlm}{\leftmargin}}{}{}
\makeatother

\usepackage{scalerel,stackengine}
\stackMath
\newcommand\widecheck[1]{%
	\savestack{\tmpbox}{\stretchto{%
			\scaleto{%
				\scalerel*[\widthof{\ensuremath{#1}}]{\kern1.2pt\bigwedge\kern1.2pt}%
				{\rule[-\textheight/2]{1ex}{\textheight}}%WIDTH-LIMITED BIG WEDGE
			}{\textheight}% 
		}{0.5ex}}%
	\stackon[1pt]{#1}{\scalebox{-1}{\tmpbox}}%
}

\newtheorem{remark}{Remark}
\newtheorem{lemma}{Lemma}
\usepackage{amssymb}
\usepackage{amsmath}
\journal{Journal of Computational Physics}

\begin{document}

\begin{frontmatter}



\title{Learning Singularity-Encoded Green’s Functions with Application to Iterative Methods}

\tnotetext[t1]{Q. Sun's work is partially supported by National Natural Science Foundation of China under grant number 12201465 and Science and Technology Commission of Shanghai Municipality under grant number 23JC1400502. X. Xu's work is partially supported by National Natural Science Foundation of China under grant numbers 12071350 and 12331015.}

\author[inst1,inst2]{Qi Sun \corref{cor1}} %% Author name
\ead{qsun\_irl@tongji.edu.cn}
\cortext[cor1]{Corresponding author}

\author[inst1]{Shengyan Li} %% Author name
\ead{2410285@tongji.edu.cn}

\author[inst3,inst4]{Bowen Zheng} %% Author name
\ead{zhengbowen@lsec.cc.ac.cn}

\author[inst5]{Lili Ju} %% Author name
\ead{ju@math.sc.edu}

\author[inst1,inst2]{Xuejun Xu} %% Author name
\ead{xuxj@tongji.edu.cn}

\affiliation[inst1]{organization={School of Mathematical Sciences, Tongji University},%Department and Organization
            %addressline={}, 
            city={Shanghai},
            postcode={200092}, 
            %state={},
            country={China}}

\affiliation[inst2]{organization={Key Laboratory of Intelligent Computing and Applications (Ministry of Education),  Tongji University},%Department and Organization
            %addressline={}, 
            city={Shanghai},
            postcode={200092}, 
            country={China}}
            
\affiliation[inst3]{organization={LSEC, ICMSEC, Academy of Mathematics and Systems Science, Chinese Academy of Sciences},%Department and Organization
            %addressline={}, 
            city={Beijing},
            postcode={100190}, 
            country={China}}    
            
\affiliation[inst4]{organization={School of Mathematical Sciences, University of Chinese Academy of Sciences},%Department and Organization
            %addressline={}, 
            city={Beijing},
            postcode={100190}, 
            country={China}}                                                           

\affiliation[inst5]{organization={Department of Mathematics, University of  South Carolina},%Department and Organization
            %addressline={}, 
            city={Columbia},
            state={SC},
            postcode={29208}, 
            country={USA}}

%%%%%%%%%%%%%%%%%%%%%%%%%%%%%%%%%%%%%%%%%%%%%%%%%%%%%
%%%%%%%%%%%%%%%%%%%%%%%%%%%%%%%%%%%%%%%%%%%%%%%%%%%%%
%% Abstract
\begin{abstract}
Green's function provides an inherent connection between theoretical analysis and numerical methods for elliptic partial differential equations, and general absence of its closed-form expression  necessitates surrogate modeling to guide the design of effective solvers. Unfortunately, numerical computation of Green's function remains challenging due to its doubled dimensionality and intrinsic singularity. In this paper, we present a novel singularity-encoded learning approach to resolve these problems in an unsupervised fashion. Our method embeds the Green's function within a one-order higher-dimensional space by encoding its prior estimate as an augmented variable, followed by a neural network parametrization to manage the increased dimensionality. By projecting the trained neural network solution back onto the original domain, our deep surrogate model exploits its spectral bias to accelerate conventional iterative schemes, serving either as a preconditioner or as part of a hybrid solver. The effectiveness of our proposed method is empirically verified through numerical experiments  with two and four dimensional Green's functions, achieving satisfactory resolution of singularities and acceleration of iterative solvers.		
\end{abstract}
%%%%%%%%%%%%%%%%%%%%%%%%%%%%%%%%%%%%%%%%%%%%%%%%%%%%%
%%%%%%%%%%%%%%%%%%%%%%%%%%%%%%%%%%%%%%%%%%%%%%%%%%%%%            

%%%%%%%%%%%%%%%%%%%%%%%%%%%%%%%%%%%%%%%
% removed by Qi

%%%%%%%%%%%%%%%%%%%%%%%%%%%%%%%%%%%%%%%
% removed by Qi

%%%%%%%%%%%%%%%%%%%%%%%%%%%%%%%%%%%%%%%%%%%%%%%%%%%%%
%%%%%%%%%%%%%%%%%%%%%%%%%%%%%%%%%%%%%%%%%%%%%%%%%%%%%
%% Keywords
%%%%%%%%%%%%%%%%%%%%%%%%%%%%%%%%%%%%%%%%%%%%%%%%%%%%%
%%%%%%%%%%%%%%%%%%%%%%%%%%%%%%%%%%%%%%%%%%%%%%%%%%%%%

\end{frontmatter}

%%%%%%%%%%%%%%%%%%%%%%%%%%%%%%%%%%%%%%%%%%%%%%%%%%%%%
%%%%%%%%%%%%%%%%%%%%%%%%%%%%%%%%%%%%%%%%%%%%%%%%%%%%%
\section{Introduction}
	
	Machine learning has achieved remarkable success across a wide range of engineering applications and scientific disciplines, emerging as an innovative methodology toward modern scientific computing \cite{lecun2015deep,goodfellow2016deep,karniadakis2021physics,azizzadenesheli2024neural}. Owing to its universal approximation capability and mesh-free discretization \cite{lagaris1998artificial,raissi2019physics}, the neural network-based method has become increasingly popular for the numerical solution of partial differential equations \cite{evans2010partial,leveque2007finite}, especially in high dimensions \cite{han2018solving} where classical mesh-based methods exhibit limitations. Scientific machine learning approaches within this field can be broadly categorized into two paradigms: solution learning and operator learning.
	
	To be specific, we consider the benchmark elliptic boundary value problem of the form \cite{evans2010partial,gilbarg1977elliptic}
	\begin{equation}
		\begin{array}{rll}
			-\nabla \cdot \left(c(\bm{x}) \nabla u(\bm{x})\right) \!\!\!\!\! &= f(\bm{x}), \ &\ \ \ \text{in}\ \Omega,\\
			u(\bm{x}) \!\!\!\!\!  &= 0, \ &\ \ \ \text{on}\ \partial\Omega,
		\end{array}
		\label{BVP-ExactSolu}
	\end{equation}
	where $\Omega\subset\mathbb{R}^d$ $(d\in\mathbb{N}_+)$ is a bounded domain with a $C^{2,\alpha}$ boundary $\partial\Omega$ for some $\alpha \in (0,1)$, $c(\bm{x}) \in C^{1,\alpha}(\bar{\Omega})$ its coefficient function that is strictly positive and uniformly elliptic, and $f(\bm{x})\in L^2(\Omega)$ the source term. Throughout this paper, the vector variable $\bm{x}=(x_1,\cdots,x_d)$ is notationally simplified to a scalar variable $x\in\mathbb{R}$ when $d=1$. In solution learning methods such as the physics-informed neural networks \cite{raissi2019physics}, the deep Ritz method \cite{yu2018deep}, and the deep Galerkin method \cite{sirignano2018dgm}, the solution of problem \eqref{BVP-ExactSolu} is parametrized using a neural network and optimized via minimizing either the residual of governing equations or the energy functional. These approaches are designed in an instance-specific manner, requiring retraining for any changes in the forcing term or boundary conditions.
	
	On the other hand, operator learning aims to approximate mappings between function spaces \cite{lu2021learning,li2020fourier,azizzadenesheli2024neural}, such as the solution operator of problem \eqref{BVP-ExactSolu} that maps a forcing term to its corresponding solution, namely,
	\begin{equation}
		u(\bm{x}) = \int_{\Omega} G(\bm{x},\bm{y})f(\bm{y}) \,d\bm{y},     
		\label{BVP-ExactSolu-Representation-ExactGreen}	
	\end{equation}
	in terms of the $2d$-dimensional Green's function $G(\bm{x},\bm{y}): \Omega \times \Omega \to \mathbb{R}\cup\{ \infty \}$ \cite{evans2010partial,gilbarg1977elliptic}. Such an analytic expression has motivated the construction of neural network models \cite{boulle2022data,li2020fourier,gin2021deepgreen,melchers2024neural,kovachki2023neural} that directly learn the integral operator \eqref{BVP-ExactSolu-Representation-ExactGreen} in a data-driven manner, while alternative approaches \cite{lu2021learning,shukla2024comprehensive} build on the universal approximation theorems. In addition to serving as fast solvers, neural operators are being seamlessly integrated into conventional scientific computing workflows to enhance computational efficiency \cite{azizzadenesheli2024neural}. The mesh-free nature enables their deployment as effective preconditioning strategies for the discrete system of equation \eqref{BVP-ExactSolu} across different mesh sizes \cite{kopanivcakova2025deeponet}, while their inherent spectral bias facilitates hybridization with classical iterative methods \cite{zhang2024blending,cui2025hybrid,hu2025hybrid}.
	
	Although end-to-end neural operators have recently attracted considerable scientific interest, their reliance on extensive high-quality training datasets and limited physical interpretability remain challenging \cite{azizzadenesheli2024neural}. An emerging alternative is to learn the Green's function directly, either from discrete source-to-solution data pairs \cite{boulle2022data,gin2021deepgreen,lin2024green} or by exploiting its governing equations through physics-informed machine learning \cite{lin2021binet,teng2022learning,ji2023deep,hu2023solving,lin2024green}. Data-driven methods, as previously noted, often require extensive labelled datasets to attain robust generalization, making the unsupervised learning approach preferable in practice \cite{li2024physics}. To avoid the computation of the Dirac delta function, methods such as the Gaussian approximation \cite{teng2022learning,ji2023deep,hao2024multiscale}, the fundamental solution-based subtraction \cite{lin2021binet,hu2023solving} and the variational reformulation \cite{chen2023friedrichs} are developed. On the other hand, the low regularity of Green’s functions makes it necessary to establish specialized neural network architectures \cite{boulle2020rational,wimalawarne2023learning,hao2024multiscale,lin2024green} to ensure approximation accuracy. Specifically, the Green's function of our elliptic problem \eqref{BVP-ExactSolu} exhibits a piecewise smooth singularity when $d=1$, a logarithmic singularity when $d=2$, and an algebraic singularity when $d\geq 3$ \cite{gruter1982green,kim2019green,dolzmann1995estimates,mitrea2010regularity}, which is inherent to the kernel function that persists regardless of data-driven or physics-informed learning approaches.
	
	In this paper, we present a novel singularity-encoded learning framework to reconstruct the Green's function without supervised training data. Building upon the analytic estimates of its singular behavior \cite{stakgold2011green,kim2019green,mitrea2010regularity,gruter1982green}, the Green's function is embedded into a higher-dimensional space by incorporating an augmented variable, followed by the parametrization via neural networks to handle the increased dimensionality. It is noteworthy that the Dirac delta function is directly computed through integration along the singular region, rather than relying on Gaussian approximations or fundamental solutions. Besides, our explicit encoding of singularities enables a smooth neural network to reconstruct singular functions. The surrogate model for our Green's functions can then be obtained by projecting the trained network solution back onto the physical domain. By exploiting the spectral bias inherent in the training process, our singularity-encoded Green's function is utilized to accelerate classical iterative methods, either as a preconditioning matrix or a component in our hybrid solvers.
	
	The remainder of this paper is organized as follows. Section \ref{Section-RelatedWork} presents a comprehensive overview of existing methods for learning the Green's function via neural networks, together with the employment of neural operators in accelerating classical iterative methods. Our singularity-encoded learning method is presented in \autoref{section-SgEncd-Green-Function} and is illustrated through a one-dimensional benchmark problem for visualization. In \autoref{Section-Acltd-IteM}, the pre-trained model is adopted to accelerate classical iterative schemes, serving either as a preconditioner or part of our hybrid solver. Numerical validation through extra experiments are provided in \autoref{Sec-Experiments}, followed by conclusions in \autoref{Sec-Conclusion}.

%%%%%%%%%%%%%%%%%%%%%%%%%%%%%%%%%%%%%%%%%%%%%%%%%%%%%
%%%%%%%%%%%%%%%%%%%%%%%%%%%%%%%%%%%%%%%%%%%%%%%%%%%%%

%%%%%%%%%%%%%%%%%%%%%%%%%%%%%%%%%%%%%%%%%%%%%%%%%%%%%
%%%%%%%%%%%%%%%%%%%%%%%%%%%%%%%%%%%%%%%%%%%%%%%%%%%%%
\section{Related Work}\label{Section-RelatedWork}
	
%%%%%%%%%%%%%%%%%%%%%%%%%%%%%%%%%%%%%%%%%%%%%%%%%%
	
\subsection{Learning Green's Functions with Neural Networks}\label{Sec-RelatedWork-GreenFunc}
	
	For any impulse source located at $\bm{y}\in \Omega$, the Green's function associated with problem \eqref{BVP-ExactSolu} satisfies
	\begin{equation}
		\begin{array}{rll}
			-\nabla \cdot \left( c(\bm{x}) \nabla G(\bm{x},\bm{y})\right) \!\!\!\!\!\! &= \delta(\bm{x}-\bm{y}), \ &\ \ \ \text{in}\ \Omega,\\
			G(\bm{x},\bm{y}) \!\!\!\!\!\! &= 0, \ &\ \ \ \text{on}\ \partial\Omega,
		\end{array}
		\label{BVP-ExactGreen}
	\end{equation}
	where $\delta(\bm{x})$ denotes the Dirac delta function. Unfortunately, the analytic expression of Green's function is often difficult to obtain, especially for problems with complex geometries or variable coefficients \cite{evans2010partial,stakgold2011green}. Alternatively, numerical approximation of Green's function becomes requisite, but conventional mesh-based methods, such as the finite difference or finite element methods \cite{brenner2008mathematical,leveque2007finite}, are confronted with some long-standing challenges \cite{lin2021binet} as follows:
	\begin{itemize}
	\item[--] the doubled dimensionality of Green's function incurs the curse of dimensionality;
	\item[--] the complex domain makes the mesh generation for Green’s functions impractical.
	\end{itemize}
	Neural network-based learning approaches \cite{goodfellow2016deep}, being mesh-free and well-suited for high-dimensional problems, have attracted significant attention in recent years, inspiring innovations in both network architectures \cite{ji2023deep,hao2024multiscale,lin2024green,boulle2020rational} and training algorithms \cite{lin2021binet,boulle2022data,wimalawarne2023learning,gin2021deepgreen,teng2022learning,zhang2021mod}. Unfortunately, a straightforward implementation via the physics-informed machine learning \cite{lagaris1998artificial,raissi2019physics,karniadakis2021physics} is complicated by several issues, among which a limitation arises from
the distributional nature of the Dirac delta function prohibits pointwise evaluation.
To bypass this issue, existing techniques can be roughly classified into three categories: (1) replacing the Dirac delta function with a Gaussian density function \cite{teng2022learning,hao2024multiscale}, despite introducing inevitable modeling errors; (2) employing the fundamental solution to decompose the Green's function into regular and singular components \cite{evans2010partial,lin2021binet,hu2023solving}, but the explicit expression of fundamental solutions is generally unavailable \cite{stakgold2011green}; (3) reformulating the problem into a variational form, with the weak solution being learned from a minimax optimization problem \cite{chen2023friedrichs}. 

	In addition to the aforementioned limitations, another fundamental challenge emerges from the fact that standard neural networks usually lack the capacity to represent/capture singular functions, regardless of whether the training strategy is supervised \cite{boulle2022data,gin2021deepgreen,lin2024green} or unsupervised \cite{teng2022learning,zhang2021mod,ji2023deep,wimalawarne2023learning}, or whether the training loss function is constructed in a strong or variational manner \cite{teng2022learning,ji2023deep,chen2023friedrichs}. To handle the singular behavior, rational activation functions \cite{boulle2020rational,boulle2022data} are proposed to improve the approximation capacity of neural networks, while the low-rank \cite{wimalawarne2023learning}, multiscale \cite{hao2024multiscale}, and hierarchical \cite{lin2024green} structures inherent in the Green's function can also be incorporated into the modification of neural network architectures.

	Unlike prior learning approaches, our proposed method first embeds the Green's functions within a higher-dimensional space by introducing an augmented variable, after which a smooth neural network is adopted for the parametric approximation. Our strategy also facilitates the effective treatment of the Dirac delta function without relying on  Gaussian approximation or fundamental solution. %(see \autoref{section-SgEncd-Green-Function} for more details).

%%%%%%%%%%%%%%%%%%%%%%%%%%%%%%%%%%%%%%%%%%%%%%%%%%	
\subsection{Iterative Methods with Pre-Trained Neural Operators}

	Traditional numerical methods for solving problem \eqref{BVP-ExactSolu} typically involves discretization through finite difference or finite element methods \cite{leveque2007finite,brenner2008mathematical}, yielding a linear system of equations 
	\begin{equation}\label{std}
		AU=F
	\end{equation}
	where $A$ is the coefficient matrix, $U$ the vector of unknowns, and $F$ the right-hand side. Iterative solution techniques \cite{varga1962iterative,golub2013matrix} have been widely employed to solve such systems, but often exhibit deteriorating performance under mesh refinement as the coefficient matrix becomes increasingly ill-conditioned. Preconditioning is therefore essential, that is, instead of \eqref{std} one  solves 
	\begin{equation}\label{stdpre}
		BAU=BF
	\end{equation}	
	where the preconditioner $B$ should approximate $A^{-1}$ effectively so as to both cluster the eigenvalues of $BA$ and reduce its condition number \cite{golub2013matrix}. Matrix factorization-based preconditioners, such as the damped Jacobi, Gauss-Seidel, and successive over-relaxation methods, are well-established but exhibit slow convergence rates in practice \cite{trefethen2022numerical,saad2003iterative}. Multilevel approaches, such as the geometric multigrid \cite{trottenberg2001multigrid,golub2013matrix} and domain decomposition \cite{toselli2004domain} methods, achieve improved performance by addressing errors across different scales of spatial discretization.
	
	In addition to classical preconditioning strategies, pre-trained neural operators \cite{lu2021learning,li2020fourier} have recently emerged as a new class of preconditioners by exploiting their inherent spectral bias \cite{xu2025understanding} to improve the spectral properties of the preconditioned system \cite{kopanivcakova2025deeponet,cui2025hybrid,li2025neural}. More precisely, neural operators, such as the deep operator network \cite{lu2021learning} and Fourier neural operator \cite{li2020fourier}, are neural networks that are commonly trained in a supervised learning manner to approximate the inverse of the differential operator for equations \eqref{BVP-ExactSolu}. Owing to their mesh-independent nature, these pre-trained models permit seamless implementation at the discrete level, particularly when deployed as the preconditioner \cite{kopanivcakova2025deeponet,cui2025hybrid,hu2025hybrid,li2025neural} or the coarse-grid solver \cite{klawonn2024learning,kopanivcakova2025deeponet}. Despite their limited overall accuracy, the spectral bias of neural operators toward low-frequency modes makes them well-suited for multilevel preconditioners.

	Moreover, the spectral bias inherent in pre-trained neural operators also allows their integration with classical iterative methods, facilitating the development of hybrid schemes with enhanced convergence properties \cite{zhang2024blending}. To be specific, given an initial guess $U^{[1]}$, the classical iterative solver for  our discrete system \eqref{stdpre} reads
	\begin{equation}
		U^{[k+1]} = U^{[k]} + B ( f - AU^{[k]} )
	\end{equation}	
	where $B$ represents a classic  preconditioner that is effective for eliminating high-frequency error components \cite{varga1962iterative,golub2013matrix} but exhibits limited efficacy for low-frequency errors (typically referred to as the smoothing effects). The  other approach employs the pre-trained neural operator $\mathcal{B}$ to construct the following iterative update
	\begin{equation}\label{stdneu}
		U^{[k+1]} = U^{[k]} + \mathcal{B} ( f - AU^{[k]} )
	\end{equation}		
	which preferentially attenuates low-frequency errors due to the spectral bias of neural operators \cite{zhang2024blending,hu2025hybrid,cui2025hybrid}, but at the cost of introducing additional high-frequency errors owing to their limited accuracy. As a direct consequence, a natural extension is to combine both iterative schemes into a hybrid method \cite{zhang2024blending,hu2025hybrid} that improves convergence across the entire spectrum of error modes \cite{cui2025hybrid}. Moreover, the hybrid method has been successfully integrated into the multigrid framework as pre- and post-smoothing procedures \cite{wang2024label}, proving valuable for hierarchical solvers.
	
	While neural operators have demonstrated notable potential in accelerating classical iterative methods, either as preconditioners or within hybrid iterative schemes, the majority of existing implementations \cite{zhang2024blending,hu2025hybrid,cui2025hybrid} rely on data-driven techniques that raise fundamental concerns regarding generalization and interpretability \cite{azizzadenesheli2024neural}. Besides, the challenge of acquiring high-quality and diverse training data remains a critical bottleneck in the development of robust neural operators, while the lack of transparency hinders their reliability in practical applications.
	
	Unlike widely adopted data-driven models, our pre-trained neural operator is constructed by directly learning the Green’s functions \eqref{BVP-ExactGreen} within an unsupervised learning regime, thereby removing dependence on labelled data and enhancing interpretability. Accordingly, our preconditioner is explicitly expressed in a matrix representation, in contrast to the end-to-end mapping realized by the deep operator network  or Fourier neural operator. Though distinct training strategies are employed, the spectral bias inherent in trained models remains unchanged, therefore enabling the acceleration of classical iterative methods through our learned Green’s functions.
		
%%%%%%%%%%%%%%%%%%%%%%%%%%%%%%%%%%%%%%%%%%%%%%%%%%%%
%%%%%%%%%%%%%%%%%%%%%%%%%%%%%%%%%%%%%%%%%%%%%%%%%%%%

%%%%%%%%%%%%%%%%%%%%%%%%%%%%%%%%%%%%%%%%%%%%%%%%%%%%%
%%%%%%%%%%%%%%%%%%%%%%%%%%%%%%%%%%%%%%%%%%%%%%%%%%%%%
\section{Singularity-Encoded Green's Functions}\label{section-SgEncd-Green-Function}
	
	In this section, the learning algorithm for our singularity-encoded Green's function is presented, accompanied by an empirical analysis that reveals the spectral bias of trained neural network models. We begin by introducing an augmented variable to embed the governing equations into a one-order higher-dimensional space, demonstrating its effectiveness in resolving singularities via a benchmark problem. Next, the algorithmic framework is detailed, followed by a spectral analysis of our trained models that motivates the design of hybrid iterative solvers.

%%%%%%%%%%%%%%%%%%%%%%%%%%%%%%%%%%%%%%%%%%%%%%%%%%%%%
\subsection{Green's Functions Embedded in the One-order Higher Dimensional Space}\label{Sec-Green-Higher-Dimension}
	
	The primary difficulty in learning Green's function stems from its inherent singular behavior near $\bm{x} = \bm{y}$, i.e.,
	\begin{equation}
		c(x) \llbracket \nabla G(x,y) \rrbracket = -1 \ \ \ \textnormal{for}\ d= 1,
		\ \ \ \textnormal{and} \ \ \ 
		G(\bm{x}, \bm{y}) \sim \left\{
		\begin{array}{cl}
			\displaystyle \ln{\| \bm{x} - \bm{y} \|}, \ &\ \ \ \textnormal{for}\ d=2,\\
			\displaystyle \| \bm{x} - \bm{y} \|^{2-d}, \ &\ \ \ \textnormal{for}\ d\geq 3,
		\end{array}\right.
		\label{Singularity-ExactGreen}
	\end{equation}
	which challenges both the approximation capacity of neural networks \cite{boulle2020rational} and the practical training procedure via auto-differentiation \cite{baydin2018automatic}. In this work, an augmented variable $\varphi( \bm{x}, \bm{y} )$ is incorporated as an additional input feature into our representation of Green's function (referred to as the ``singularity-encoded Green's function''):
	\begin{equation}
		G( \bm{x}, \bm{y} ) = \widehat{G}( \bm{x}, \bm{y},  \varphi( \bm{x}, \bm{y} )) 
		\ \ \ \textnormal{with} \ \ \
		\varphi( \bm{x}, \bm{y} ) = \left\{
		\begin{array}{cl}
			\displaystyle |x-y|, \ &\ \ \ \mathrm{for}\ d=1,\\
			\displaystyle \ln {\| \bm{x} - \bm{y} \|}, \ &\ \ \ \textnormal{for}\ d=2,\\
			\displaystyle \| \bm{x} - \bm{y} \|^{2-d}, \ &\ \ \ \textnormal{for}\ d\geq 3,
		\end{array}\right.
		\label{Singularity-SgEncdGreen}
	\end{equation}
	thereby facilitating the resolution of singularities during the training of neural networks.	Throughout our paper, $\llbracket \cdot \rrbracket$ denotes the difference of quantity across the discontinuity interface $\Gamma = \{ x \,|\, x=y\}$ \cite{stakgold2011green}.
	
	Before introducing the neural network parametrization of our Green's function, a direct calculation implies
	\begingroup
	\renewcommand*{\arraystretch}{1.5}
	\begin{equation*}
		\begin{array}{cl}
			\displaystyle \nabla \cdot \big( c(\bm{x}) \nabla G( \bm{x}, \bm{y} ) \big) = & \!\!\!\! \nabla c(\bm{x}) \cdot \big( \nabla_{\bm{x}} \widehat{G}( \bm{x}, \bm{y},  \varphi( \bm{x}, \bm{y} ))  +  \nabla \varphi( \bm{x}, \bm{y} ) \partial_z \widehat{G}( \bm{x}, \bm{y},  \varphi( \bm{x}, \bm{y} )) \big) \\
			& \!\!\!\! +\, c(\bm{x}) \big( \Delta \varphi( \bm{x}, \bm{y} )  \partial_z \widehat{G}( \bm{x}, \bm{y},  \varphi( \bm{x}, \bm{y} ))  
			 + \, \Delta_{\bm{x}} \widehat{G}( \bm{x}, \bm{y},  \varphi( \bm{x}, \bm{y} )) \\
			 & \!\!\!\! +\, 2 \nabla \varphi( \bm{x}, \bm{y} ) \cdot \nabla_{\bm{x}} ( \partial_z \widehat{G}( \bm{x}, \bm{y},  \varphi( \bm{x}, \bm{y} )) )  + \|\nabla \varphi( \bm{x}, \bm{y} ) \|^2 \partial_{zz} \widehat{G}( \bm{x}, \bm{y},  \varphi( \bm{x}, \bm{y} )) \big) 
		\end{array}
	\end{equation*}
	\endgroup
	where $\nabla_{\bm{x}}$, $\Delta_{\bm{x}}$, and $\partial_z$ represent partial derivatives with respect to the $\bm{x}$- or $z$-variable of function $\widehat{G}( \bm{x}, \bm{y},  z)$, 
	\begin{equation*}
		\nabla  \varphi( \bm{x}, \bm{y} ) = \left\{
		\begin{array}{cl}
			\displaystyle \textnormal{sgn}(x-y), \ &\ \ \ \mathrm{for}\ d=1,\\
			\displaystyle \lVert \bm{x}-\bm{y} \rVert^{-2}(\bm{x}-\bm{y}), \ &\ \ \ \textnormal{for}\ d=2,\\
			\displaystyle (2-d) \lVert \bm{x} - \bm{y} \rVert ^{-d} (\bm{x}-\bm{y}), \ &\ \ \ \textnormal{for}\ d\geq 3,
		\end{array}\right.
		\ \ \ \textnormal{and} \ \ \ 
		\Delta  \varphi( \bm{x}, \bm{y} ) = \left\{
		\begin{array}{cl}
			\displaystyle 0, \ &\ \ \ \mathrm{for}\ d=1,\\
			\displaystyle 0, \ &\ \ \ \textnormal{for}\ d=2,\\
			\displaystyle 0, \ &\ \ \ \textnormal{for}\ d\geq 3.
		\end{array}\right.
	\end{equation*}
	
	To account for the $\delta$-function in equation \eqref{BVP-ExactGreen}, we integrate over some ball $B_\epsilon(\bm{y}) = \{ \bm{x} \, | \, \| \bm{x} -\bm{y} \| < \epsilon\}$ to obtain the so-called normalization condition for each $\bm{x}\in \Gamma = \{ \bm{x} \,|\, \bm{x}=\bm{y} \in \Omega \}$, namely,
	\begin{equation*}
		- \int_{B_\epsilon(\bm{y})} \nabla \cdot ( c(\bm{x}) \nabla G(\bm{x}, \bm{y}) ) \, d \bm{x} = - \int_{\partial B_\epsilon(\bm{y})} c(\bm{x}) \nabla G(\bm{x}, \bm{y}) \cdot \bm{n}_{\bm{x}} \, dS(\bm{x}) = 1,
	\end{equation*}
	where $ \bm{n}(\bm{x}) = \frac{\bm{x} - \bm{y}}{\lVert \bm{x} - \bm{y} \rVert}$ is the unit outer normal vector along the boundary $\partial B_\epsilon(\textbf{y})$. Accordingly, it can be recast as
	\begin{equation}
		- \int_{\partial B_\epsilon(\bm{y})} c(\bm{x}) \Big( \nabla_{\bm{x}} \widehat{G}( \bm{x}, \bm{y},  \varphi( \bm{x}, \bm{y} )) + \nabla \varphi(\bm{x},\bm{y}) \partial_z\widehat{G}(\bm{x}, \bm{y},  \varphi( \bm{x}, \bm{y} )) \Big) \cdot \bm{n}_{\bm{x}} \, dS(\bm{x}) = 1
		\label{normalization-condition-embedded}
	\end{equation}
	in terms of our singularity-encoded Green's function. It is noteworthy that for $d=1$, an equivalent expression of \eqref{normalization-condition-embedded} can be obtained by sending $\epsilon\to 0$, namely, 	
	\begin{equation}
		- 2 c(x) \partial_z \widehat{G}(x,y,\varphi(x,y)) = - 2 c(x) \partial_z \widehat{G}(x,y,|x-y|)  = 1, \ \ \ \textnormal{for}\ x\in\Gamma=\{ x \,|\, x = y \},
		\label{normalization-condition-embedded-1D}
	\end{equation}
	which allows the derivative jump of Green’s function \eqref{Singularity-ExactGreen} to be accurately represented by our singularity-encoded Green's function. In other words, the piecewise-continuous Green's function $G(x,y)$ is mapped onto two disjoint hyperplanes within a 3-dimensional space (see \autoref{Experiment-1D-ToyModel} for instance), which admits a continuous extension to the entire space $G(x,y,z)$ as ensured by the Whitney extension theorem \cite{tseng2023cusp,hu2025discontinuity}.
	
	To sum up, the problem satisfied by our singularity-encoded Green's function $\widehat{G} (\bm{x}, \bm{y}, \varphi( \bm{x}, \bm{y}) ) $ is defined as
	\begin{equation}
		\left\{
		\begin{array}{cl}
			- \nabla c \cdot \Big( \nabla_{\bm{x}} \widehat{G} + \nabla \varphi \partial_z \widehat{G} \Big) - c \Big( \Delta_{\bm{x}} \widehat{G} + 2 \nabla \varphi \cdot \nabla_{\bm{x}} ( \partial_z \widehat{G} )  + \|\nabla \varphi \|^2 \partial_{zz} \widehat{G} \Big) = 0, & \mathrm{for}\ \bm{x} \in \Omega\setminus\Gamma,\\
			\displaystyle - \int_{\partial B_\epsilon(\bm{y})} c \Big( \nabla_{\bm{x}} \widehat{G} + \nabla \varphi \partial_z\widehat{G} \Big) \cdot \bm{n}_{\bm{x}} \, dS(\bm{x}) = 1, & \textnormal{for}\ \bm{x}\in \Gamma= \{ \bm{x} \,|\, \bm{x}=\bm{y} \},\\
			\widehat{G} = 0, & \mathrm{for}\ \bm{x} \in \partial\Omega,
		\end{array}\right.
		\label{BVP-ExactGreen-Encoded}
	\end{equation}
	for any fixed $\bm{y}\in\Omega$ and $\epsilon >0$, which is embedded within a one-order higher-dimensional space compared to its original form \eqref{BVP-ExactGreen}. Here, and in what follows, we employ the general form \eqref{normalization-condition-embedded} for clarity, reserving the specialized case \eqref{normalization-condition-embedded-1D} whenever necessary. 
	
	In contrast to the derivative jump in \eqref{Singularity-ExactGreen}, the singularity inherent to the Green’s function becomes more severe in higher dimensions, posing challenges for conventional numerical methods due to the divergence of derivatives as $\bm{x}\to\bm{y}$. Physics-informed neural networks \cite{lagaris1998artificial,raissi2019physics,karniadakis2021physics}, which employ the automatic differentiation to calculate diverging derivatives near $\bm{x}=\bm{y}$, are no exception. Thanks to the inclusion of domain knowledge via augmenting the input variables \eqref{Singularity-SgEncdGreen}, our computation of singular first-order derivatives in \eqref{normalization-condition-embedded} is realized through
	\begin{equation}
		\nabla G(\bm{x},\bm{y}) = \nabla \widehat{G}(\bm{x},\bm{y}, \varphi(\bm{x},\bm{y})) = \nabla_{\bm{x}} \widehat{G}(\bm{x},\bm{y},\varphi(\bm{x},\bm{y})) + \nabla \varphi(\bm{x},\bm{y}) \partial_z \widehat{G}(\bm{x},\bm{y},\varphi(\bm{x},\bm{y})),
		\label{SgEncd-Green-1st-Derivatives}
	\end{equation} 
	in which $\nabla \varphi(\bm{x},\bm{y})$ grows at least at the same rate as $\nabla G(\bm{x},\bm{y})$. As a result, satisfactory accuracy could be attained for ensuring the constraint \eqref{normalization-condition-embedded} with derivatives $\nabla_{\bm{x}} \widehat{G}(\bm{x},\bm{y},\varphi(\bm{x},\bm{y}))$ and $\partial_z \widehat{G}(\bm{x},\bm{y},\varphi(\bm{x},\bm{y}))$ remaining bounded.
		
	\begin{remark}
	In fact, constructing the augmented variable for $d=2, 3$ requires only a rough upper-bound estimation of singularities, which is sufficient to bound partial derivatives of our singularity-encoded Green's functions \eqref{SgEncd-Green-1st-Derivatives}. For instance, in $d=2$ cases where the exact Green's function exhibits a logarithmic singularity \eqref{Singularity-ExactGreen}, our method remains effective when using $\varphi( \bm{x}, \bm{y} ) = \lVert \bm{x} - \bm{y} \rVert ^{-0.2}$ as numerically reported in \ref{Appendix-1}.
	\end{remark}
	
%%%%%%%%%%%%%%%%%%%%%%%%%%%%%%%%%%%%%%%%%%%%%%%%%%%%%
\subsubsection{Example: Poisson's Equation in One Dimension}\label{Section-1D-ToyModel}
		
	Before detailing the learning algorithms,  we  present a benchmark problem to visualize the proposed singularity-encoded Green's function \eqref{BVP-ExactGreen-Encoded} within a one-order higher-dimensional space. To be specific, we consider the following one-dimensional Poisson problem:
	\begin{equation}
		- u''(x) = f(x), \ \ \ \textnormal{for}\ x\in\Omega = (0,1),
		\label{BVP-1D-ToyModel}
	\end{equation}
	with homogeneous Dirichlet boundary conditions $u(0) =  u(1) = 0$ as the illustrative example. Then the boundary value problem governing its exact Green's function $G(x,y)$ is characterized by
	\begin{equation}
		\begin{array}{cl}
			-\partial_{xx} G(x,y) = 0, \ \ \ &\ \textnormal{for}\ x\in\Omega\setminus\Gamma, \\
			\llbracket G(x,y) \rrbracket = 0,\ \ \ \llbracket \partial_x G(x,y) \rrbracket = -1, \ \ \ &\ \textnormal{for}\ x\in\Gamma = \{ x \,|\, x=y\},\\
			G(x,y) = 0, \ \ \ &\ \textnormal{for}\ x\in\partial\Omega, \\
		\end{array}
		\label{BVP-1D-ToyModel-ExactGreen}
	\end{equation}
	for any fixed $y\in(0,1)$, while that of singularity-encoded Green's function $\widehat{G}(x,y,\varphi(x,y))$ takes on the form
	\begin{equation}
		\begin{array}{cl}
			-(\partial_{xx} + \partial_{zz}) \widehat{G}(x,y,\varphi(x,y)) - 2 \textnormal{sgn}(x-y) \partial_{x z}\widehat{G}(x,y,\varphi(x,y)) = 0, \ &\ \ \ \mathrm{for}\ x \in \Omega\setminus\Gamma,\\
			\displaystyle -2\partial_{z}\widehat{G}(x,y,\varphi(x,y)) = 1, \ &\ \ \ \textnormal{for}\ x\in \Gamma=\{ x \,|\, x = y \},\\
			\widehat{G} (x, y, \varphi( x, y) ) = 0, \ &\ \ \ \mathrm{for}\ x \in \partial\Omega.
		\end{array}
		\label{BVP-1D-ToyModel-SgEncdGreen}
	\end{equation}
	Notably, the exact Green's function associated with problem \eqref{BVP-1D-ToyModel} is given by
	\begin{equation*}
		G(x,y) = \left\{
		\begin{array}{cl}
			x(1-y),\ &\ \ \ \textnormal{for}\ 0\leq x\leq y,\\
			y(1-x),\ &\ \ \ \textnormal{for}\ y \leq x \leq 1.
		\end{array}\right.
	\end{equation*}
	while the constraint $0 = \llbracket G(x,y) \rrbracket = \llbracket \widehat{G}(x,y,\varphi(x,y)) \rrbracket$ for $x\in\Gamma$ is naturally enforced using $\varphi(x,y) = |x-y|$.	
		
	% figure
	%%--------------------------------------%
	\begin{figure}[!t]
		\centering
		\begin{subfigure}[b]{0.31\textwidth}
			\centering
			\includegraphics[width=\textwidth]{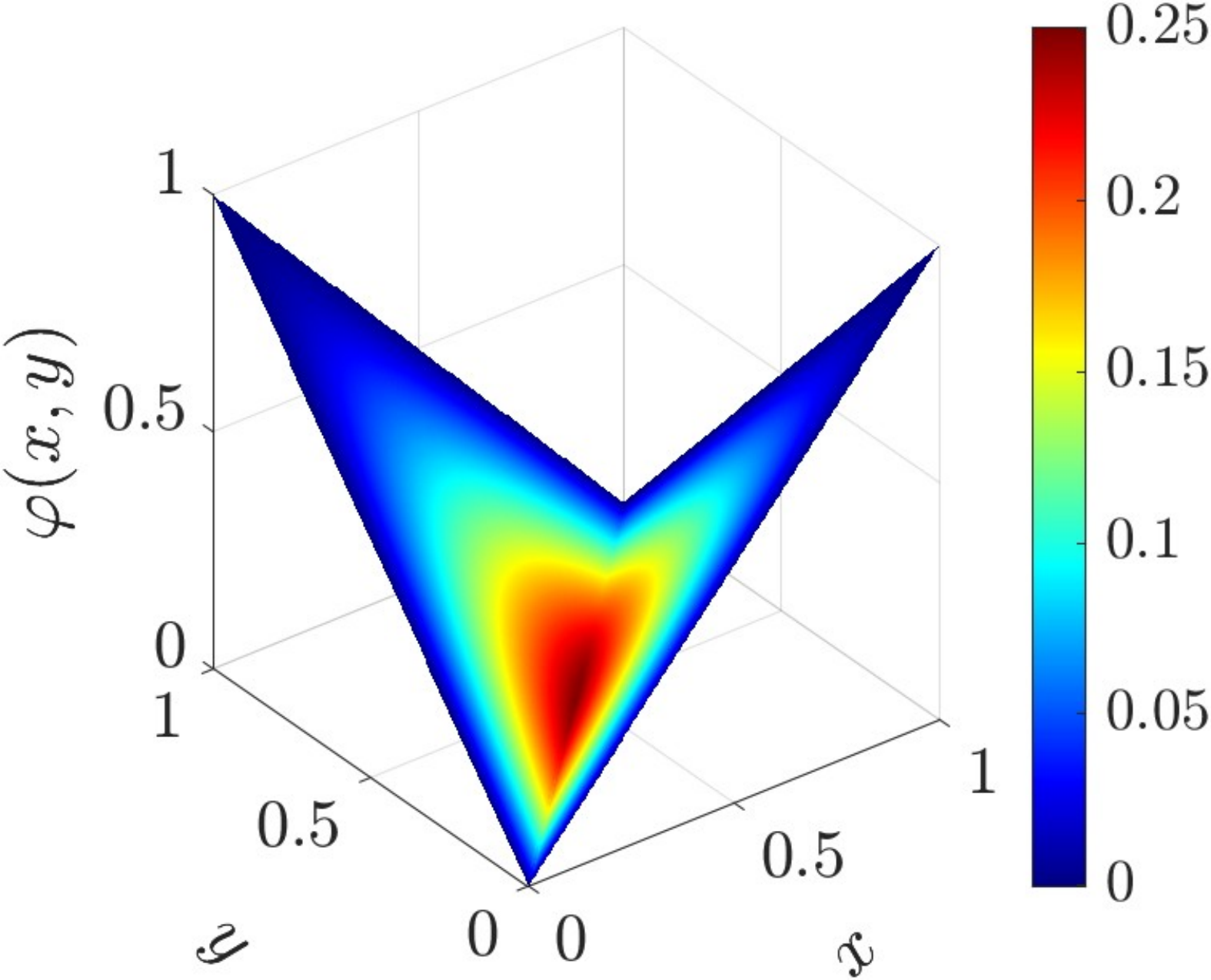}
			\caption{$\widehat{G}(x,y,\varphi(x,y))$}
			\label{fig-SgEncdGreen-1D-ToyModel-3D}
		\end{subfigure}
		\hspace{0.25cm}
		\begin{subfigure}{0.303\textwidth}
			\centering
			\includegraphics[width=\textwidth]{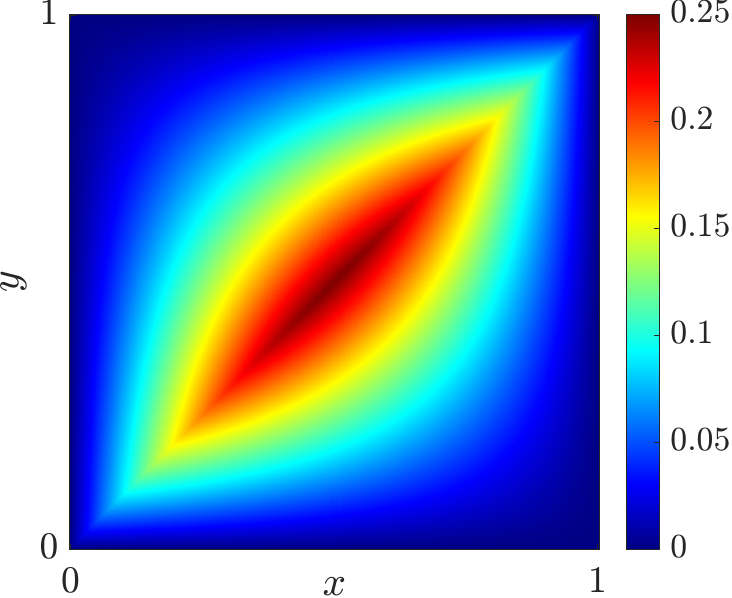}
			\caption{$\widecheck{G}(x,y)$}
			\label{fig-SgEncdGreen-1D-ToyModel-2D}
		\end{subfigure}
		\hspace{0.25cm}
		\begin{subfigure}{0.312\textwidth}
			\centering
			\includegraphics[width=\textwidth]{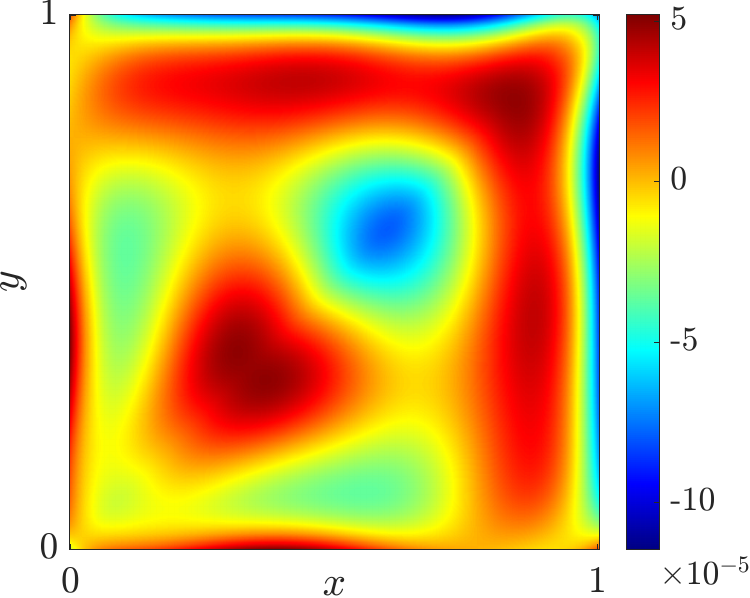}
			\caption{$G(x,y) - \widecheck{G}(x,y)$}
			\label{fig-SgEncdGreen-1D-ToyModel-2D-PtErr}
		\end{subfigure}
		
		\vspace{0.3cm}
		
		\begin{subfigure}{0.29\textwidth}
			\centering
			\includegraphics[width=\textwidth]{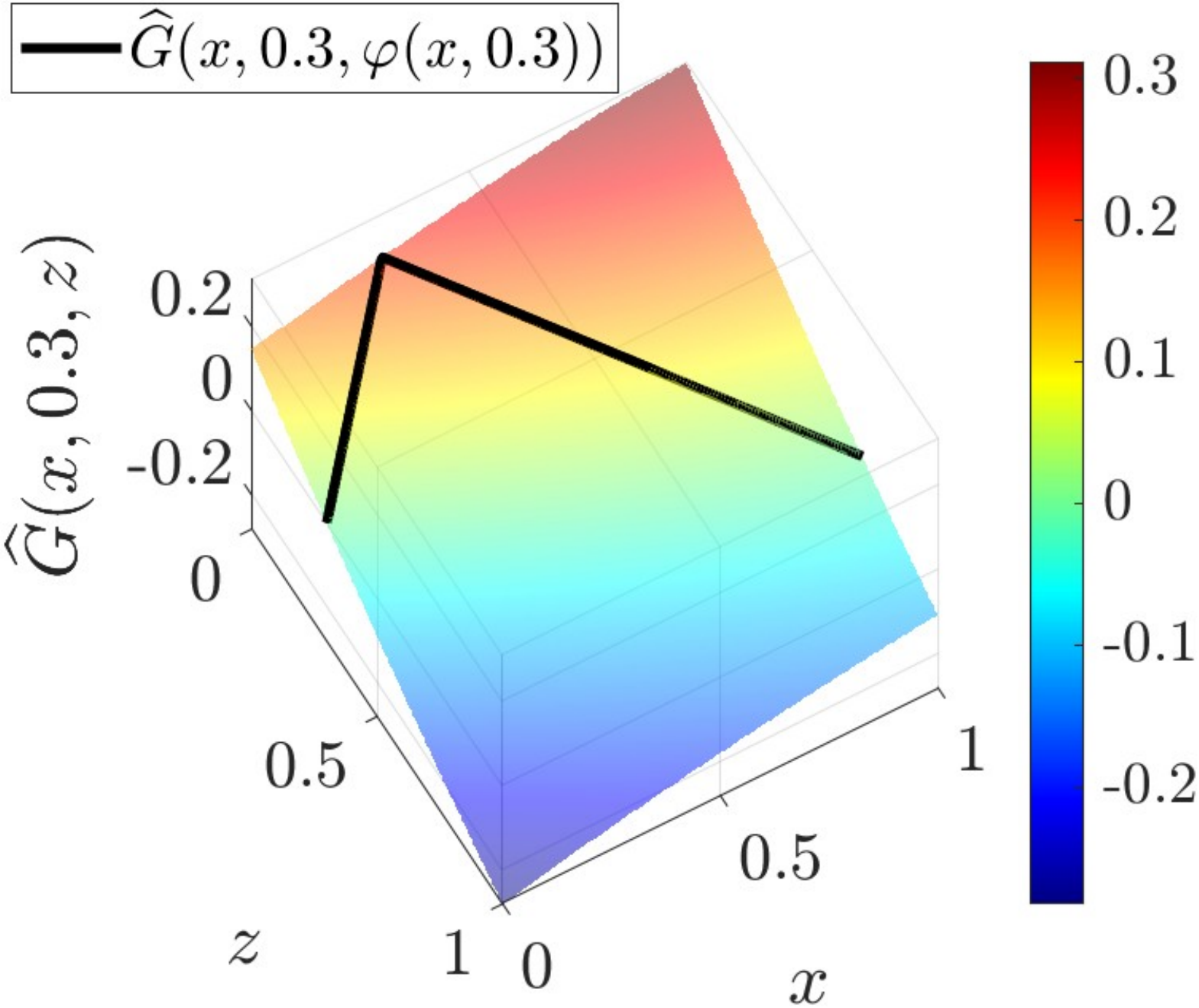}
			\caption{$\widehat{G}(x, 0.3, z)$}
			\label{fig-SgEncdGreen-1D-ToyModel-3D-Slice-x1}
		\end{subfigure}
		\hspace{0.25cm}
		\begin{subfigure}{0.30\textwidth}
			\centering
			\includegraphics[width=\textwidth]{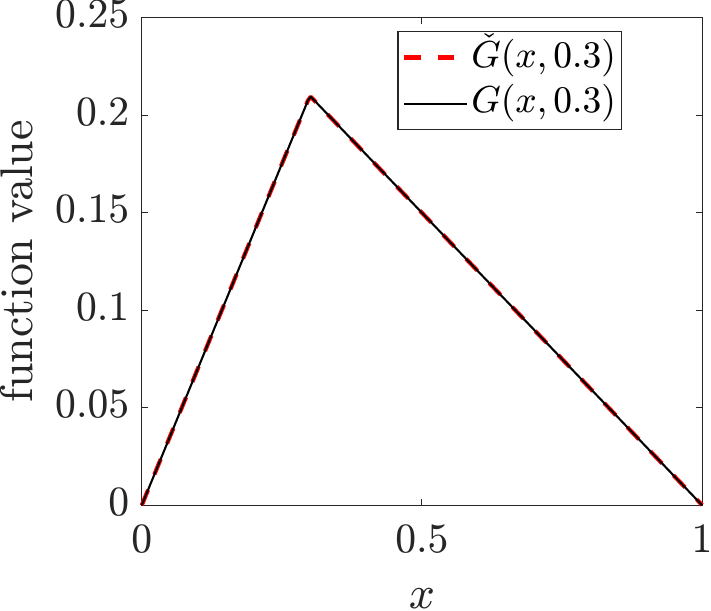}
			\caption{$\widecheck{G}(x,0.3)$ and $G(x, 0.3)$}
			\label{fig-SgEncdGreen-1D-ToyModel-2D-Slice-x1}
		\end{subfigure}
		\hspace{0.25cm}
		\begin{subfigure}{0.302\textwidth}
			\centering
			\includegraphics[width=\textwidth]{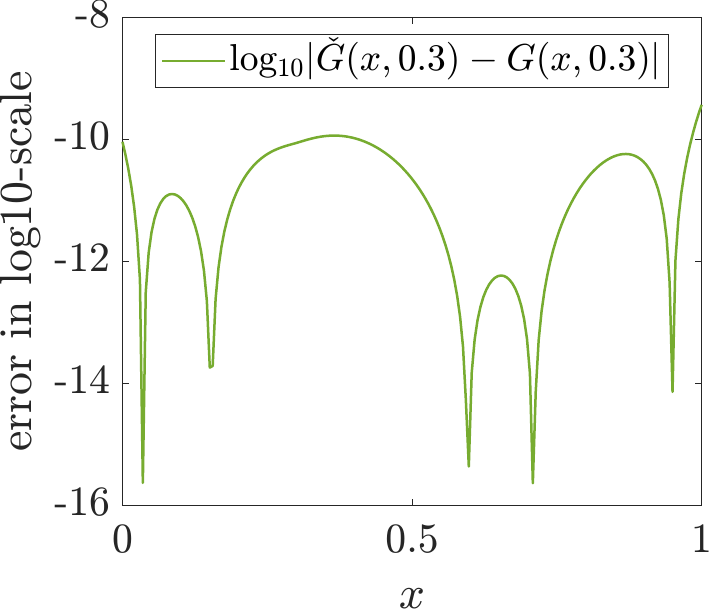}
			\caption{$\textnormal{log}_{10}|\widecheck{G}(x, 0.3) - G(x, 0.3)|$}
			\label{fig-SgEncdGreen-1D-ToyModel-2D-Slice-x1-PtErr}	
		\end{subfigure}
		
		\caption{Numerical visualizations of our singularity-encoded Green's function for the illustrative example \eqref{BVP-1D-ToyModel}.}
		\label{Experiment-1D-ToyModel}
		\vspace{-0.2cm}
	\end{figure}
	%%--------------------------------------%	
			
	Next, the unknown solution to \eqref{BVP-1D-ToyModel-SgEncdGreen} is parametrized using a fully-connected neural network \cite{hornik1989multilayer}, followed by a standard training process under the pointwise residual-minimization framework \cite{lagaris2000neural,raissi2019physics}. The trained model (still denoted as $\widehat{G} (x, y, \varphi( x, y) )$ for simplicity) and its projection (denoted as $\widecheck{G} (x, y)$) back onto the lower-dimensional plane are depicted in \autoref{fig-SgEncdGreen-1D-ToyModel-3D} and \autoref{fig-SgEncdGreen-1D-ToyModel-2D}, while \autoref{fig-SgEncdGreen-1D-ToyModel-2D-PtErr} shows its pointwise error profile compared to the true Green's function. Though the exact Green's function $G(x, 0.3)$ exhibits a jump in its derivative as depicted in \autoref{fig-SgEncdGreen-1D-ToyModel-2D-Slice-x1}, its numerical counterpart $\widecheck{G}(x,0.3)$ could be embedded into a smooth surface $\widehat{G}(x, 0.3, z)$ within a one-order higher-dimensional space (see \autoref{fig-SgEncdGreen-1D-ToyModel-3D-Slice-x1}), thereby allowing both our trained model and numerical solution to achieve satisfactory accuracy as demonstrated in \autoref{fig-SgEncdGreen-1D-ToyModel-2D-Slice-x1-PtErr} and \autoref{fig-SgEncdGreen-1D-ToyModel-2D-PtErr}.

%%%%%%%%%%%%%%%%%%%%%%%%%%%%%%%%%%%%%%%%%%%%%%%%%%%%%
\subsection{Singularity-Encoded Neural Networks}
	
	Thanks to the enhanced regularity achieved via embedding within a one-order higher-dimensional space, a fully-connected neural network \cite{goodfellow2016deep} with smooth activation functions, e.g., $\sigma(x) = \textnormal{tanh}(x)$, is deployed to parametrize our singularity-encoded Green's function (not relabelled for notational simplicity), namely,
	\begin{equation}
		\widehat{G}( \bm{x}, \bm{y}, z;\theta) = (L_{D+1} \circ \sigma \circ L_D \circ \cdots \circ \sigma \circ L_1 \circ L_0)( \bm{x}, \bm{y}, z)
		\label{SgEncdGreen-NN-Model}
	\end{equation}
	where $D\in\mathbb{N}_+$ denotes the number of hidden layers, $L_\ell$ the $\ell$-th linear layer that applies an affine linear transformation to the incoming data, and $\theta$ the collection of all trainable parameters. 
	
	Next, our embedded problem \eqref{BVP-ExactGreen-Encoded} is solved under an unsupervised learning framework \cite{raissi2019physics}, starting with the generation of input points $Y_\Omega = \{ \bm{y}_m \in \Omega \}_{m=1}^{M}$ that are sampled uniformly at random inside the physical domain. Then, for every source point $\bm{y}_m\in\Omega$, the corresponding training datasets in $\bm{x}$-variable
	\begin{equation*}
		X_{\textnormal{Reglr}} = \big\{ \bm{x}_n^{\Omega} \big\}_{n=1}^{N_{\textnormal{R}}},\ \ \ X_{\textnormal{Bndry}} = \big\{ \bm{x}_n^{\partial\Omega} \big\}_{n=1}^{N_{\textnormal{B}}},\ \ \  \textnormal{and} \ \ \ X_{\textnormal{Snglr}} = \big\{ \bm{x}_n^{\Gamma} \big\}_{n=1}^{N_{\textnormal{S}}}, 
	\end{equation*}
	are sampled uniformly at random from the regular region $\Omega \setminus \Gamma$, boundaries $\partial\Omega$ and $\partial B_\epsilon(\bm{x}) $, respectively. Here, $N_{\textnormal{R}}$, $N_{\textnormal{B}}$, and $N_{\textnormal{S}}$ indicate the batch sizes of training datasets $X_{\textnormal{Reglr}}$, $X_{\textnormal{Bndry}}$, and $X_{\textnormal{Snglr}}$, respectively. Notably, in the one-dimensional case \eqref{normalization-condition-embedded-1D}, we have $X_{\textnormal{Snglr}} = Y_\Omega$ by sending the radius parameter $\epsilon\to 0$ in \eqref{normalization-condition-embedded}.

	\begin{remark}
		Note that the value of our augmented variable $z=\varphi( \bm{x}, \bm{y} )$ in the singularity-encoded Green's function is uniquely determined at each collocation point $( \bm{x}, \bm{y} )$, which implies that our representation of Green's function is embedded on hyperplanes within a one-order higher dimensional space rather than being globally defined (see \autoref{fig-SgEncdGreen-1D-ToyModel-3D} for example). As a result, the total number of input data points (or, equivalently, the execution times of forward and backward propagations) remains unchanged regardless of the increased dimensionality.
	\end{remark}
	
	%% algorithm %
	%%--------------------------------------%
	\begin{figure}[t!]
		\vspace*{-0.8cm}
		\begin{algorithm}[H]
			\caption{Singularity-Encoded Neural Networks for Learning Green's Functions}
			\fontsize{10}{12}\selectfont
			\begin{algorithmic}
				\STATE{\% \textit{Preparation} }
				\STATE{-- generate training datasets $Y_\Omega$, $X_{\textnormal{Reglr}}$, $X_{\textnormal{Bndry}}$, and $X_{\textnormal{Snglr}}$;}
				\STATE{-- calculate the augmented input $\varphi(\bm{x}, \bm{y})$ for each sample point;}
				%\STATE{-- compute the value of augmented variable for each sample;}
				
				\STATE{\% \textit{Training Process}}
				\STATE{-- construct and initialize the neural network model $\widehat{G}( \bm{x}, \bm{y}, z; \theta)$;}
				\WHILE{the maximum number of epochs is not reached}
				\STATE{-- draw mini-batch data uniformly at random from training datasets;}
				\STATE{-- network training on shuffled datasets with a suitable learning rate and penalty coefficients: }
				\STATE{
					\vspace{-.28cm}
					\begin{equation*}
						\theta^* = \operatorname*{arg\, min}_\theta\  L_{\textnormal{Reglr}} (\theta) + \beta_{\textnormal{Snglr}} L_{\textnormal{Snglr}}(\theta) + \beta_{\textnormal{Bndry}} L_{\textnormal{Bndry}} (\theta) + \beta_{\textnormal{Symtr}} L_{\textnormal{Symtr}}( \theta );
					\end{equation*}
					\vspace{-0.4cm}
				}
				\ENDWHILE
				
				\STATE{\% \textit{Testing Process} }
				\STATE{-- forward pass of the trained neural network model:}
				\STATE{
					\vspace{-.28cm}
					\begin{equation*}
						\widecheck{G}(\bm{x},\bm{y} ) = \frac12  \widehat{G}( \bm{x},\bm{y}, \varphi(\bm{x},\bm{y}); \theta^*) + \frac12  \widehat{G}( \bm{y}, \bm{x}, \varphi(\bm{y}, \bm{x}); \theta^*).
					\end{equation*}
					\vspace{-0.32cm}
				}
			\end{algorithmic}
			\label{Algorithm-SgEncdGreen}
		\end{algorithm}
		\vspace{-0.6cm}
	\end{figure}
	%%--------------------------------------%
	
	According to the residual of equations \eqref{BVP-ExactGreen-Encoded}, our empirical loss functions are defined as follows
	\begin{equation*}
		\begin{array}{c}
			\displaystyle  L_{\textnormal{Reglr}} (\theta) = \frac{1}{M} \sum_{m=1}^M \bigg( \frac{1}{N_{\textnormal{R}}} \sum_{n=1}^{N_{\textnormal{R}}} 
			\bigg| \nabla c(\bm{x}_n^{\Omega}) \cdot \big( \nabla_{\bm{x}} + \nabla \varphi( \bm{x}_n^{\Omega}, \bm{y}_m ) \partial_z \big) \widehat{G}( \bm{x}_n^{\Omega}, \bm{y}_m, \varphi( \bm{x}_n^{\Omega}, \bm{y}_m ) ;\theta) + c(\bm{x}_n^{\Omega}) \Big( \Delta_{\bm{x}} \, + \ \ \ \\[-0.1cm]
			\displaystyle \qquad \qquad \qquad \qquad \qquad \qquad 2 \nabla \varphi( \bm{x}_n^{\Omega}, \bm{y}_m ) \cdot \nabla_{\bm{x}} \, \partial_z  + \|\nabla \varphi( \bm{x}_n^{\Omega}, \bm{y}_m ) \|^2 \partial_{zz}  \Big) \widehat{G}( \bm{x}_n^{\Omega}, \bm{y}_m, \varphi( \bm{x}_n^{\Omega}, \bm{y}_m ) ;\theta)  \bigg|^2 \bigg), \\[0.3cm]
			\displaystyle L_{\textnormal{Snglr}}^{d\geq 2} (\theta) = \frac{1}{M} \sum_{m=1}^M \bigg( 1 + \frac{2\pi \epsilon}{N_{\textnormal{S}}} \sum_{n=1}^{N_{\textnormal{S}}}  c(\bm{x}_n^{\Gamma}) \big( \nabla_{\bm{x}} + \nabla \varphi( \bm{x}_n^{\Gamma}, \bm{y}_m ) \partial_z \big) \widehat{G}( \bm{x}_n^{\Gamma}, \bm{y}_m, \varphi( \bm{x}_n^{\Gamma}, \bm{y}_m ) ;\theta) \cdot \frac{\bm{x}_n^\Gamma - \bm{y}_m}{ \lVert \bm{x}_n^\Gamma - \bm{y}_m \rVert}  \bigg)^2, \\
			\displaystyle L_{\textnormal{Bndry}} (\theta) =  \frac{1}{M} \sum_{m=1}^M \bigg( \frac{1}{N_{\textnormal{B}}} \sum_{n=1}^{N_{\textnormal{B}}} \big|  \widehat{G}( \bm{x}_n^{\partial \Omega}, \bm{y}_m, \varphi( \bm{x}_n^{\partial\Omega}, \bm{y}_m ) ;\theta) \big|^2 \bigg),
		\end{array}
	\end{equation*}
	where the normalization condition \eqref{normalization-condition-embedded} induced loss function can be simplified in the one dimension, namely,
	\begin{equation*}
		L_{\textnormal{Snglr}}^{d=1} (\theta) = \frac{1}{M} \sum_{m=1}^M \left( \big| 2 c(y_m) \partial_z\widehat{G}( y_m, y_m, \varphi(y_m, y_m);\theta) + 1 \big|^2 \right).
	\end{equation*}
	Moreover, by introducing a soft-constraint term to enforce the symmetry of Green's function away from $\bm{x}=\bm{y}$,
	\begin{equation*}
		L_{\textnormal{Symtr}}( \theta ) = \frac{1}{M} \sum_{m=1}^{M}  \Bigg( \frac{1}{N_{\textnormal{R}}} \sum_{n=1}^{N_{\textnormal{R}}} \left| \widehat{G}( \bm{x}_n^{\Omega}, \bm{y}_m, \varphi( \bm{x}_n^{\Omega}, \bm{y}_m) ) - \widehat{G}( \bm{y}_m, \bm{x}_n^{\Omega}, \varphi( \bm{y}_m, \bm{x}_n^{\Omega} ) ) \right|^2 \Bigg),
	\end{equation*}
	the learning task associated with our singularity-encoded Green's function takes on the form
	\begin{equation*}
		\theta^* = \operatorname*{arg\, min}_\theta\  L_{\textnormal{Reglr}} (\theta) + \beta_{\textnormal{Snglr}} L_{\textnormal{Snglr}}(\theta) + \beta_{\textnormal{Bndry}} L_{\textnormal{Bndry}} (\theta) + \beta_{\textnormal{Symtr}} L_{\textnormal{Symtr}}( \theta ),
	\end{equation*}
	where $\beta_{\textnormal{Snglr}}$, $\beta_{\textnormal{Bndry}}$, and $\beta_{\textnormal{Symtr}}>0$ are user-defined penalty parameters. Upon the completion of training phase, our approximate Green's function, denoted as $\widecheck{G}(\bm{x},\bm{y})$ throughout this paper, is obtained by projecting the trained neural network model \eqref{SgEncdGreen-NN-Model} back onto the original physical domain, that is,
	\begin{equation}
		\widecheck{G}(\bm{x},\bm{y} ) = \frac12  \widehat{G}( \bm{x},\bm{y}, \varphi(\bm{x},\bm{y}); \theta^*) + \frac12  \widehat{G}( \bm{y}, \bm{x}, \varphi(\bm{y}, \bm{x}); \theta^*)
		\label{SgEncdGreen-Projected}
	\end{equation}
	where the averaging operation is utilized as a hard constraint to guarantee the symmetry of Green's function. The implementation details are presented in Algorithm \ref{Algorithm-SgEncdGreen}, with hyperparameter configurations reported in \autoref{Sec-Experiments}. 
	
\begin{figure}[!t]
		\begin{subfigure}{0.30\textwidth}
			\centering
			\includegraphics[width=\textwidth]{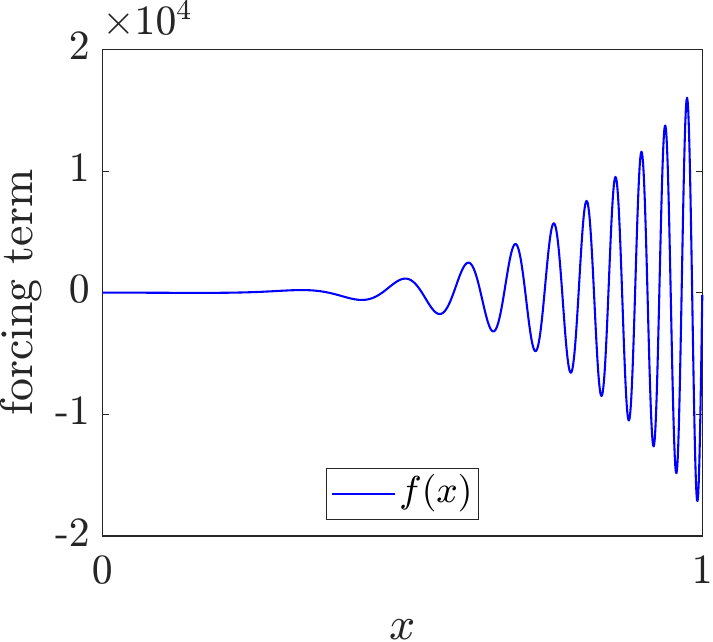}
			\caption{$f(x)$}
			\label{fig-SgEncdGreen-1D-ToyModel-forcing-term}
		\end{subfigure}
		\hspace{0.25cm}
		\begin{subfigure}{0.298\textwidth}
			\centering
			\includegraphics[width=\textwidth]{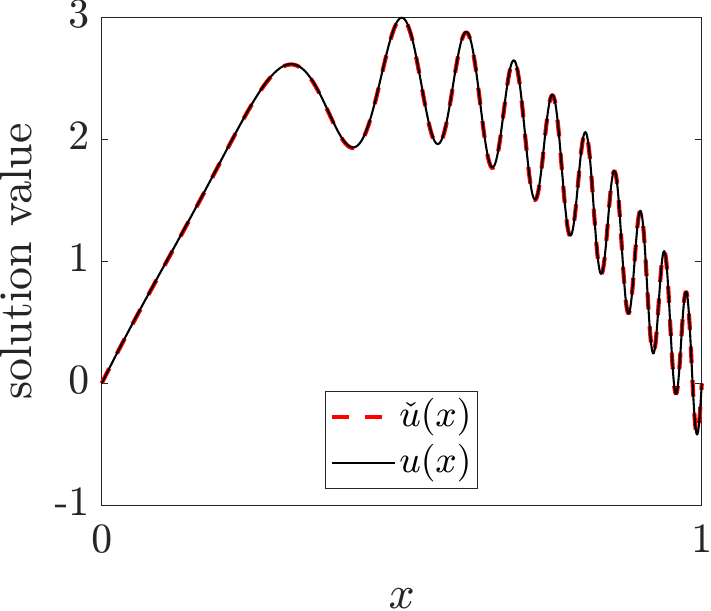}
			\caption{$u(x)$ and $\check{u}(x)$}
			\label{fig-SgEncdGreen-1D-ToyModel-solution-u}
		\end{subfigure}
		\hspace{0.25cm}
		\begin{subfigure}{0.297\textwidth}
			\centering
			\includegraphics[width=\textwidth]{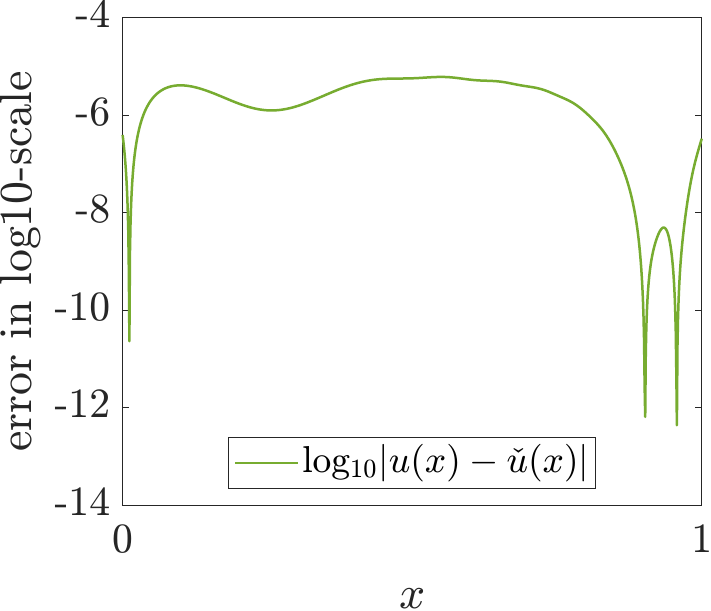}
			\caption{$\textnormal{log}_{10}|u(x) - \check{u}(x)|$}
			\label{fig-SgEncdGreen-1D-ToyModel-PtErr-u}	
		\end{subfigure}	
		
		\caption{Application of the singularity-encoded Green's function as a fast solver for the illustrative example \eqref{BVP-1D-ToyModel}.}
		\label{Experiment-1D-ToyModel-fast-solver}
		\vspace{-0.2cm}
	\end{figure}

	As a direct result, our approximate Green's function \eqref{SgEncdGreen-Projected} can be used for computing the solution of problem \eqref{BVP-ExactSolu} for arbitrary forcing terms, namely, through the numerical integration of \eqref{BVP-ExactSolu-Representation-ExactGreen}
	\begin{equation}
		u(\bm{x}) \approx \check{u}( \bm{x} ) = \sum_{T_p\in\mathcal{T}_h}   \sum_{(\omega_y,\bm{y}_q)} \omega_y f(\bm{y}_q) \widecheck{G}( \bm{x}, \bm{y}_q) ,
		\label{BVP-NumSolu-Representation-SgEncdGreen}	
	\end{equation}
	in which $\mathcal{T}_h$ indicates a triangulation of domain $\Omega$ with the mesh size $h>0$, and $(\omega_y,\bm{y}_q)$ the quadrature weights and points associated with the $p$-th element $T_p\in\mathcal{T}_h$. Notably, \eqref{BVP-NumSolu-Representation-SgEncdGreen} only requires a forward pass of our trained neural network model, enabling real-time efficiency in the numerical solution of boundary value problem \eqref{BVP-ExactSolu}. For instance, we choose the forcing term $f(x)=20-60\pi x\cos(20\pi x^3)+1800\pi^2 x^4\sin(20\pi x^3)$ for our model problem \eqref{BVP-1D-ToyModel}, where the exact solution is given by $u(x)=10x-10x^2+0.5\sin(20\pi x^3)$. The computational results shown in \autoref{fig-SgEncdGreen-1D-ToyModel-solution-u} and \autoref{fig-SgEncdGreen-1D-ToyModel-PtErr-u} demonstrate that our approach \eqref{BVP-NumSolu-Representation-SgEncdGreen} achieves good approximation accuracy without requiring the numerical discretization of the differential equation \eqref{BVP-1D-ToyModel}.
	
	% figure
	%%--------------------------------------%
	%%--------------------------------------%	
	
%%%%%%%%%%%%%%%%%%%%%%%%%%%%%%%%%%%%%%%%%%%%%%%%%%%%%
\subsection{Spectral Bias of Singularity-Encoded Green's Functions}\label{Section-Spectral-Bias}
	
	Despite its computational efficiency as a fast solver \eqref{BVP-NumSolu-Representation-SgEncdGreen} for elliptic boundary value problems with diverse source terms, the spectral bias \cite{rahaman2019spectral,xu2024overview} of our singularity-encoded Green's functions \eqref{SgEncdGreen-Projected} is addressed in this section, which motivates the design of hybrid iterative method in the next section. To begin with, we revisit the Mercer's theorem \cite{hsing2015theoretical} for representing measurable kernels, which is also known as the Karhunen-Lo$\acute{\textnormal{e}}$ve expansion \cite{frauenfelder2005finite}.
		
	\begin{lemma} \label{Lemma-Mercer}
		Suppose that the kernel function $K(\bm{x},\bm{y})\in L^2(\Omega\times\Omega)$ is symmetric and non-negative definite, then its corresponding integral operator
		\begin{equation*}
			\mathcal{K}: L^2(\Omega)\to L^2(\Omega),\ \ \ (\mathcal{K}\phi)(\bm{x}) := \int_{\Omega} K(\bm{x},\bm{y})\phi(\bm{y})\, d \bm{y}\ \ \ \textnormal{for any}\ \phi\in L^2(\Omega)
		\end{equation*}
		is continuous, self-adjoint and Hilbert-Schmidt, admiting a countable sequence $\{ \mu_j, \phi_j(\bm{x}) \}_{j\geq 1}$ of eigenpairs with $\mathbb{R} \ni \mu_j \searrow 0$ as $j \to \infty$. Moreover, the kernel $K(\bm{x},\bm{y})$ can be represented in the series expansion 
		\begin{equation}
			K(\bm{x},\bm{y}) = \sum_{j=1}^\infty \mu_j \phi_j(\bm{x}) \phi_j(\bm{y})
			\label{Thm-Mercer-Expansion}
		\end{equation}
		that converges in the $L^2$ norm \cite{steinwart2012mercer}, namely, the truncation of expansion \eqref{Thm-Mercer-Expansion} after first $N\in\mathbb{N}_+$ terms (indicated by $K_N(\bm{x},\bm{y})$ in what follows) satisfies the optimality condition
		\begin{equation*}
			\left\lVert K(\bm{x},\bm{y}) - K_N(\bm{x},\bm{y}) \right\rVert^2_{L^2(\Omega\times\Omega)} = \Big\lVert K(\bm{x},\bm{y}) - \sum_{j=1}^N \mu_j \phi_j(\bm{x}) \phi_j(\bm{y}) \Big\rVert^2_{L^2(\Omega\times\Omega)} = \sum_{j=N+1}^\infty \mu_j^2.
		\end{equation*}
	\end{lemma}
	
	Recall that for dimensions $1\leq d\leq 3$, our singularity-encoded Green’s function is symmetric by construction \eqref{SgEncdGreen-Projected} and $L^2$-integrable due to our inclusion of weakly singular augmented variables \eqref{Singularity-SgEncdGreen}. Then, by Lemma \ref{Lemma-Mercer}, the eigenvalue problem admits a well-posed variational form: find $0\neq \mu_j\in \mathbb{R}$ and $0 \neq \phi_j(\bm{x}) \in L^2(\Omega)$ such that
	\begin{equation*}
		\int_\Omega \int_\Omega \widecheck{G}(\bm{x},\bm{y}) \phi_j(\bm{y}) \psi(\bm{x})\, d\bm{y}d\bm{x} = \mu_j \int_\Omega \phi_j(\bm{x}) \psi(\bm{x})\, d\bm{x}\ \ \ \textnormal{for any}\ \psi(\bm{x})\in L^2(\Omega).
	\end{equation*}
	This enables spectral analysis of our singularity-encoded Green’s function through the Galerkin approximation: find $0\neq \check{\mu}_j\in\mathbb{R}$ and $\check{\phi}_j(\bm{x}) \in V_h\subset L^2(\Omega)$ such that
	\begin{equation}
		\int_\Omega \int_\Omega \widecheck{G}(\bm{x},\bm{y}) \check{\phi}_j(\bm{y}) \psi_h(\bm{x})\, d\bm{y}d\bm{x} = \check{\mu}_j \int_\Omega \check{\phi}_j(\bm{x}) \psi_h(\bm{x})\, d\bm{x}\ \ \ \textnormal{for any}\ \psi_h(\bm{x})\in V_h,
		\label{EigenProblem-SgEncdGreen}
	\end{equation}
	where $V_h$ denotes the finite element space, e.g., the space of linear or quadratic Lagrange elements. Here, the singular value $\widecheck{G}(\bm{x},\bm{x})$ is approximated by averaging over some balls centered at $\bm{x}\in\Omega$, in accordance with the Lebesgue differentiation theorem \cite{folland1999real} for the local integrable Green's functions ($d=2$ and $3$).
	
	% figure
	%%--------------------------------------%
	\begin{figure}[!t]
		\centering
		\begin{subfigure}[b]{0.31\textwidth}
			\centering
			\includegraphics[width=\textwidth]{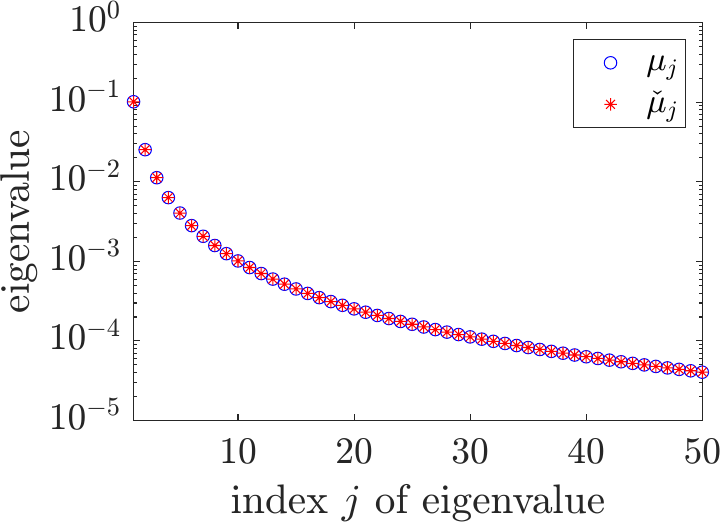}
			\caption{$\mu_j$, $\check{\mu}_j$ with $1\leq j\leq 50$}
			\label{fig-SgEncdGreen-1D-ToyModel-EigVal-1-50}	
		\end{subfigure}
		\hspace{0.25cm}
		\begin{subfigure}{0.312\textwidth}
			\centering
			\includegraphics[width=\textwidth]{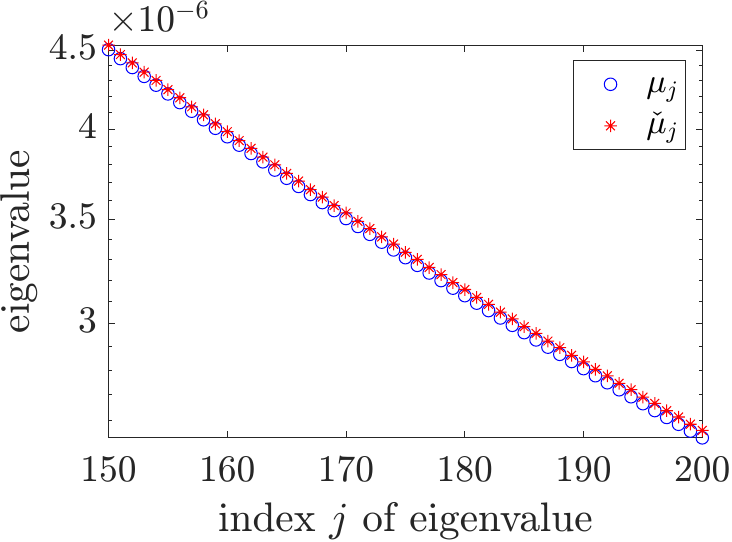}
			\caption{$\mu_j$, $\check{\mu}_j$ with $201\leq j\leq 250$}
			\label{fig-SgEncdGreen-1D-ToyModel-EigVal-201-250}	
		\end{subfigure}
		\hspace{0.25cm}
		\begin{subfigure}{0.303\textwidth}
			\centering
			\includegraphics[width=\textwidth]{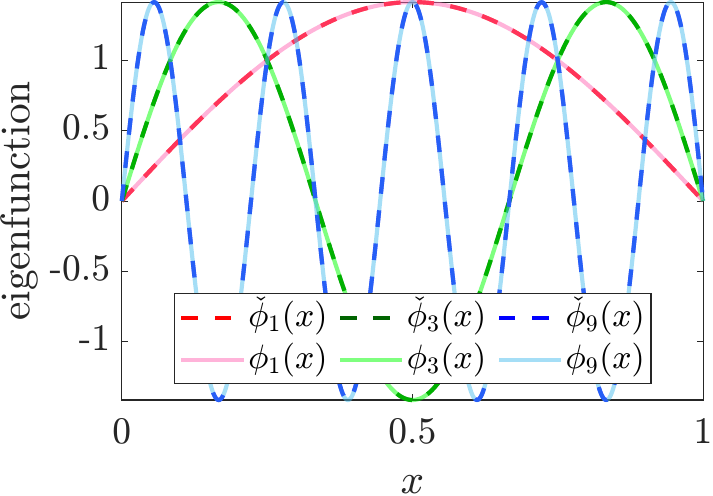}
			\caption{$\phi_j(x)$, $\check{\phi}_j(x)$ with $j=1$, $3$, $9$ }
			\label{fig-SgEncdGreen-1D-ToyModel-EigFct}	
		\end{subfigure}
		
		\vspace{0.3cm}
		
		\begin{subfigure}{0.30\textwidth}
			\centering
			\includegraphics[width=\textwidth]{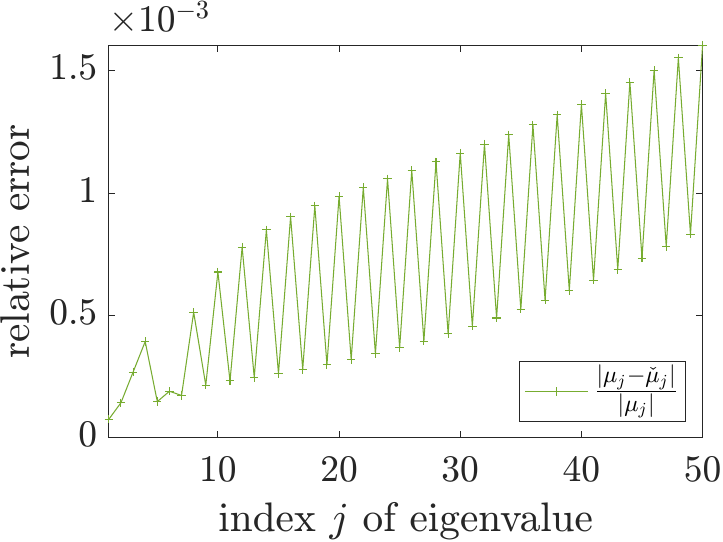}
			\caption{$\delta_{\mu_j}$ with $1\leq j\leq 50$}
			\label{fig-SgEncdGreen-1D-ToyModel-EigVal-RelErr-1-50}	
		\end{subfigure}
		\hspace{0.45cm}
		\begin{subfigure}{0.30\textwidth}
			\centering
			\includegraphics[width=\textwidth]{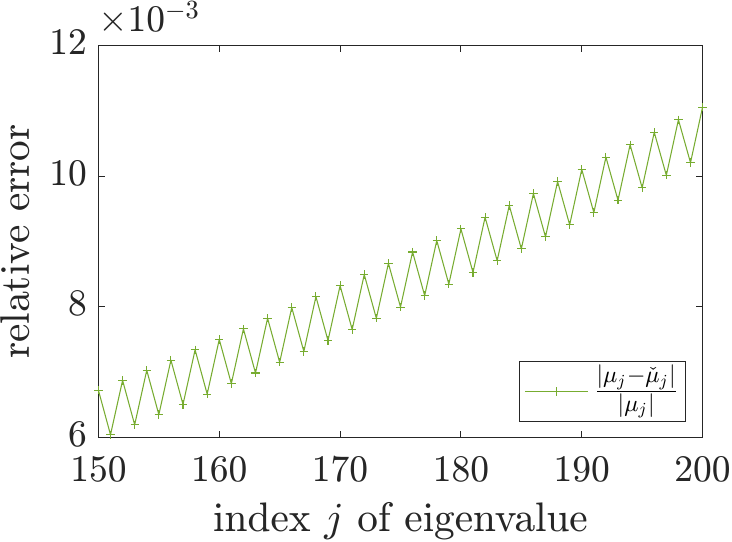}
			\caption{$\delta_{\mu_j}$ with $150\leq j\leq 200$}
			\label{fig-SgEncdGreen-1D-ToyModel-EigVal-RelErr-201-250}		
		\end{subfigure}
		\hspace{0.45cm}
		\begin{subfigure}{0.302\textwidth}
			\centering
			\includegraphics[width=\textwidth]{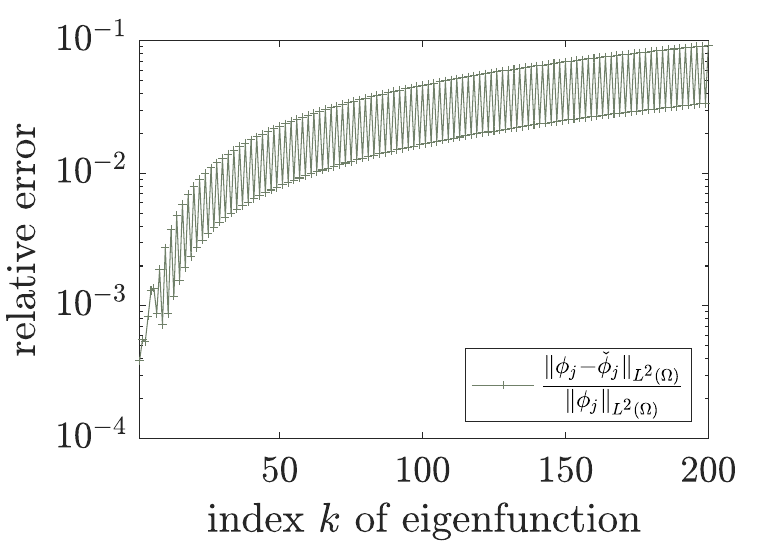}
			\caption{$\delta_{\phi_j}$ with $1\leq j\leq 200$}	
			\label{fig-SgEncdGreen-1D-ToyModel-EigFct-RelErr}		
		\end{subfigure}
		
		\caption{Numerical eigenpairs of our singularity-encoded Green's function for the illustrative example \eqref{BVP-1D-ToyModel}.}
		\label{Experiment-1D-ToyModel-Eigenpairs}
		
		\vspace{-0.2cm}
	\end{figure}
	%%--------------------------------------%
	
	The spectral properties of the singularity-encoded Green’s function, as characterized numerically through the eigenvalue problem \eqref{EigenProblem-SgEncdGreen} of our illustrative example \eqref{BVP-1D-ToyModel} in \autoref{Experiment-1D-ToyModel-Eigenpairs}, exhibit notable consistency across variable-coefficients and higher-dimensional settings (see \autoref{Sec-Experiments}). Specifically, we adopt quadratic elements with regular mesh width $h=2^{-10}$ for solving the eigenvalue problem \eqref{EigenProblem-SgEncdGreen}, with relative errors being defined as 
	\begin{equation*}
		\delta_\mu = \frac{ |  \check{\mu}_j - \mu_j | }{ |\mu_j|} \ \ \ \textnormal{and}\ \ \ \delta_\phi = \frac{ \lVert \check{\phi}_j(x) - \phi_j(x) \rVert_{L^2(\Omega)} }{ \lVert \phi_j(x) \rVert_{L^2(\Omega)} }.
	\end{equation*}
	As evidenced in 	\autoref{fig-SgEncdGreen-1D-ToyModel-EigVal-RelErr-1-50}, \autoref{fig-SgEncdGreen-1D-ToyModel-EigVal-RelErr-201-250}, and \autoref{fig-SgEncdGreen-1D-ToyModel-EigFct-RelErr}, our singularity-encoded Green's function exhibits markedly better accuracy in the low-frequency regimes compared to its high-frequency counterparts, which aligns with the spectral bias \cite{rahaman2019spectral} (or frequency principle \cite{xu2025understanding}) observed and studied across diverse applications. 
	
	While the spectral bias inherent in our singularity-encoded Green's function might appear as a deficiency, we show in section \ref{Sec-Hybrid-IteM} that it can be utilized to accelerate the convergence of classical iterative methods, particularly for low-frequency components. In fact, the presence of spectral bias persists in the deep operator network \cite{lu2021learning} or the Fourier neural operator \cite{li2020fourier}, motivating the design of efficient numerical methods through hybrid iterative schemes \cite{zhang2024blending} or preconditioning strategies \cite{kopanivcakova2025deeponet}.		
	
%%%%%%%%%%%%%%%%%%%%%%%%%%%%%%%%%%%%%%%%%%%%%%%%%%%%%
%%%%%%%%%%%%%%%%%%%%%%%%%%%%%%%%%%%%%%%%%%%%%%%%%%%%%

%%%%%%%%%%%%%%%%%%%%%%%%%%%%%%%%%%%%%%%%%%%%%%%%%%%%%
%%%%%%%%%%%%%%%%%%%%%%%%%%%%%%%%%%%%%%%%%%%%%%%%%%%%%
\section{Iterative Methods Enhanced by Singularity-Encoded Green's Functions}\label{Section-Acltd-IteM}
	
	In this section, we focus on the discrete system of elliptic problem \eqref{BVP-ExactSolu}, where the pre-trained singularity-encoded Green's function is utilized to improve the convergence speed of iterative solvers, either by constructing effective preconditioners in section \ref{Sec-Neural-Preconditioner} or designing efficient hybrid iterative schemes in section \ref{Sec-Hybrid-IteM}. On the one hand, our trained surrogate model inspires the development of both dense preconditioning matrices for small-scale systems and the two-level overlapping additive Schwarz preconditioners for large-scale linear systems. On the other hand, hybrid iterative strategies are presented by leveraging the spectral bias of trained models, followed by integrating into the geometric multigrid framework to accelerate convergence. The effectiveness of our proposed methods is validated via the one-dimensional example \eqref{BVP-1D-ToyModel-SgEncdGreen} for illustration, while more examples are presented in \autoref{Sec-Experiments}.	

%%%%%%%%%%%%%%%%%%%%%%%%%%%%%%%%%%%%%%%%%%%%%%%%%%%%%	
\subsection{Preconditioned System via Singularity-Encoded Green's Functions}\label{Sec-Neural-Preconditioner}
	
	Recall that the exact Green's function serves as the inverse of a differential operator in the continuous setting \eqref{BVP-ExactSolu-Representation-ExactGreen}. Discretizing the elliptic boundary value problem \eqref{BVP-ExactSolu} via finite difference or finite element methods \cite{leveque2007finite,brenner2008mathematical} on a domain triangulation $\mathcal{T}_h$ with mesh size $h>0$ gives rise to a linear algebraic system 
	\begin{equation}
		AU=F
		\label{BVP-MatrixForm}
	\end{equation}	
	where $A\in\mathbb{R}^{n\times n}$ denotes a symmetric positive definite matrix, $F\in\mathbb{R}^n$ the right-hand side, and $U\in\mathbb{R}^n$ the vector of unknowns. To improve the convergence of iterative solvers, it is common to consider a preconditioned system 
	\begin{equation*}
		BAU = BF
	\end{equation*}		
	where a preconditioner $B\in\mathbb{R}^{n\times n} \approx A^{-1}$ is incorporated to improve the spectral properties of \eqref{BVP-MatrixForm}. Note that $A^{-1}$ corresponds to a discrete counterpart of the exact Green's function \cite{leveque2007finite}, hence a promising preconditioner can be constructed by evaluating our singularity-encoded Green's functions \eqref{SgEncdGreen-Projected} at the grid points $\{\bm{x}_i\}_{i=1}^n$
	\begin{equation}
		(\widecheck{B})_{ij} = \widecheck{G}(\bm{x}_i, \bm{x}_j),\quad \text{for}\ 1\le i, j\le n.
		\label{DiscreteSgnEcdGreen}
	\end{equation} 
	Though the construction of preconditioner \eqref{DiscreteSgnEcdGreen} remains valid under mesh refinement, the global support property of the Green's function \eqref{BVP-ExactGreen} would result in a dense matrix \cite{hartmann2012green} that may hinder computational efficiency \cite{varga1962iterative}. As such, the preconditioning matrix \eqref{DiscreteSgnEcdGreen} may not be suitable for solving large-scale systems, as both the computational cost and the accumulation of round-off errors become difficulty to manage.
		
	To address large-scale linear systems \eqref{BVP-MatrixForm}, we consider an additive preconditioner (not relabelled for simplicity) within the framework of two-level overlapping Schwarz methods \cite{toselli2004domain,kopanivcakova2025deeponet}. To be specific, let $\{ \Omega_\ell \}_{\ell=1}^L$ denote an overlapping decomposition of the computational domain and let $R_\ell$ be the matrix representation of restriction operator on the subdomain $\Omega_\ell$ \cite{toselli2004domain,kopanivcakova2025deeponet}, the local matrix associated with $\Omega_\ell$ can then be formulated as 
	\begin{equation*}
		A_\ell = R_\ell A R_\ell^T
	\end{equation*}	
	for $1\leq \ell\leq L$. On the other hand, given a shape-regular coarse mesh $\mathcal{T}_H$ on the domain, the restriction matrix $R_0$ for its corresponding coarse space can be defined analogously \cite{toselli2004domain,golub2013matrix}, as well as the coarse matrix $A_0 = R_0 A R_0^T$. It is also noteworthy that the coarse matrix can be built by employing discretization methods on the coarse mesh. Consequently, the two-level overlapping additive Schwarz preconditioner \cite{toselli2004domain} takes on the form	
	\begin{equation}
		B = R_0^T A_0^{-1} R_0 + \sum_{\ell=1}^L R_\ell^T A_\ell^{-1}R_\ell,
	\end{equation}
	in which the coarse solver $A_0^{-1}$ provides a global correction mechanism that complements local solvers $\{ A_\ell^{-1} \}_{\ell=1}^L$ by handling low-frequency errors that extend across multiple subdomains.
	
	Particularly, our coarse solver is constructed by applying the same routine \eqref{DiscreteSgnEcdGreen} but for the coarse grid $\{\bm{x}_i'\}_{i=1}^{n'}$ of triangulation $\mathcal{T}_H$, namely, $(\widecheck{B}_0)_{ij} = \widecheck{G}(\bm{x}_i', \bm{x}_j')$ for $1\le i,j\le n'$, giving rise to an additive preconditioner of the form
	\begin{equation}
		\widecheck{B} = R_0^T \widecheck{B}_0 R_0 + \sum_{\ell=1}^L R_\ell^T A_\ell^{-1}R_\ell,
		\label{Preconditioner-LargeScale}
	\end{equation}	
	which mitigates both computational complexity and storage requirements when handling large-scale systems. 
	
	For the ease of illustration, we refer to \eqref{DiscreteSgnEcdGreen} and \eqref{Preconditioner-LargeScale} as the neural preconditioner \cite{cui2025hybrid} throughout this work, as they are not only established through a forward pass of the trained neural network model \eqref{SgEncdGreen-Projected} but also inherit the spectral bias of our singularity-encoded Green's functions. Finally, the preconditioned system
	\begin{equation}
		\widecheck{B}AU=\widecheck{B}F 
		\label{PreconditionedSystem}
	\end{equation}
is solved using the biconjugate gradient method \cite{varga1962iterative}, with implementation details provided in Algorithm \autoref{Alg-SgEncd-Green-Preconditioner}.
	
	%--------------------------%		
	\begin{algorithm}[t!] 
		\caption{Preconditioned System via Singularity-Encoded Green's Functions}
		\fontsize{10}{12}\selectfont
		\begin{algorithmic}[1]
			\STATE{\% \textit{Preparation} }			
			\STATE{ -- discretize problem \eqref{BVP-ExactSolu} over the triangulation $\mathcal{T}_h$ to obtain a discrete linear system $A U = F$;}
			\STATE{ -- employ Algorithm \autoref{Algorithm-SgEncdGreen} to learn the singularity-encoded Green's function for problem \eqref{BVP-ExactSolu};}
			
			\STATE{\% \textit{Neural Preconditioner} }									
			\IF{small-scale linear systems}
			\STATE{ -- forward pass of the trained model with input samples $\{\bm{x}_i\}_{i=1}^n$ taken from the fine mesh $\mathcal{T}_h$:}			
			\STATE{\vspace{0.1cm}				
				\hfill $(\widecheck{B})_{ij} = \widecheck{G}(\bm{x}_i, \bm{x}_j),\quad \text{for}\ 1\le i,j\le n;$ \hfill
			\vspace{0.1cm}	}
			\ELSIF{large-scale linear systems}
			\STATE{ -- start with a coarse mesh $\mathcal{T}_H$, triangulate the domain $\Omega$ into overlapping subdomains $\{ \Omega_\ell \}_{\ell=1}^L$;}	
			\STATE{ -- assemble local coefficient matrices $\{ A_\ell \}_{i=1}^L$ and their associated restriction matrices $\{ R_\ell \}_{\ell=0}^L$;}		
			\STATE{ -- forward pass of the trained model with input samples $\{\bm{x}_i'\}_{i=1}^{n'}$ taken from the coarse mesh $\mathcal{T}_H$:}			
			\STATE{\vspace{0.04cm}				
				\hfill $\displaystyle \widecheck{B} = R_0^T \widecheck{B}_0 R_0 + \sum_{\ell=1}^L R_\ell^T A_\ell^{-1}R_\ell\ \ $ with $\ \ (\widecheck{B}_0)_{ij} = \widecheck{G}(\bm{x}_i', \bm{x}_j')\ \text{for}\ 1\le i,j\le n';$ \hfill
			\vspace{0.04cm}	}
			\ENDIF						
			\STATE{\% \textit{Preconditioned System} }			
			\STATE{ -- employ the biconjugate gradient method to solve the preconditioned linear system:}			
			\STATE{\vspace{0.15cm}				
				\hfill $\widecheck{B}AU=\widecheck{B}F.$ \hfill
			\vspace{0.08cm}	}
		\end{algorithmic}
		\label{Alg-SgEncd-Green-Preconditioner}
	\end{algorithm}	
	%--------------------------%
	
%%%%%%%%%%%%%%%%%%%%%%%%%%%%%%%%%%%%%%%%%%%%%%%%%%%%%
\subsubsection{Example: Poisson's Equation in One Dimension}\label{Section-1D-ToyModel-Preconditioner}
	
	To validate the effectiveness of our neural preconditioner \eqref{DiscreteSgnEcdGreen} and \eqref{Preconditioner-LargeScale}, we study in this section the linear system arising from the discretization of our illustrative example \eqref{BVP-1D-ToyModel} on a uniform grid $\displaystyle \{x_i = ih\}_{i=0}^{n+1}$ with the mesh size $h=\frac{1}{n+1}$ \cite{leveque2007finite}. For instance, when using linear finite elements \cite{brenner2008mathematical}, the discrete system \eqref{BVP-MatrixForm} takes on the form
		\begin{equation}
		A = \frac{1}{h}
		\begin{bmatrix}
		2 & -1 & & & \\
		-1 & 2 & -1 & & \\
		& \ddots & \ddots & \ddots & & \\
		& & -1 & 2 & -1 \\
		& & & -1 & 2
		\end{bmatrix},
		\qquad
		U = \begin{bmatrix}
		U_1\\
		U_2\\
		\vdots\\
		U_{n-1}\\
		U_n
		\end{bmatrix},
		\qquad
		F = h \begin{bmatrix}
		F_1 \\
		F_2\\
		\vdots\\
		F_{n-1}\\
		F_n
		\end{bmatrix},
		\label{BVP-1D-ToyExmp-MatrixForm}			
		\end{equation}
		where $\displaystyle F_i = \int_0^1 f(x)\phi_i(x) dx$ for $1\leq i\leq n$. To mitigate the ill-conditioning of the stiffness matrix in \eqref{BVP-1D-ToyExmp-MatrixForm}, a variety of preconditioning strategies have been developed in the literature. It is also noteworthy that for this example, the inverse of the stiffness matrix coincides with the evaluation of the exact Green’s function at grid points \cite{leveque2007finite}.				
		
	Here, we employ the singularity-encoded Green's function pre-trained in section \ref{Section-1D-ToyModel} to construct our neural preconditioner for linear system \eqref{BVP-1D-ToyExmp-MatrixForm}. Spectral properties of our preconditioned matrix $\widecheck{B}A$ are quantified by the ratio of its largest to smallest eigenvalues, namely,
	\begin{equation*}
		\kappa(\widecheck{B}A)=||\widecheck{B}A||_2 ||(\widecheck{B}A)^{-1}||_2 = \dfrac{\lambda_{\text{max}}}{\lambda_{\text{min}}},		
	\end{equation*}
	which reveals the clustering behavior of eigenvalues that governs the convergence speed of iterative methods.
	
	%-------------------------%
	\begin{table}[!t]
	\centering
	\renewcommand\arraystretch{2} 
	\begin{tabular}{|c|c|c|c|c|}
		\hline
		mesh size & eigenvalues of preconditioned matrix & condition number & $\#$ iterations & $\ell_2$-norm error \\
		\hline
		$h = 2^{-6}$ &
		\begin{minipage}[b]{0.3\columnwidth}
			\centering
			\raisebox{-.4\height}{\includegraphics[width=\linewidth]{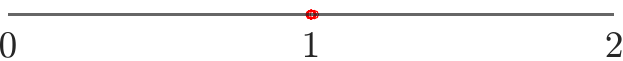}}			
		\end{minipage}
		& \makecell[c]{$\kappa (\widecheck{B}A) = 1.01$\\ $\kappa (A) = 1.66 \times 10^3$} &\makecell[c]{$6$\\$(63)$} &\makecell[c]{$6.54 \times 10^{-15}$\\$(8.78 \times 10^{-15})$} \\
		\hline
		$h = 2^{-8}$ &
		\begin{minipage}[b]{0.3\columnwidth}
			\centering
			\raisebox{-.4\height}{\includegraphics[width=\linewidth]{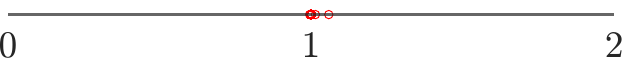}}
		\end{minipage}
		& \makecell[c]{$\kappa (\widecheck{B}A) = 1.06$ \\ $ \kappa (A) = 2.66 \times 10^4$} &\makecell[c]{$6$\\$(255)$} &\makecell[c]{$3.27 \times 10^{-14}$\\$(1.63 \times 10^{-14})$}\\
		\hline
		$h = 2^{-10}$ &
		\begin{minipage}[b]{0.3\columnwidth}
			\centering
			\raisebox{-.4\height}{\includegraphics[width=\linewidth]{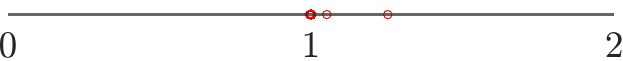}}
		\end{minipage}
		& \makecell[c]{$\kappa (\widecheck{B}A) = 1.26$ \\ $ \kappa (A) = 4.25 \times 10^5$} &\makecell[c]{$6$\\$(1023)$} &\makecell[c]{$7.60 \times 10^{-13}$\\$(6.91 \times 10^{-13})$}\\
		\hline
		$h = 2^{-12}$ &
		\begin{minipage}[b]{0.3\columnwidth}
			\centering
			\raisebox{-.4\height}{\includegraphics[width=\linewidth]{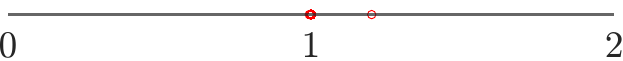}}
		\end{minipage}
		& \makecell[c]{$\kappa (\widecheck{B}A) = 2.04$ \\ $ \kappa (A) = 6.80 \times 10^6$} &\makecell[c]{$6$\\$(4095)$} &\makecell[c]{$4.89 \times 10^{-12}$\\$(4.82 \times 10^{-12})$}\\
		\hline
	\end{tabular}
	\caption{Comparison of the biconjugate gradient method with and without neural preconditioning for the discrete linear system \eqref{BVP-MatrixForm} with stiffness matrix \eqref{BVP-1D-ToyExmp-MatrixForm} and $f(x)=20-60\pi x\cos(20\pi x^3)+1800\pi^2 x^4\sin(20\pi x^3)$.}
	\label{GreenPreconditioner-Poisson1DExp1}
	\vspace{-0.3cm}
	\end{table}
	%-------------------------%
	
	For small-scale systems, our neural preconditioner is constructed directly through formula \eqref{DiscreteSgnEcdGreen}, with \autoref{GreenPreconditioner-Poisson1DExp1} comparing the performance of the biconjugate gradient method with and without preconditioning. Obviously, the condition number of coefficient matrix in problem \eqref{BVP-1D-ToyExmp-MatrixForm} scales as $\kappa(A) = \mathcal{O}(h^2)$, while that of our preconditioned matrix is close to $\kappa(\widecheck{B}A)\approx \mathcal{O}(1)$ regardless of the mesh refinement. As a result, the number of iterations required for convergence becomes mesh-independent after preconditioning, in contrast to the original problems that grows as the mesh gets refined (see the last two columns in \autoref{GreenPreconditioner-Poisson1DExp1}, where parentheses denotes results obtained without preconditioning). Our iteration for preconditioned system is terminated when the error, defined by the $\ell_2$-norm of the difference between the iterative and exact solutions, matches that of the vanilla biconjugate gradient method.
	
	%-------------------------%
	\begin{table}[!b]
	\centering
	\renewcommand\arraystretch{2} 
	\begin{tabular}{|c|c|c|c|c|c|}
		\hline
		mesh size & coarse grid & overlap & condition number & $\#$ iterations & $\ell_2$-norm error  \\
		\hline
		$h=2^{-16}$ & $H=2^{-10}$ & $\delta = 30h$ & $\kappa (\widecheck{B}A) = 7.26$ &$23$ &$7.73 \times 10^{-10}$\\
		\hline
		$h=2^{-18}$ & $H=2^{-10}$ & $\delta = 30h$ & $\kappa (\widecheck{B}A) = 5.51$ &$23$ &$2.5 \times 10^{-8}$\\
		\hline
	\end{tabular}
	\caption{Numerical results of our two-level additive Schwarz preconditioner \eqref{Preconditioner-LargeScale} for the discrete linear preconditioned system \eqref{PreconditionedSystem} with stiffness matrix \eqref{BVP-1D-ToyExmp-MatrixForm}. }
	\label{GreenPreconditioner-Twolevel-Poisson1DExp1}
	\vspace{-0.3cm}	
	\end{table}
	%-------------------------%	
	
	For large-scale systems, the use of a dense preconditioning matrix \eqref{DiscreteSgnEcdGreen} would render both computational cost and memory demands prohibitive, along with a	 notable accumulation of round-off errors. Therefore, the additive Schwarz preconditioner \eqref{Preconditioner-LargeScale} is deployed, in which the evaluation of singularity-encoded Green's function \eqref{SgEncdGreen-Projected} is restricted on the coarse grid to obtain $\widecheck{B}_0$ while the local problems $\{ A_\ell^{-1} \}_{\ell=1}^L$ are of small size and directly solved. Numerical results are displayed in \autoref{GreenPreconditioner-Twolevel-Poisson1DExp1}, which demonstrate the effectiveness of our neural preconditioner \eqref{Preconditioner-LargeScale} as the mesh gets further refined (the unpreconditioned counterpart is omitted as it becomes impractical to solve). 	
%%%%%%%%%%%%%%%%%%%%%%%%%%%%%%%%%%%%%%%%%%%%%%%%%%%%%

%%%%%%%%%%%%%%%%%%%%%%%%%%%%%%%%%%%%%%%%%%%%%%%%%%%%%
\subsection{Hybrid Iterative Methods with Classical and Neural Preconditioners}\label{Sec-Hybrid-IteM}
	
	Motivated by the spectral bias/frequency principle of singularity-encoded Green's functions (see section \ref{Section-Spectral-Bias}), we construct hybrid strategies that enhance classical iterative methods by effectively attenuating their low-frequency errors. More specifically, our hybrid iterative methods are formulated in a residual-correction form \cite{cui2025hybrid,hu2025hybrid}: 		
	\begin{equation}
		U^{[k+1]} = \left\{
		\begin{array}{cl}			
			\displaystyle  U^{[k]} + \widecheck{B} ( F - A U^{[k]} ), \ \ & \ \ \textnormal{if}\ k\equiv 0\ (\textnormal{mod}\ K), \\
			\displaystyle  U^{[k]} + B ( F - A U^{[k]} ), \ \ & \ \ \textnormal{otherwise}, \\
		\end{array}\right.
		\label{Hybrid-IteM-General}
	\end{equation}
	where $k\in\mathbb{N}_+$ is the iteration index, $B\in\mathbb{R}^{n\times  n}$ the conventional preconditioner \cite{varga1962iterative} designed to approximate $A^{-1}$ and $\widecheck{B}\in\mathbb{R}^{n\times  n}$ the preconditioning matrix deduced from our singularity-encoded Green's function \eqref{DiscreteSgnEcdGreen}. Here, the parameter $K$ serves as the switching period, specifying how often the iteration alternates between the two update schemes presented in \eqref{Hybrid-IteM-General}.
	
	As a direct consequence of its construction \eqref{DiscreteSgnEcdGreen}, the neural preconditioner $\widecheck{B}$ used in \eqref{Hybrid-IteM-General} inherits the spectral bias of our pre-trained singularity-encoded Green's function, which is promising in reducing low-frequency error components compared to its high-frequency counterparts. Then, our neural preconditioner is periodically applied to attenuate low-frequency errors that are difficult for conventional preconditioners to eliminate \cite{cui2025hybrid,zhang2024blending}, while the high-frequency error introduced by our neural preconditioner can be rapidly reduced through classical iterations.						
	
	Accordingly, the error propagation associated with our hybrid iterative method \eqref{Hybrid-IteM-General} can be expressed as
	\begin{equation*}
	U - U^{[\ell K]} = (I - BA)^{K-1} (I - \widecheck{B}A) ( U - U^{[(\ell-1) K]} ) = \big[ (I - BA)^{K-1} (I - \widecheck{B}A) \big]^\ell ( U - U^{[0]} )
	\end{equation*}
	for any $\ell\in\mathbb{N}$, which decays to zero if and only if the amplification matrix satisfies
	\begin{equation*}
	\rho\big( (I - BA)^{K-1} (I - \widecheck{B}A) \big) < 1,
	\end{equation*}
	where $\rho(\cdot)$ represents the spectral radius \cite{hu2025hybrid,brenner2008mathematical}. Therefore, a formal convergence analysis of our hybrid iterative scheme \eqref{Hybrid-IteM-General} can be established, while implementation details are summarized in Algorithm \autoref{Alg-SgEncd-Green-Hybrid-IteM}.
	
	%--------------------------%		
	\begin{algorithm}[t!] 
		\caption{Hybrid Iterative Methods with Classical and Neural Preconditioners}
		\fontsize{10}{12}\selectfont
		\begin{algorithmic}[1]
			\STATE{\% \textit{Preparation} }			
			\STATE{ -- discretize problem \eqref{BVP-ExactSolu} over the triangulation $\mathcal{T}_h$ to obtain a discrete linear system $A U = F$;}
			\STATE{ -- employ Algorithm \autoref{Algorithm-SgEncdGreen} to learn the singularity-encoded Green's function for problem \eqref{BVP-ExactSolu};}
			\STATE{ -- construct the preconditioner \eqref{DiscreteSgnEcdGreen} or \eqref{Preconditioner-LargeScale} via a forward pass of the trained neural network;}			

			\STATE{\% \textit{Iterative Procedure} }
			\STATE{ -- set the initial guess $U^{[1]}\in\mathbb{R}^{n}$ to a vector of zeros;}
			\FOR{$k = 1, 2, 3, \cdots$ }
			\WHILE{ $U^{[k]}$ does not satisfy the stopping criteria}
			\IF{$ k \equiv 0$ (mod $K=2$)}
			\STATE{\vspace{0.08cm}				
				$U^{[k+1]} = U^{[k]} + \widecheck{B} ( F - A U^{[k]})$			
			\vspace{0.02cm}	} \COMMENT{\textit{reduce low-frequency errors by our neural preconditioner}}
			\ELSE
			\STATE{\vspace{0.08cm}				
				$U^{[k+1]} = U^{[k]} + B ( F - A U^{[k]})$			
			\vspace{0.02cm}	} \COMMENT{\textit{reduce high-frequency errors by classical preconditioners}}
			\ENDIF						
			\ENDWHILE			
			\ENDFOR
		\end{algorithmic}
		\label{Alg-SgEncd-Green-Hybrid-IteM}
	\end{algorithm}	
	%--------------------------%
		
	\begin{lemma} \textnormal{\cite{hu2025hybrid}} \label{Thm-Hybrid-IteM-Convergence}
	There exists a correction period $K\in\mathbb{N}_+$ large enough such that the hybrid iterative method \eqref{Hybrid-IteM-General} converges, if the applied original iterative method converges, i.e., $\rho(I-BA)<1$.
	\end{lemma}			
	
	Notably, the switching or correction period $K$ is assigned relatively large values in existing studies \cite{cui2025hybrid,zhang2024blending,hu2025hybrid}, which seems to be consistent with the convergence result established in Lemma \ref{Thm-Hybrid-IteM-Convergence}. However, low-frequency error components are not addressed until $K-1$ iterations of the traditional method have been completed in each cycle, by which time the solution has already been smoothed by its preceding updates. To alleviate this inefficiency, we propose adopting $K=2$ in our hybrid iterative method \eqref{Hybrid-IteM-General}, in which the classical and neural preconditioners are applied alternatively in a sequential manner. The effectiveness of our statement is further validated via numerical experiments in section \ref{Section-1D-ToyModel-Hybrid-IteM} and \autoref{Sec-Experiments}, where $K=2$ shows superior acceleration compared to other values.	
	
	In addition to accelerating basic iterative schemes such as the damped Jacobi method \cite{varga1962iterative} (see section \ref{Section-1D-ToyModel-Hybrid-IteM}), our hybrid approach \eqref{Hybrid-IteM-General} can also be integrated into the well-established geometric multigrid framework \cite{trottenberg2001multigrid}. To be specific, geometric multigrid methods recursively deploy fine-grid relaxation (e.g., damped Jacobi method) to attenuate high-frequency errors with coarse-grid correction to eliminate low-frequency errors. Although efficient on structured grids, its direct extension to unstructured meshes and anisotropic problems is often hindered by the complexity of generating hierarchical meshes \cite{hu2025hybrid,zhang2024blending}. Alternatively, by substituting the pre- and post-smoothing procedures with our hybrid iterative approach \eqref{Hybrid-IteM-General}, comparable accuracy can be achieved with a reduced number of grid levels, thereby avoiding the construction of nested hierarchical meshes (see section \ref{Section-1D-ToyModel-Hybrid-IteM} and section \ref{Sec-Experiments}).
			 
%%%%%%%%%%%%%%%%%%%%%%%%%%%%%%%%%%%%%%%%%%%%%%%%%%%%%
\subsubsection{Example: Poisson's Equation in One Dimension}\label{Section-1D-ToyModel-Hybrid-IteM}
				
	To demonstrate the efficacy our proposed hybrid iterative scheme \eqref{Hybrid-IteM-General}, we consider its application to the discrete linear system \eqref{BVP-1D-ToyExmp-MatrixForm} arising from the discretization of the illustrative example \eqref{BVP-1D-ToyModel} and compare results with that of the damped Jacobi method \cite{varga1962iterative}. More specifically, we employ the frequently used preconditioner
	\begin{equation*}
		B = \omega D^{-1}
	\end{equation*}
	in which $\omega = \frac12$ denotes a relaxation parameter \cite{varga1962iterative} and $D$ the diagonal part of our stiffness matrix $A$ using mesh size $h=2^{-8}$, to generate iterative solutions of the classical damped Jacobi method, namely,
	\begin{equation}
		U^{[k+1]}_{\textnormal{Jacobi}} = U^{[k]}_{\textnormal{Jacobi}} + \omega D^{-1} ( F - A U^{[k]}_{\textnormal{Jacobi}} ).
		\label{1D-Poisson-Damped-Jacobi}
	\end{equation}
	On the other hand, the same preconditioner is utilized in a hybrid manner within our iterative scheme \eqref{Alg-SgEncd-Green-Hybrid-IteM}
	\begin{equation}
		U^{[k+1]}_{\textnormal{Hybrid}} = \left\{
		\begin{array}{cl}			
			\displaystyle  U^{[k]}_{\textnormal{Hybrid}} + \widecheck{B} ( F - A U^{[k]}_{\textnormal{Hybrid}} ), \ \ & \ \ \textnormal{if}\ k\ \textnormal{is even},\\
			\displaystyle  U^{[k]}_{\textnormal{Hybrid}} + \omega D^{-1} ( F - A U^{[k]}_{\textnormal{Hybrid}} ), \ \ & \ \ \textnormal{if}\ k\ \textnormal{is odd}.
		\end{array}\right.
		\label{1D-Poisson-Hybrid-IteM}
	\end{equation}
	
	%--------------- 1DPoisson-Example iterative methods ---------------%
	\begin{figure}[!b]
		\centering
		\begin{subfigure}{0.60\textwidth}
			\centering
			\includegraphics[width=1\textwidth]{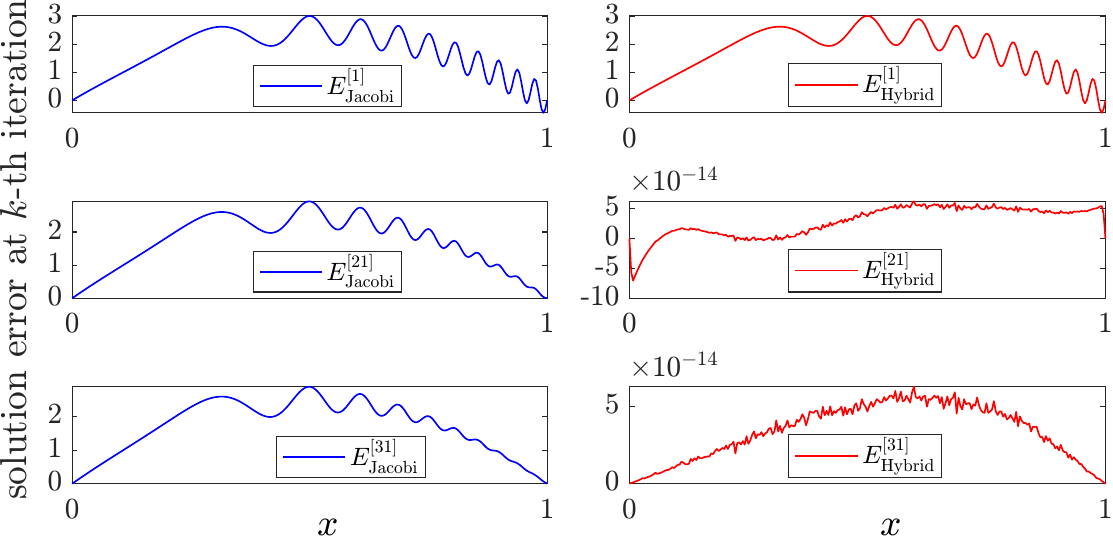}
			\caption{comparison of smoothing effects between methods \eqref{1D-Poisson-Damped-Jacobi} and \eqref{1D-Poisson-Hybrid-IteM}}
			\label{fig-1D-Poisson-Removing-Smoothing-Effects}
		\end{subfigure}
		\hspace{0.3cm}
		\begin{subfigure}{0.34\textwidth}
			\centering
			\includegraphics[width=1\textwidth]{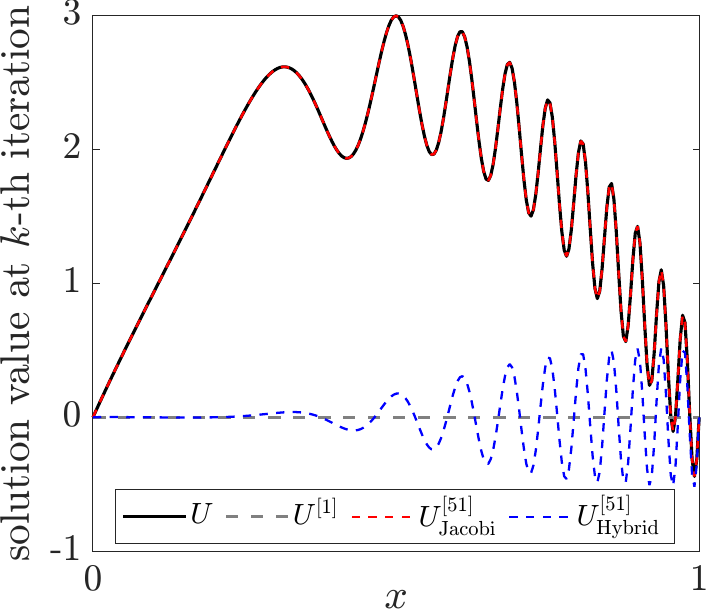}
			\caption{comparison of iterative solutions}
			\label{fig-1D-Poisson-Comparison-Ite50}
		\end{subfigure}	
		
		\vspace{-0.1cm}
			
		\caption{Elimination of low-frequency errors for the discrete linear system \eqref{BVP-1D-ToyExmp-MatrixForm} of the illustrative problem \eqref{BVP-1D-ToyModel}.}
		\label{Experiments-Poisson1D-IteM-1}
		\vspace{-0.3cm}
	\end{figure}
	%--------------- 1DPoisson-Example iterative methods ---------------%		
	
	By defining the solution error, residual error, and mode-wise error at the $k$-th iteration as 
	\begin{equation*}
		E^{[k]} = U - U^{[k]}, \ \ \ R^{[k]} = F - A U^{[k]},\ \ \ \textnormal{and}\ \ \ M^{[k],j} = | \langle E^{[k]}, \xi_j \rangle |,
	\end{equation*}
	where $\xi_j$ denotes the $j$-th eigenvector of the amplification matrix $I\!-\!BA$, we now present a quantitative comparison of \eqref{1D-Poisson-Damped-Jacobi} and \eqref{1D-Poisson-Hybrid-IteM}. As shown in \autoref{fig-1D-Poisson-Removing-Smoothing-Effects}, the solution error computed by the damped Jacobi method \eqref{1D-Poisson-Damped-Jacobi} becomes smooth within the first few iterations, which is known as the smoothing effects whereby high-frequency error components are rapidly eliminated while its low-frequency counterparts persist. Conversely, the solution error of our hybrid iterative scheme \eqref{1D-Poisson-Hybrid-IteM} retains oscillatory behavior as displayed in \autoref{fig-1D-Poisson-Removing-Smoothing-Effects}, owing to the extra high-frequency error introduced by our neural preconditioner. As a direct result, superior alignment with the true solution is observed in the early iterations, as demonstrated in \autoref{fig-1D-Poisson-Comparison-Ite50}.
	
	%--------------- 1DPoisson-Example iterative methods ---------------%
	\begin{figure}[!t]
		\centering		
		\begin{subfigure}{0.47\textwidth}
			\centering
			\includegraphics[width=1\textwidth]{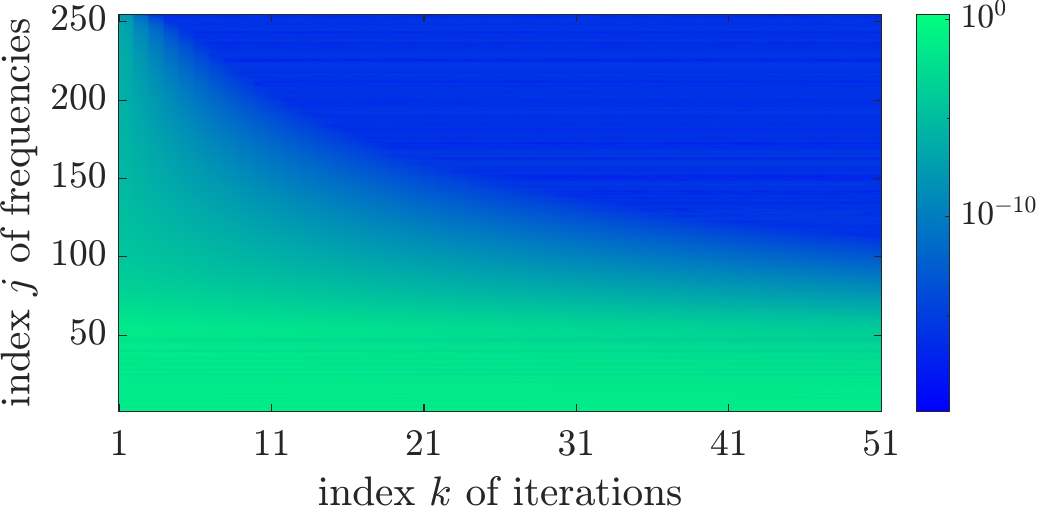}
			\caption{the spectral error of damped Jacobi method \eqref{1D-Poisson-Damped-Jacobi}}
			\label{fig-1D-Poisson-Jacobi-ModeErr}
		\end{subfigure}
		\hspace{0.2cm}		
		\begin{subfigure}{0.47\textwidth}
			\centering
			\includegraphics[width=1\textwidth]{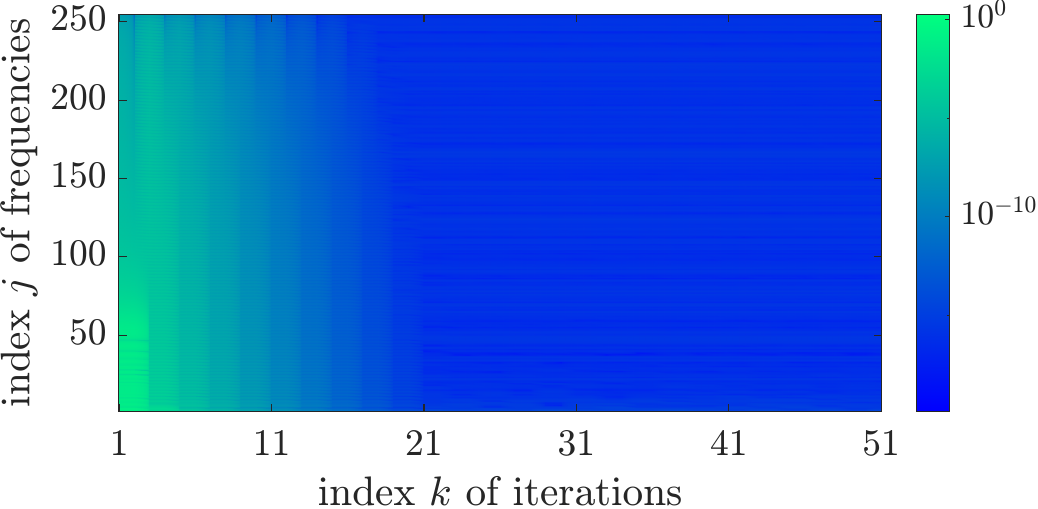}
			\caption{the spectral error of hybrid iterative method \eqref{1D-Poisson-Hybrid-IteM}}
			\label{fig-1D-Poisson-Hybrid-ModeErr}
		\end{subfigure}
						
		\centering
		\begin{subfigure}{0.31\textwidth}
			\centering
			\includegraphics[width=.95\textwidth]{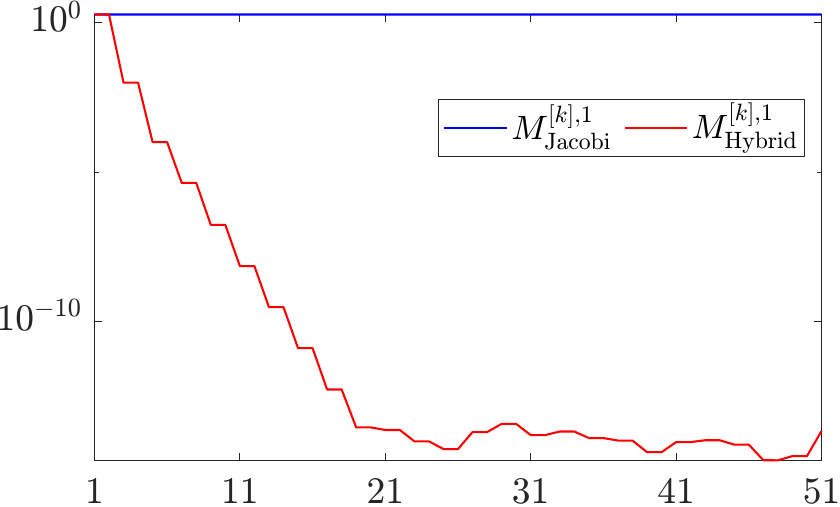}
			\caption{the lowest-frequency error component}
			\label{fig-1D-Poisson-Hybrid-ModeErr-j1} 
		\end{subfigure}
		\hspace{0.05cm}
		\begin{subfigure}{0.31\textwidth}
			\centering
			\includegraphics[width=.95\textwidth]{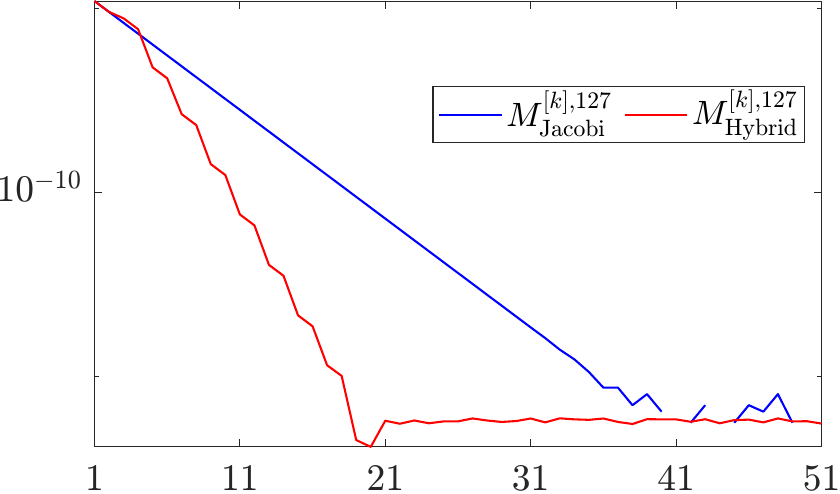}
			\caption{the mid-frequency error component}
			\label{fig-1D-Poisson-Hybrid-ModeErr-j127} 
		\end{subfigure}
		\hspace{0.05cm}
		\begin{subfigure}{0.31\textwidth}
			\centering
			\includegraphics[width=.95\textwidth]{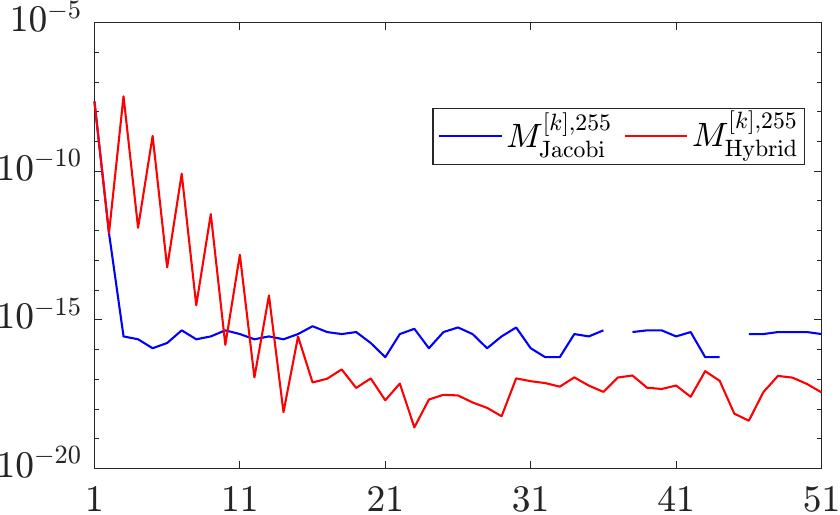}
			\caption{the highest-frequency error component}
			\label{fig-1D-Poisson-Hybrid-ModeErr-j255} 
		\end{subfigure}

		\vspace{-0.1cm}
								
		\caption{Comparison of spectral errors for solving the discrete linear system \eqref{BVP-1D-ToyExmp-MatrixForm} of the illustrative problem \eqref{BVP-1D-ToyModel}.}
		\label{Experiments-Poisson1D-IteM-2}
		\vspace{-0.3cm}
	\end{figure}
	%--------------- 1DPoisson-Example iterative methods ---------------%		
	
	Moreover, \autoref{fig-1D-Poisson-Jacobi-ModeErr} and \autoref{fig-1D-Poisson-Hybrid-ModeErr} depict the mode-wise error $M^{[k],j}$ across all frequencies for both methods, while the lowest-, intermediate-, and highest-frequency error components are individually displayed in \autoref{fig-1D-Poisson-Hybrid-ModeErr-j1}, \autoref{fig-1D-Poisson-Hybrid-ModeErr-j127}, and \autoref{fig-1D-Poisson-Hybrid-ModeErr-j255}, respectively. Though the inclusion of neural preconditioner introduces additional high- frequency errors (see \autoref{fig-1D-Poisson-Hybrid-ModeErr-j255} for instance) as a consequence of its spectral bias, they can be suppressed by our subsequent execution of Jacobi method. After approximately 20 updates, our proposed approach almost achieves machine accuracy, whereas the traditional method still exhibits notable low-frequency errors.
		
\begin{figure}[!b]								
		\centering
		\begin{subfigure}{0.46\textwidth}
			\centering
			\includegraphics[width=1\textwidth]{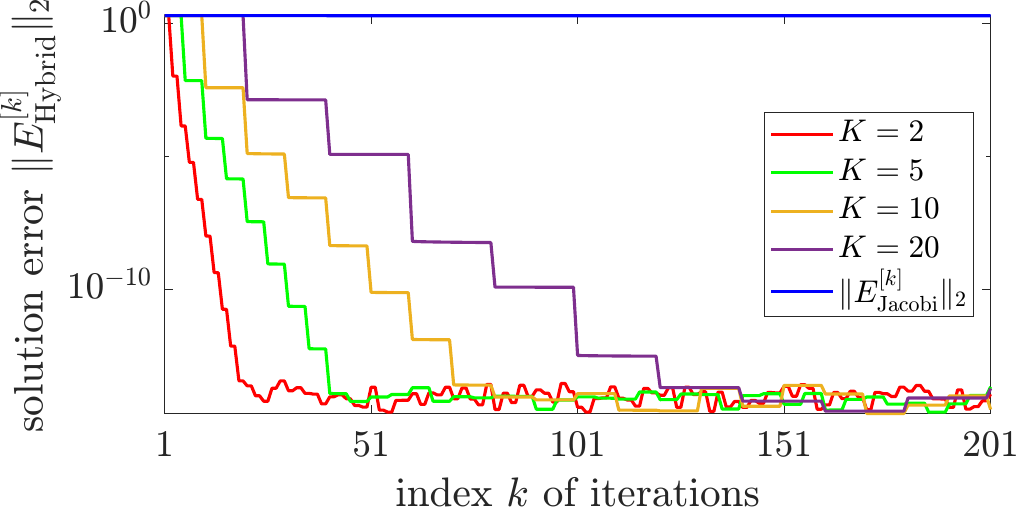}
			\caption{solution error using different switching periods}
			\label{fig-1D-Poisson-Solu-Err} 
		\end{subfigure}
		\hspace{0.3cm}
		\begin{subfigure}{0.46\textwidth}
			\centering
			\includegraphics[width=1\textwidth]{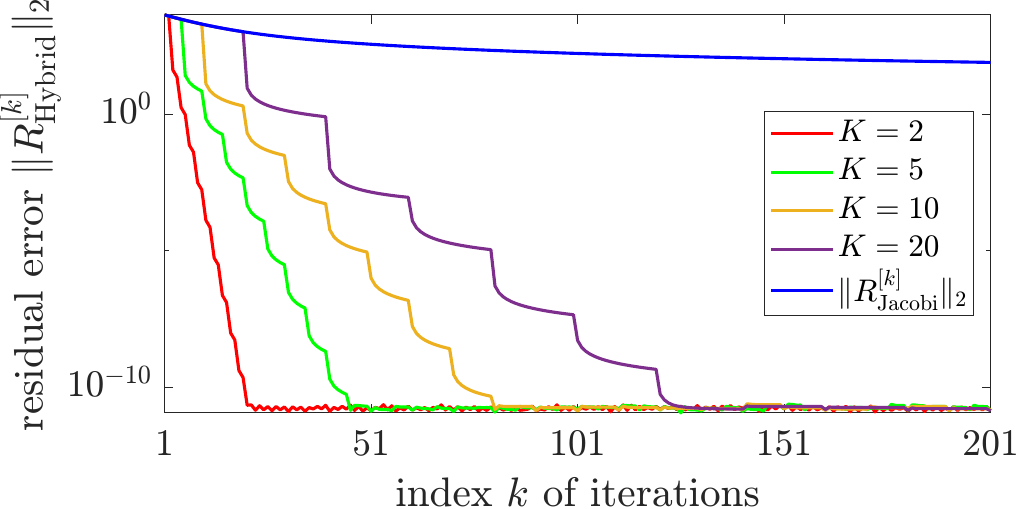}
			\caption{residual error using different switching periods}
			\label{fig-1D-Poisson-Res-Err}
		\end{subfigure}
			
		\caption{Hybrid iteration \eqref{1D-Poisson-Hybrid-IteM} with different switching periods for the discrete system \eqref{BVP-1D-ToyExmp-MatrixForm} of the illustrative problem \eqref{BVP-1D-ToyModel}.}
		\label{Experiments-Poisson1D-IteM-3}
		\vspace{-0.3cm}
	\end{figure}

	In addition, \autoref{fig-1D-Poisson-Solu-Err} and \autoref{fig-1D-Poisson-Res-Err} show the solution error and the residual error for both methods, in which different values of the switching period are utilized in Algorithm \ref{Alg-SgEncd-Green-Hybrid-IteM}. The most rapid decay is attained with $K=2$, wherein the traditional and neural preconditioners are applied alternatively during the iteration. In contrast, using a relatively larger value of $K$, as is common in the existing literature \cite{cui2025hybrid,zhang2024blending,hu2025hybrid}, impedes the prompt reduction of low-frequency errors and thereby degrades the speed of convergence.
	
		%--------------- 1DPoisson-Example iterative methods ---------------%
	%--------------- 1DPoisson-Example iterative methods ---------------%
	
\begin{figure}[!t]
		\centering
		\begin{subfigure}{0.5\textwidth}
			\centering
			\includegraphics[width=\textwidth]{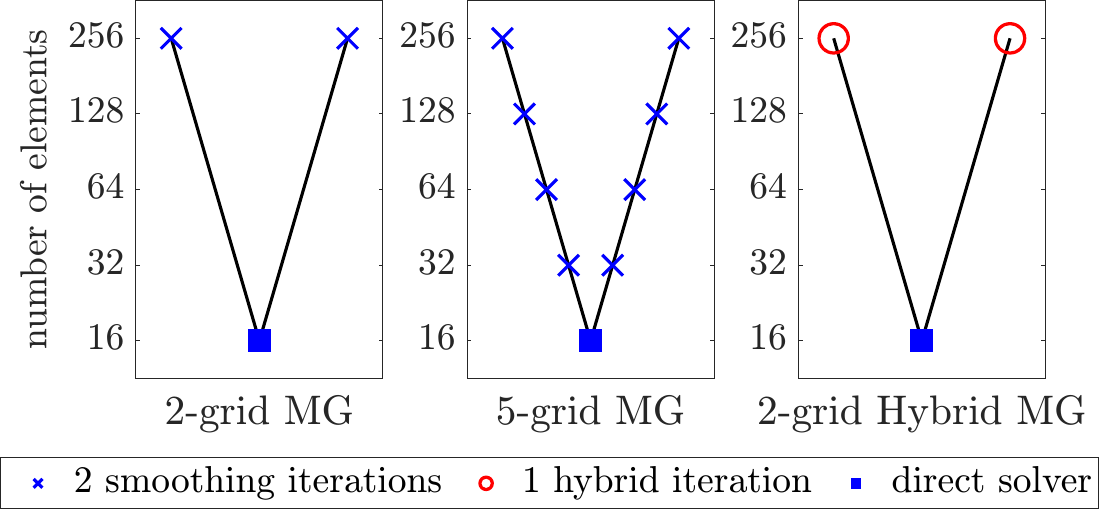}
			\caption{V-cycle structures for different approaches}
			\label{fig-1D-Poisson-MG-V-Cycle}
		\end{subfigure}
		\hspace{0.25cm}
		\begin{subfigure}{0.45\textwidth}
			\centering
			\includegraphics[width=\textwidth]{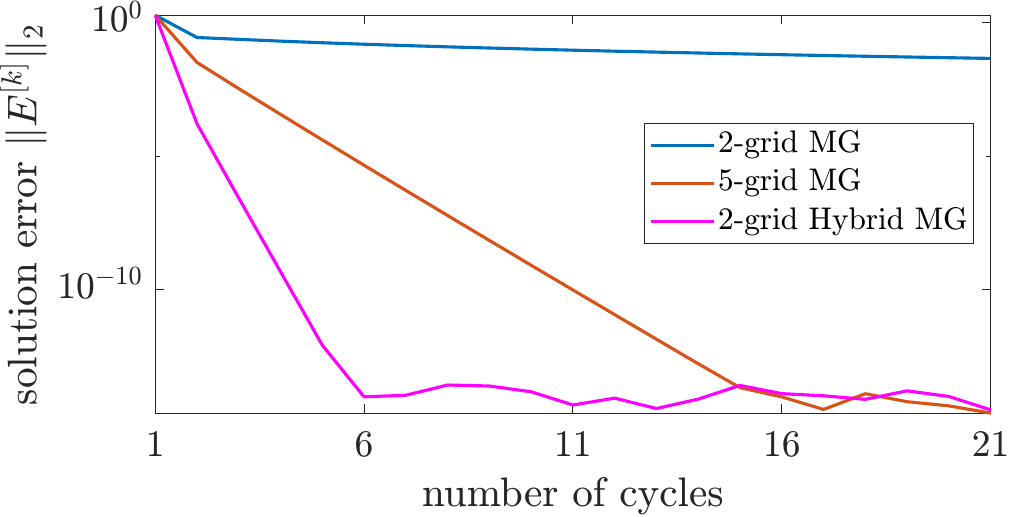}
			\caption{solution errors over iterative cycles}
			\label{fig-1D-Poisson-MG-Solu-Err}
		\end{subfigure}		
		\caption{Comparison of multigrid methods with and without hybrid iterations for solving the discrete linear system \eqref{BVP-1D-ToyExmp-MatrixForm} of the illustrative problem \eqref{BVP-1D-ToyModel}.}
		\label{Experiment-Poisson1D-HybridMG}
		\vspace{-0.3cm}
	\end{figure}

	In addition to accelerating simple iterative methods, our hybrid iterative scheme \eqref{Hybrid-IteM-General} could also be integrated into the multigrid framework \cite{varga1962iterative}, as discussed at the end of section \ref{Sec-Hybrid-IteM}. Specifically, we consider the numerical solution of linear system \eqref{BVP-1D-ToyExmp-MatrixForm} using a mesh size of $h=2^{-8}$, while the coarsest mesh size is set to $h=2^{-4}$ so that a direct solver can be applied. Our comparative study begins with a standard two-grid multigrid approach, whose V-cycle structure is visualized in \autoref{fig-1D-Poisson-MG-V-Cycle} (left). Unfortunately, such a notable coarser gird may fail to represent low-frequency errors inherited from the fine grid, resulting in performance deterioration as depicted in \autoref{fig-1D-Poisson-MG-Solu-Err} \cite{trottenberg2001multigrid}. Therefore, a geometric multigrid approach with a hierarchy of five nested grids (see \autoref{fig-1D-Poisson-MG-V-Cycle}) is employed instead, achieving satisfactory accuracy within 16 V-cycles as shown in \autoref{fig-1D-Poisson-MG-Solu-Err}.
	
	%---------------------------%
	%---------------------------%	
	
	While uniformly refined nested meshes are ideal for multigrid solvers, their construction becomes non-trivial when applied to problems with anisotropic coefficients or complex geometries \cite{zhang2024blending,brenner2008mathematical}. Thanks to the capability of our neural preconditioner in eliminating low-frequency errors, \autoref{fig-1D-Poisson-MG-Solu-Err} implies that comparable accuracy could be achieved within a two-grid framework by employing our hybrid iterative method \eqref{Hybrid-IteM-General} ($K=2$) during the pre- and post-smoothing procedures (see \autoref{fig-1D-Poisson-MG-V-Cycle}). Besides, our hybrid strategy enables convergence to the machine precision within 6 updates of V-cycles, demonstrating accelerated convergence while reducing its dependence on deeply nested grid hierarchies.	
	
%%%%%%%%%%%%%%%%%%%%%%%%%%%%%%%%%%%%%%%%%%%%%%%%%%%%%
%%%%%%%%%%%%%%%%%%%%%%%%%%%%%%%%%%%%%%%%%%%%%%%%%%%%%

%%%%%%%%%%%%%%%%%%%%%%%%%%%%%%%%%%%%%%%%%%%%%%%%%%%%%
%%%%%%%%%%%%%%%%%%%%%%%%%%%%%%%%%%%%%%%%%%%%%%%%%%%%%
\section{Numerical Experiments}\label{Sec-Experiments}		

	%--------------------------------------%
	\begin{table}[!b]
		\renewcommand{\arraystretch}{1.4}
		\setlength{\tabcolsep}{4.5pt}
		\fontsize{9}{12}\selectfont
		\centering
		\begin{tabular}{cccccc}
			\hline 
			& \makecell{Neural Network\\ (Depth, Width)} & \makecell{Penalty Coeff.\\ ($\beta_{\textnormal{Snglr}}$, $\beta_{\textnormal{Bndry}}$, $\beta_{\textnormal{Symtr}}$)} & \makecell{Train Data Size\\ ($N_R$, $N_S$, $N_B$, $M$)} & \makecell{Initial Learning Rate\\ (Decay Milestones)}  \\ \hline
			Poisson Equation \eqref{BVP-1D-ToyModel-ExactGreen} & (2, 40) & (400, 400, 400) & (160, 500, 2, 500) & 0.001 (12$k$, 22$k$)  \\ \cline{1-5}
			Helmholtz Equation \eqref{BVP-1D-Helmholtz} & (2, 40) & (400, 400, 400) & (500, 500, 2, 500) & 0.001 (8$k$, 15$k$) \\ \cline{1-5}
			Poisson Equation \eqref{BVP-2D-Laplacian} & (6, 40) & (400, 400, 400) & (640, 200, 640, 160) & 0.001 (15$k$, 25$k$) \\ \hline 
		\end{tabular}	
		\vspace*{-0.15cm}
		\caption{List of hyper-parameter configurations for learning the Green's function in our experiments.}
		\label{table-hyper-parameters}
		\vspace*{-0.3cm}
	\end{table}
	%%--------------------------------------%
	
	Beyond our illustrative example discussed in section \ref{Section-1D-ToyModel}, \ref{Section-1D-ToyModel-Preconditioner}, and \ref{Section-1D-ToyModel-Hybrid-IteM}, this section reports numerical studies on additional boundary value problems to demonstrate the effectiveness of our proposed methods. To be specific, each experiment starts with the learning of the true Green's function via Algorithm \ref{Algorithm-SgEncdGreen}, followed by its application to accelerating conventional iterative methods, either by serving as a preconditioning matrix in Algorithm \ref{Alg-SgEncd-Green-Preconditioner} or as a component of our hybrid iterative method in Algorithm \ref{Alg-SgEncd-Green-Hybrid-IteM}. Particularly, the indefinite Helmholtz problem is also considered, as it poses greater challenges for classical iterative methods compared to the positive definite cases. 
	
	The introduction of an augmented variable \eqref{Singularity-SgEncdGreen}, coupled with the embedding of the Green's function \eqref{BVP-ExactGreen} into a one-order higher-dimensional space, facilitates the effective representation of singularity during training. Though our singularity-encoded Green's function \eqref{Singularity-SgEncdGreen} manifests increased dimensionality, the curse of dimensionality can be mitigated through the parametrization using fully-connected neural networks. Besides, the generation of input collocation points remains unaffected, as their associated augmented variables can be directly computed through formula \eqref{Singularity-SgEncdGreen}. It is also noteworthy that the hyperbolic tangent activation function is adopted throughout this work, since our proposed methods enable the usage of smooth neural network models to reconstruct singular functions. As before, our trained model is expressed as $\widehat{G} (\bm{x}, \bm{y}, \varphi( \bm{x}, \bm{y}) )$, while its projection \eqref{SgEncdGreen-Projected} is denoted by $\widecheck{G} (\bm{x}, \bm{y})$.	 

	All neural networks are trained via the ADAMW optimizer \cite{kingma2014adam}, using a step decay schedule that reduces the initial learning rate by a factor of 10 after specific milestones. The key hyperparameter configurations adopted in our practical implementations are summarized in \autoref{table-hyper-parameters}, and all experiments are conducted through PyTorch on the Nvidia GeForce RTX 4090 GPU cards. Source codes will be made publicly available upon the acceptance of this paper (\url{https://github.com/q1sun/SgEncd_Green_Function}). % Source codes are publicly accessible on GitHub\footnote{\url{https://github.com/q1sun/SgEncd_Green_Function}}.	
	
%%%%%%%%%%%%%%%%%%%%%%%%%%%%%%%%%%%%%%%%%%%%%%%%%%%%%
\subsection{Green's Function for the Helmholtz Equation with a Variable Coefficient}	\label{Sec-Helmholtz1D}
	
	%%--------------------------------------%
	\begin{figure}[b!]
	\centering
		\begin{subfigure}[b]{0.31\textwidth}
			\centering
			\includegraphics[width=\textwidth]{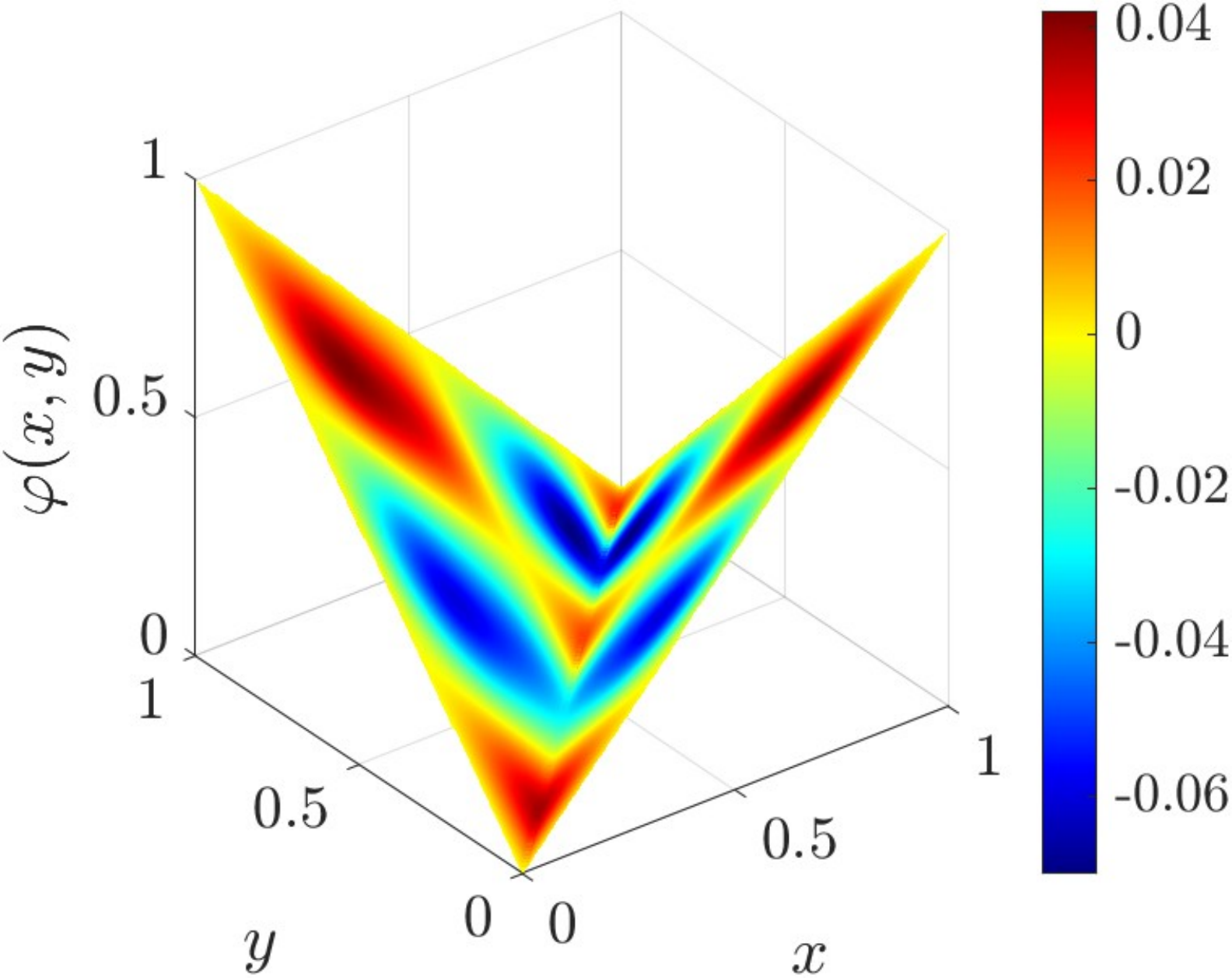}
			\caption{$\widehat{G}(x,y,\varphi(x,y))$}
			\label{fig-SgEncdGreen-1D-VarbHelmholtz-3D}
		\end{subfigure}
		\hspace{0.25cm}
		\begin{subfigure}{0.31\textwidth}
			\centering
			\includegraphics[width=\textwidth]{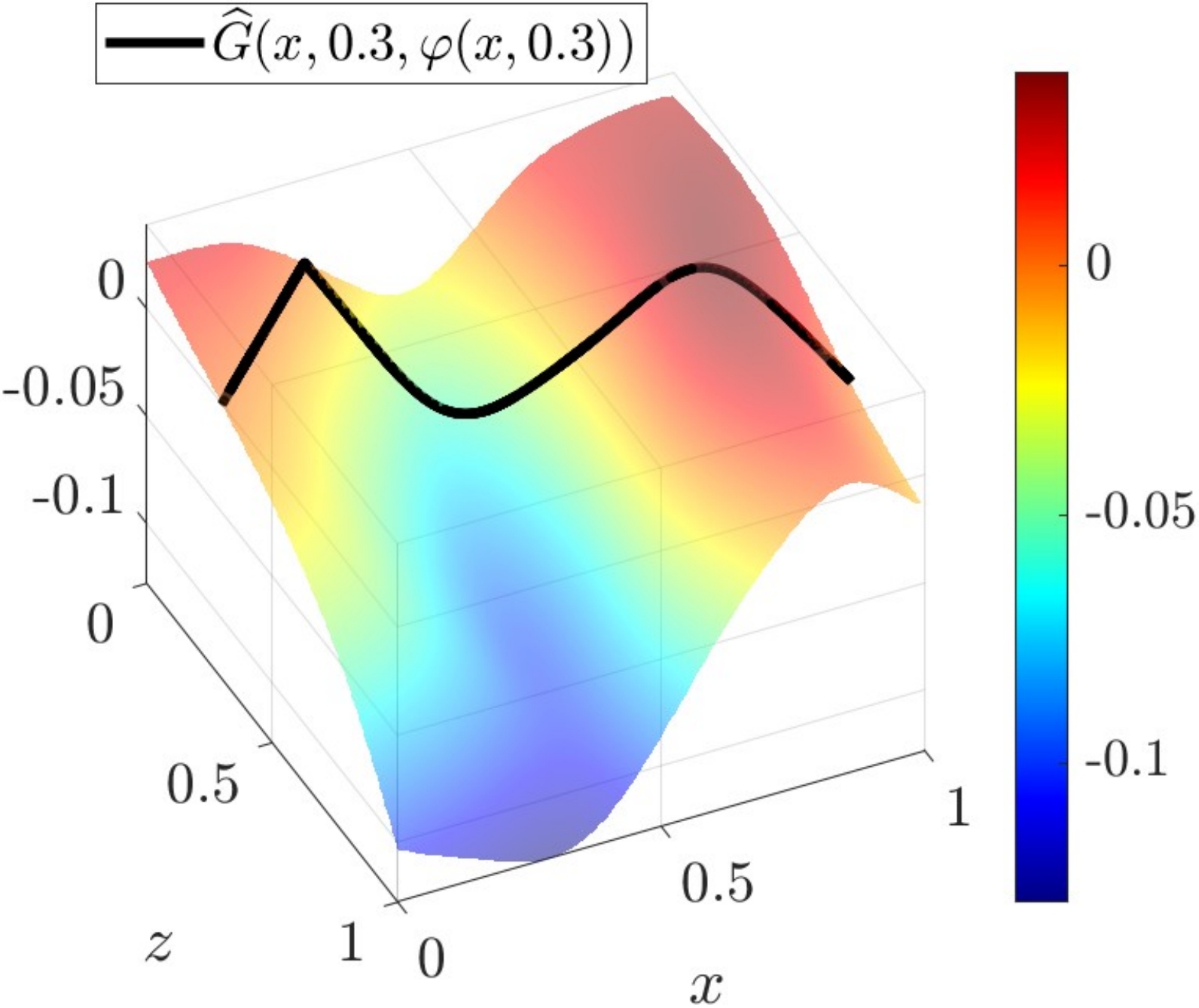}
			\caption{$\widehat{G}(x,0.3,\varphi(x,0.3))$}
			\label{fig-SgEncdGreen-1D-VarbHelmholtz-yfixed-3D}
		\end{subfigure}
		\hspace{0.25cm}
		\begin{subfigure}{0.294\textwidth}
			\centering
			\includegraphics[width=\textwidth]{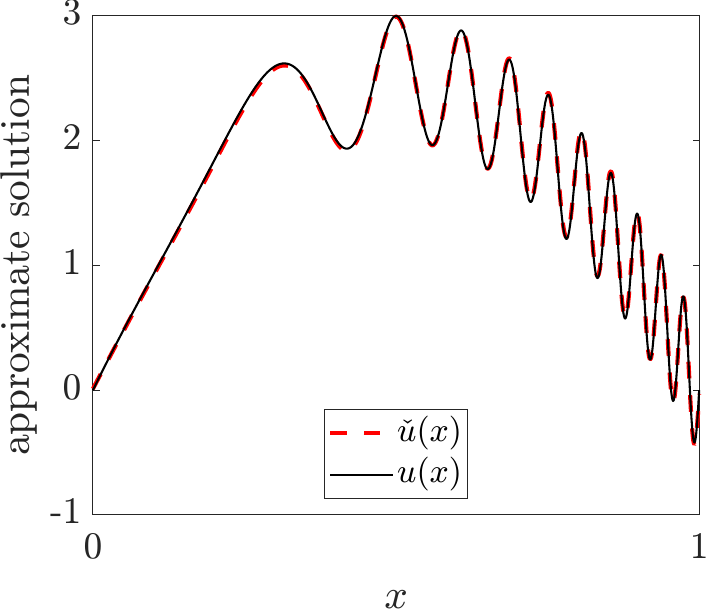}
			\caption{$u(x)$ and $\check{u}(x)$}
			\label{fig-SgEncdGreen-1D-VarbHelmholtz-u}
		\end{subfigure}
		
		\vspace{0.1cm}
		
		\begin{subfigure}{0.31\textwidth}
			\centering
			\includegraphics[width=\textwidth]{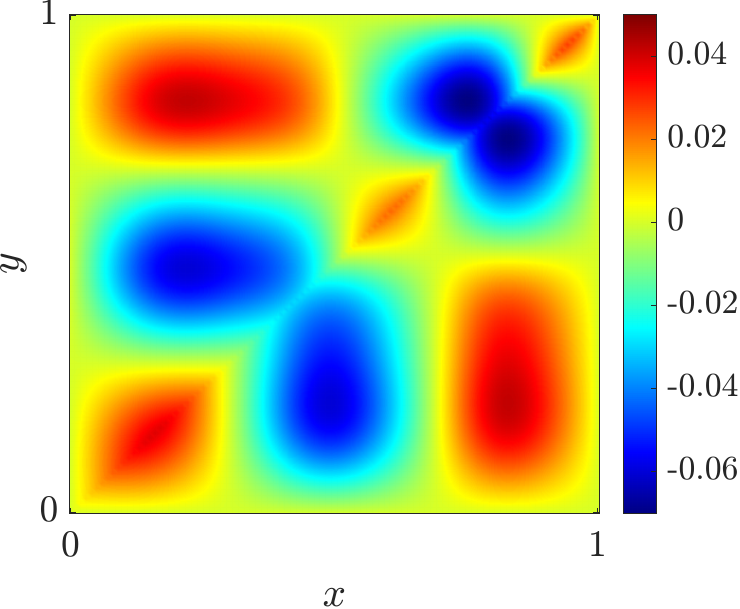}
			\caption{$\widecheck{G}(x,y)$}
			\label{fig-SgEncdGreen-1D-VarbHelmholtz-2D}
		\end{subfigure}
		\hspace{0.25cm}
		\begin{subfigure}{0.2894\textwidth}
			\centering
			\includegraphics[width=\textwidth]{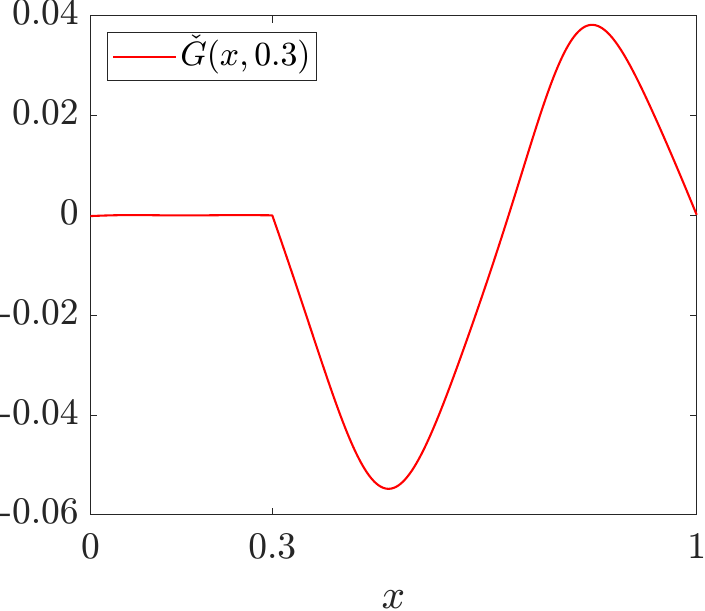}
			\caption{$\widecheck{G}(x,0.3)$}
			\label{fig-SgEncdGreen-1D-VarbHelmholtz-yfixed-2D}
		\end{subfigure}
		\hspace{0.5cm}
		\begin{subfigure}{0.288\textwidth}
			\centering
			\includegraphics[width=\textwidth]{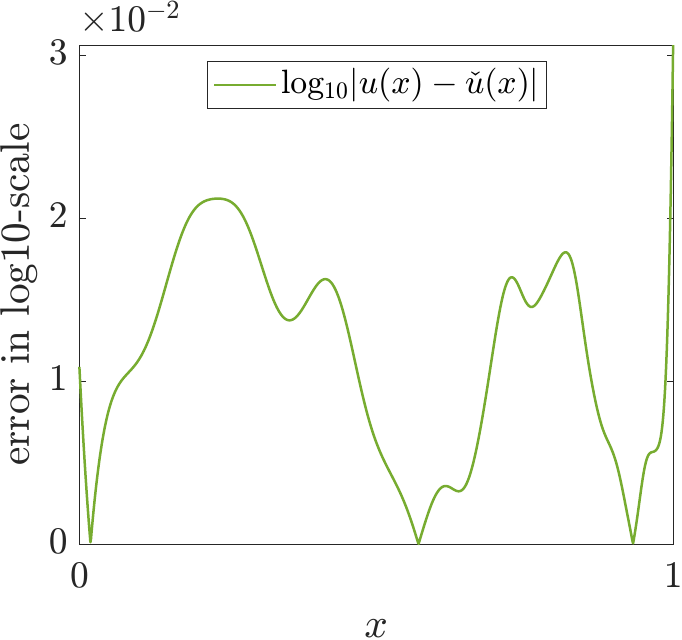}
			\caption{$\textnormal{log}_{10}|u(x) - \check{u}(x)|$}
			\label{fig-SgEncdGreen-1D-VarbHelmholtz-u-PtErr}	
		\end{subfigure}
		
		\vspace{-0.1cm}
		
		\caption{Numerical validation of our singularity-encoded Green's function for the Helmholtz equation \eqref{BVP-1D-Helmholtz}.}
		\label{Experiment-1D-VarbHelmholtz}
		
		\vspace{-0.3cm}
	\end{figure}
	%%--------------------------------------%
	
	Although initially developed for positive-definite problems \eqref{BVP-ExactSolu}, our singularity-encoded learning method remains effective when generalized to indefinite problems, as demonstrated through the Helmholtz equations
	\begin{equation}
		\begin{array}{cl}
			-(c(x)u'(x))' - k^2(x) u(x) = f(x), \ \ \ &\ \textnormal{for}\ x\in \Omega = (0,1),
		\end{array}
		\label{BVP-1D-Helmholtz}
	\end{equation}
	with Dirichlet boundary conditions $u(0) = u(1) = 0$, variable coefficients $c(x) = (x-2)^2$ and $k(x) = 15\sin(10x)$. Accordingly, the exact Green's function $G(x,y)$ satisfies, for any $y\in\Omega$,
	\begin{equation}
		\begin{array}{cl}
			-c(x)\partial_{xx}G(x,y) - c'(x)\partial_x G(x,y) - k^2(x)G(x,y) = 0, \ \ \ &\ \textnormal{for}\ x\in\Omega\setminus\Gamma, \\
			\llbracket G(x,y) \rrbracket = 0,\ \ \ c(x)\llbracket \partial_x G(x,y) \rrbracket = -1, \ \ \ &\ \textnormal{for}\ x\in\Gamma = \{ x \,|\, x=y\},\\
			G(x,y) = 0, \ \ \ &\ \textnormal{for}\ x\in\partial\Omega. \\
		\end{array}
		\label{BVP-1D-Helmholtz-ExactGreen}
	\end{equation}
	Unlike our illustrative example \eqref{BVP-1D-ToyModel-ExactGreen}, the jump of derivative in \eqref{BVP-1D-Helmholtz-ExactGreen} now incorporates a variable coefficient, which would complicate the derivation of the analytic Green's function \cite{stakgold2011green,evans2010partial}. Nevertheless, the exact Green's function maintains continuity across the interface $\Gamma$ while exhibits jump in its derivative, hence our prescribed augmented variable $\varphi(x,y)=|x-y|$ remains valid in capturing the singular behavior governed by equations \eqref{BVP-1D-Helmholtz-ExactGreen}.  
	
	Then, by \eqref{BVP-ExactGreen-Encoded} and \eqref{normalization-condition-embedded-1D}, the singularity-encoded Green's function $\widehat{G}(x,y,\varphi(x,y))$ for problem \eqref{BVP-1D-Helmholtz-ExactGreen} satisfies 	
		\begin{equation}
		\begin{array}{cl}
			\big[ c(x)\big(\partial_{xx}+2\textnormal{sgn}(x-y)\partial_{xz} + \partial_{zz} \big) + c'(x)\big(\partial_x + \textnormal{sgn}(x-y)\partial_z\big) + k^2(x) \big] \widehat{G}(x,y,\varphi(x,y)) = 0, & \mathrm{for}\ x \in \Omega\setminus\Gamma,\\
			\displaystyle -2c(x)\partial_{z}\widehat{G}(x,y,\varphi(x,y)) = 1, & \textnormal{for}\ x\in \Gamma,\\
			\widehat{G} (x, y, \varphi( x, y) ) = 0, & \mathrm{for}\ x \in \partial\Omega.
		\end{array}
		\label{BVP-1D-Helmholtz-SgEncdGreen}
	\end{equation}
	for any fixed $y\in(0,1)$, while the constraint $0 = \llbracket G(x,y) \rrbracket = \llbracket \widehat{G}(x,y,\varphi(x,y)) \rrbracket$ for $x\in\Gamma$ is automatically enforced owing to the continuity of our augmented variable $\varphi(x,y)=|x-y|$ along the interface $\Gamma = \{ x \,|\, x=y\}$.
	
	Next, the singularity-encoded Green's function is parametrized through a fully-connected neural network and learned via our proposed learning algorithm \autoref{Algorithm-SgEncdGreen}. The trained model $\widehat{G} (x, y, \varphi( x, y) )$ and its projection back onto the original 2-dimensional plane are shown in \autoref{fig-SgEncdGreen-1D-VarbHelmholtz-3D} and \autoref{fig-SgEncdGreen-1D-VarbHelmholtz-2D}, respectively. Moreover, by choosing $y=0.3$, our approximate Green's function $\widecheck{G} (x, y=0.3)$ manifests a discontinuity in its first-order derivative as displayed in \autoref{fig-SgEncdGreen-1D-VarbHelmholtz-yfixed-2D}, which can be effectively resolved due to its embedding into a smooth hyperplane $\widehat{G} (x, 0.3, z )$ within a one-order higher-dimensional space (see \autoref{fig-SgEncdGreen-1D-VarbHelmholtz-yfixed-3D}).
	
	As 	the close-form Green's function is generally unavailable for problems with variable coefficients, we adopt the numerical integration \eqref{BVP-NumSolu-Representation-SgEncdGreen} to evaluate the precision of our singularity-encoded Green's function. Specifically, with the forcing term $f(x)$ derived from the exact solution  
	\begin{equation*}
		u(x)=10x-10x^2+0.5 \sin(20\pi x^3)
	\end{equation*}		
	of Helmholtz equations \eqref{BVP-1D-Helmholtz}, we use $h=2^{-10}$ in \eqref{BVP-NumSolu-Representation-SgEncdGreen} and compare our reconstructed solution $\check{u}(x)$ with the exact solution $u(x)$. \autoref{fig-SgEncdGreen-1D-VarbHelmholtz-u} and \autoref{fig-SgEncdGreen-1D-VarbHelmholtz-u-PtErr} demonstrate good agreement between numerical and true solutions, which implicitly implies that our approximate Green's function achieves satisfactory accuracy.
	
%%%%%%%%%%%%%%%%%%%%%%%%%%%%%%%%%%%%%%%%%%%%%%%%%%%%%
		
%%%%%%%%%%%%%%%%%%%%%%%%%%%%%%%%%%%%%%%%%%%%%%%%%%%%%
\subsubsection{Neural Preconditioning Matrix}	
	
	%-------------------------%
	\begin{table}[b!]
		\centering
		\renewcommand\arraystretch{2} 
		\begin{tabular}{|c|c|c|c|c|}
			\hline
			mesh size & eigenvalues of preconditioned matrix & condition number & $\#$ iterations & $\ell_2$-norm error \\
			\hline
			$h = 2^{-6}$ &
			\begin{minipage}[b]{0.3\columnwidth}
				\centering
				\raisebox{-.4\height}{\includegraphics[width=\linewidth]{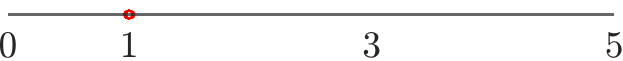}}			
			\end{minipage}
			& $\kappa (\widecheck{B}A) = 1.03$ &\makecell[c]{$4$\\$(30)$} &\makecell[c]{$8.15 \times 10^{-10}$\\$(4.92 \times 10^{-1})$} \\
			\hline
			$h = 2^{-8}$ &
			\begin{minipage}[b]{0.3\columnwidth}
				\centering
				\raisebox{-.4\height}{\includegraphics[width=\linewidth]{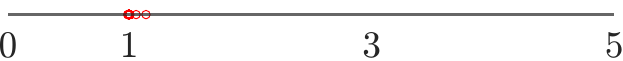}}
			\end{minipage}
			& $\kappa (\widecheck{B}A) = 1.15$  &\makecell[c]{$5$\\$(252)$} &\makecell[c]{$3.85 \times 10^{-10}$\\$(6.52 \times 10^{-3})$}\\
			\hline
			$h = 2^{-10}$ &
			\begin{minipage}[b]{0.3\columnwidth}
				\centering
				\raisebox{-.4\height}{\includegraphics[width=\linewidth]{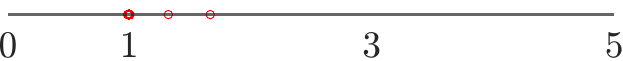}}
			\end{minipage}
			& $\kappa (\widecheck{B}A) = 1.69$  &\makecell[c]{$5$\\$(1014)$} &\makecell[c]{$6.82 \times 10^{-9}$\\$(2.24 \times 10^{-4})$}\\
			\hline
			$h = 2^{-12}$ &
			\begin{minipage}[b]{0.3\columnwidth}
				\centering
				\raisebox{-.4\height}{\includegraphics[width=\linewidth]{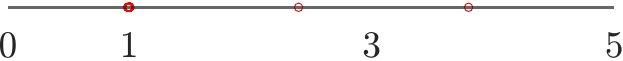}}
			\end{minipage}
			& $\kappa (\widecheck{B}A) = 3.85$  &\makecell[c]{$6$\\$(4058)$} &\makecell[c]{$1.55 \times 10^{-10}$\\$(1.98 \times 10^{-4})$}\\
			\hline
		\end{tabular}
		\caption{Comparative study of GMRES with and without neural preconditioning for solving the discrete linear system of the Helmholtz equation \eqref{BVP-1D-Helmholtz}.}
		\label{GreenPreconditioner-Helmholtz1DExp1}
		\vspace{-0.3cm}
	\end{table}
	%-------------------------%	
	
	The Helmholtz equations \eqref{BVP-1D-Helmholtz} are discretized using a central difference scheme \cite{leveque2007finite}, leading to the linear system $AU=F$ (not relabelled for simplicity). Unlike the previous example \eqref{BVP-1D-ToyExmp-MatrixForm}, the discrete system becomes indefinite that undermines the effectiveness of traditional iterative methods (e.g., the damped Jacobi method). To guarantee convergence for indefinite systems, the generalized minimal residual method (abbreviated as GMRES) \cite{saad1986gmres,golub2013matrix} is typically utilized with certain restart strategies. Nevertheless, ensuring efficient convergence remains challenging \cite{saad1986gmres}, necessitating the development and deployment of reliable preconditioners.	
	
	Building on the singularity-encoded Green's function offline-trained in our previous section, a preconditioner $\widecheck{B}$ can be constructed through Algorithm \ref{Alg-SgEncd-Green-Preconditioner} for solving our indefinite systems. To be specific, the spectral property of our preconditioned matrix $\widecheck{B}A$ is computationally reported in \autoref{GreenPreconditioner-Helmholtz1DExp1}, including the distribution of eigenvalues and the condition number across different mesh sizes. Our neural preconditioner, built from a surrogate model of its exact Green's function, produces clustered eigenvalues for the preconditioned system, which is crucial for the convergence of iterative methods. Subsequently, we evaluate the performance of GMRES (without restarting) for both the original (displayed in parentheses) and preconditioned systems, with comparative results reported in the final two columns of \autoref{GreenPreconditioner-Helmholtz1DExp1}. Clearly, our preconditioned system can achieve satisfactory accuracy within only a few iterations, while the original system exhibits persistent error even after many iterations.

%%%%%%%%%%%%%%%%%%%%%%%%%%%%%%%%%%%%%%%%%%%%%%%%%%%%%

%%%%%%%%%%%%%%%%%%%%%%%%%%%%%%%%%%%%%%%%%%%%%%%%%%%%%
\subsubsection{Hybrid Iterative Methods}	
	
\begin{figure}[t!]
		\centering
		\begin{subfigure}{0.47\textwidth}
			\centering
			\includegraphics[width=1\textwidth]{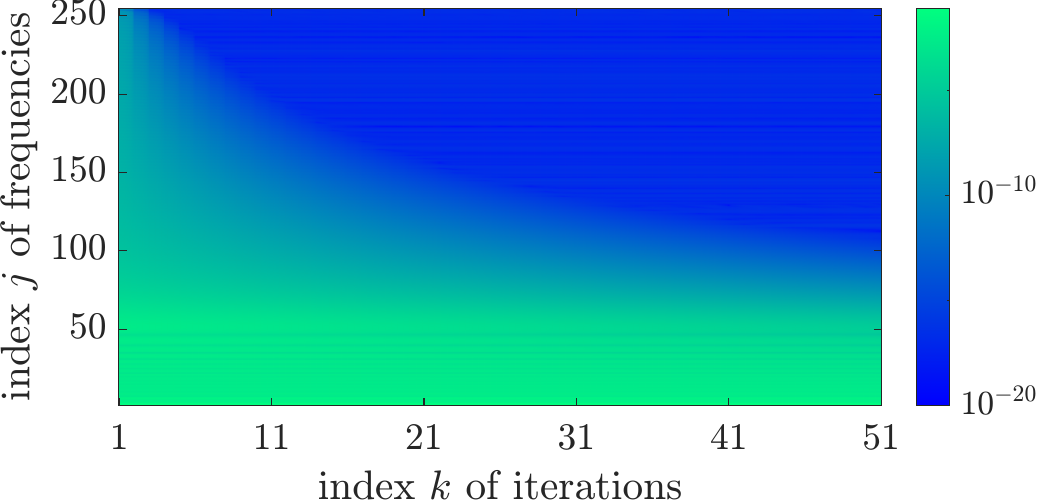}
			\caption{the spectral error of damped Jacobi method \eqref{1D-Poisson-Damped-Jacobi}}
			\label{fig-1D-VarbHelmholtz-Jacobi-ModeErr}
		\end{subfigure}
		\hspace{0.2cm}		
		\begin{subfigure}{0.47\textwidth}
			\centering
			\includegraphics[width=1\textwidth]{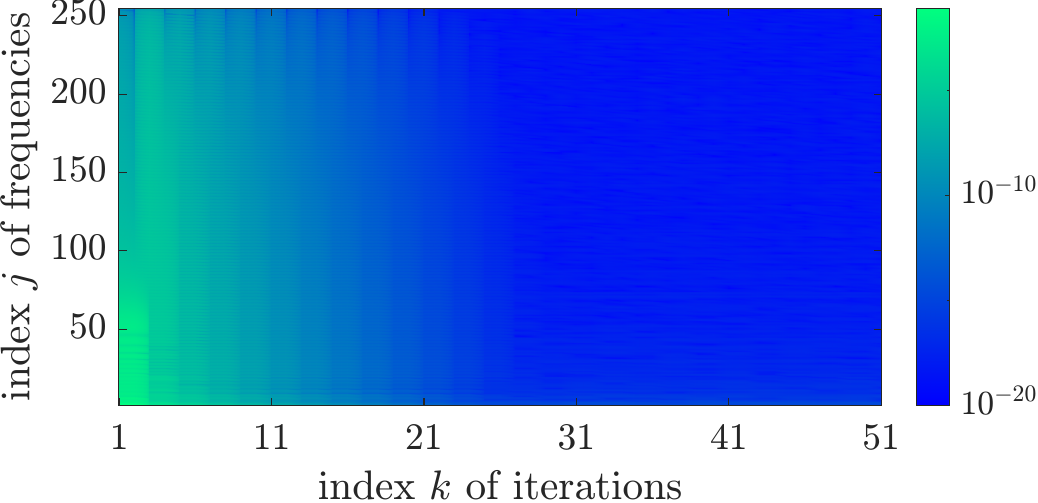}
			\caption{the spectral error of hybrid iterative method \eqref{1D-Poisson-Hybrid-IteM}}
			\label{fig-1D-VarbHelmholtz-Hybrid-ModeErr}
		\end{subfigure}
		
		\centering
		\vspace{0.3cm}
		\begin{subfigure}{0.31\textwidth}
			\centering
			\includegraphics[width=1\textwidth]{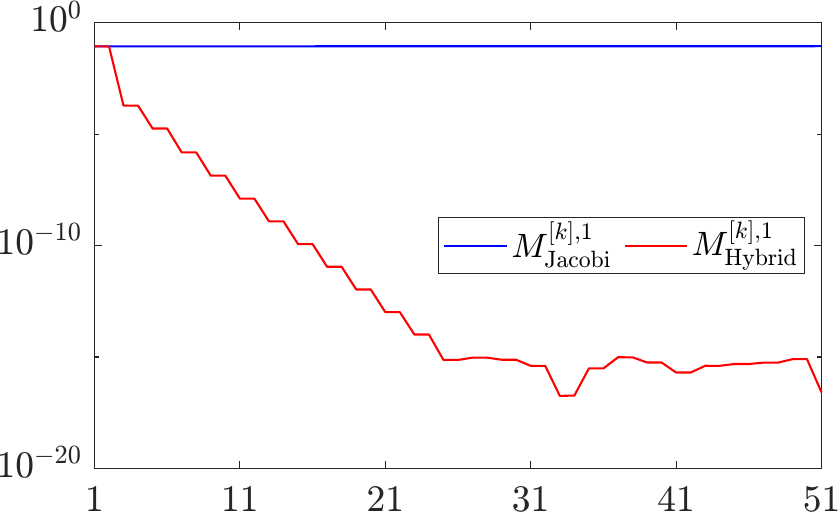}
			\caption{the lowest-frequency error component}
			\label{fig-1D-VarbHelmholtz-Hybrid-ModeErr-j1} 
		\end{subfigure}
		\hspace{0.1cm}
		\begin{subfigure}{0.31\textwidth}
			\centering
			\includegraphics[width=1\textwidth]{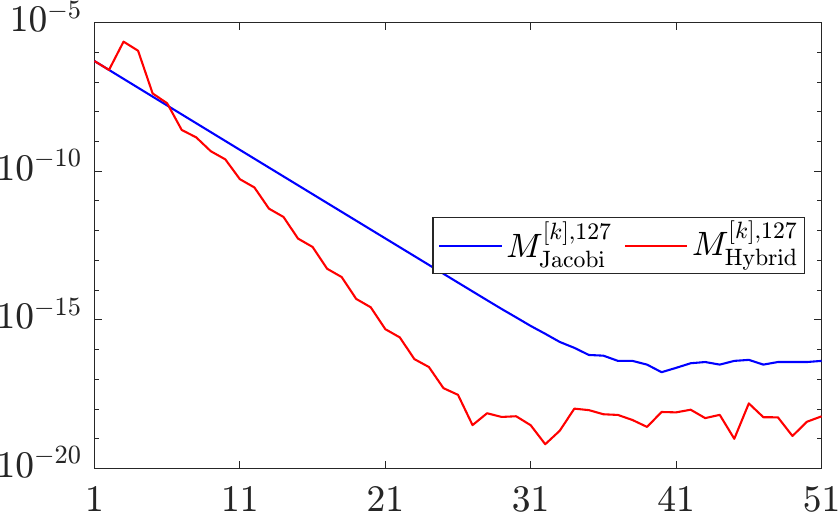}
			\caption{the mid-frequency error component}
			\label{fig-1D-VarbHelmholtz-Hybrid-ModeErr-j127} 
		\end{subfigure}
		\hspace{0.1cm}
		\begin{subfigure}{0.31\textwidth}
			\centering
			\includegraphics[width=1\textwidth]{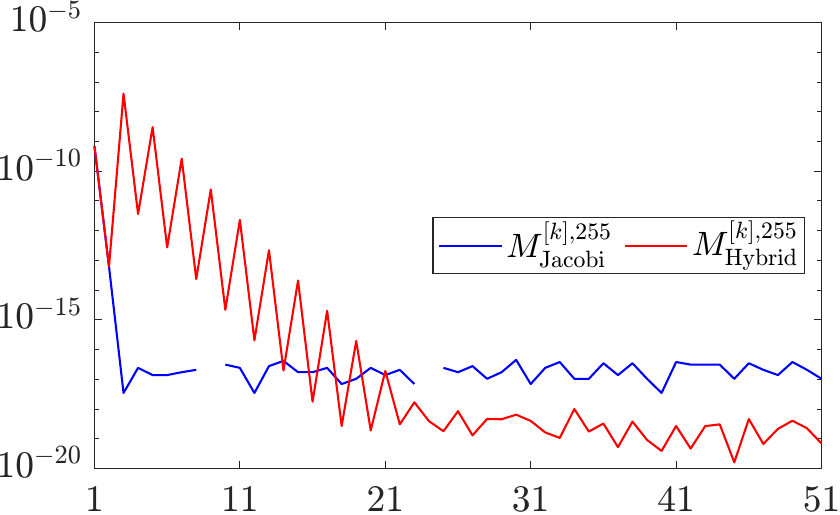}
			\caption{the highest-frequency error component}
			\label{fig-1D-VarbHelmholtz-Hybrid-ModeErr-j255} 
		\end{subfigure}
		
		\caption{Comparison of spectral errors for solving the discrete linear system of the Helmholtz equation \eqref{BVP-1D-Helmholtz}.}
		\label{Experiments-VarbHelmholtz1D-IteM-1}
		\vspace{-0.3cm}
	\end{figure}

	Note that the indefiniteness of problem \eqref{BVP-1D-Helmholtz} causes many classical schemes, such as the Jacobi and Gauss-Seidel methods \cite{golub2013matrix}, to diverge when applied to solving its discrete linear system. Even for convergent components, the slow reduction of low-frequency errors undermines the overall efficiency. Here, we consider a linear system with mesh size $h=2^{-8}$ for illustration, adopting the damped Jacobi method \eqref{1D-Poisson-Damped-Jacobi} with a relaxation parameter $\omega=0.5$, and depict in \autoref{fig-1D-VarbHelmholtz-Jacobi-ModeErr} the mode-wise error across all frequencies.
	
	%--------------- 1DVarbHelmholtz-Example iterative methods ---------------%
	%--------------- 1DVarbHelmholtz-Example iterative methods ---------------%

	%--------------- 1DVarbHelmholtz-Example iterative methods ---------------%
	\begin{figure}[b!]		
		\centering
		\begin{subfigure}{0.46\textwidth}
			\centering
			\includegraphics[width=1\textwidth]{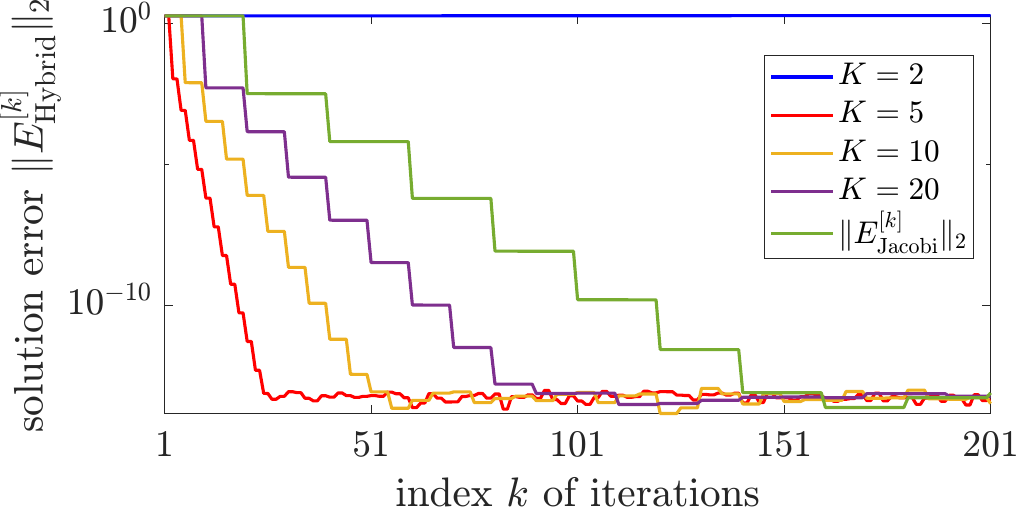}
			\caption{comparison of the solution error over iterations}
			\label{fig-1D-VarbHelmholtz-Solu-Err} 
		\end{subfigure}
		\hspace{0.3cm}
		\begin{subfigure}{0.46\textwidth}
			\centering
			\includegraphics[width=1\textwidth]{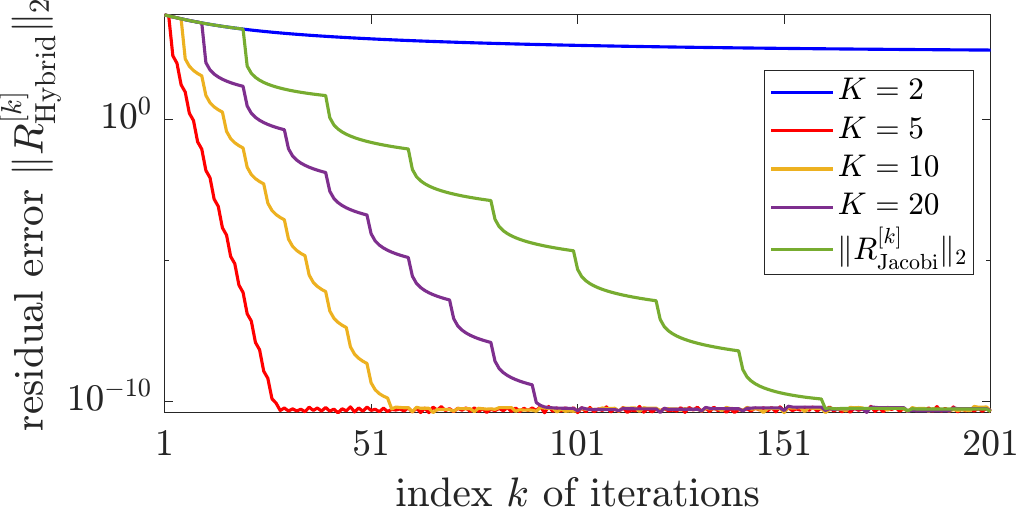}
			\caption{comparison of the residual error over iterations}
			\label{fig-1D-VarbHelmholtz-Res-Err}
		\end{subfigure}
		
		\caption{Hybrid iteration \eqref{1D-Poisson-Hybrid-IteM} with different switching periods for the discrete system of the Helmholtz equation \eqref{BVP-1D-Helmholtz}.}
		\label{Experiments-VarbHelmholtz1D-IteM-2}
		\vspace{-0.3cm}
	\end{figure}
	%--------------- 1DVarbHelmholtz-Example iterative methods ---------------%

	As can be seen from \autoref{fig-1D-VarbHelmholtz-Hybrid-ModeErr-j1}, the mode-wise error $M^{[k],1}$ associated with the damped Jacobi method fails to converge due to the indefiniteness. Nevertheless, as shown in \autoref{fig-1D-VarbHelmholtz-Hybrid-ModeErr-j255}, the rapid convergence of high-frequency errors remains unaffected. In contrast, our hybrid iterative method \eqref{1D-Poisson-Hybrid-IteM} successfully resolves these limitations, as demonstrated in \autoref{fig-1D-VarbHelmholtz-Hybrid-ModeErr} and \autoref{fig-1D-VarbHelmholtz-Hybrid-ModeErr-j1}. While the spectral bias inherent in our neural preconditioner introduces extra high-frequency errors (see \autoref{fig-1D-VarbHelmholtz-Hybrid-ModeErr-j255}), these artifacts can be efficiently eliminated via traditional smoothing iterations. Moreover, empirical results in \autoref{fig-1D-VarbHelmholtz-Solu-Err} and \autoref{fig-1D-VarbHelmholtz-Res-Err} imply that the switching period $K=2$ achieves the minimum number of iterations required to attain convergence, whereas larger values impede the efficiency.
	
\begin{figure}[!t]
		\centering
		\begin{subfigure}{0.5\textwidth}
			\centering
			\includegraphics[width=\textwidth]{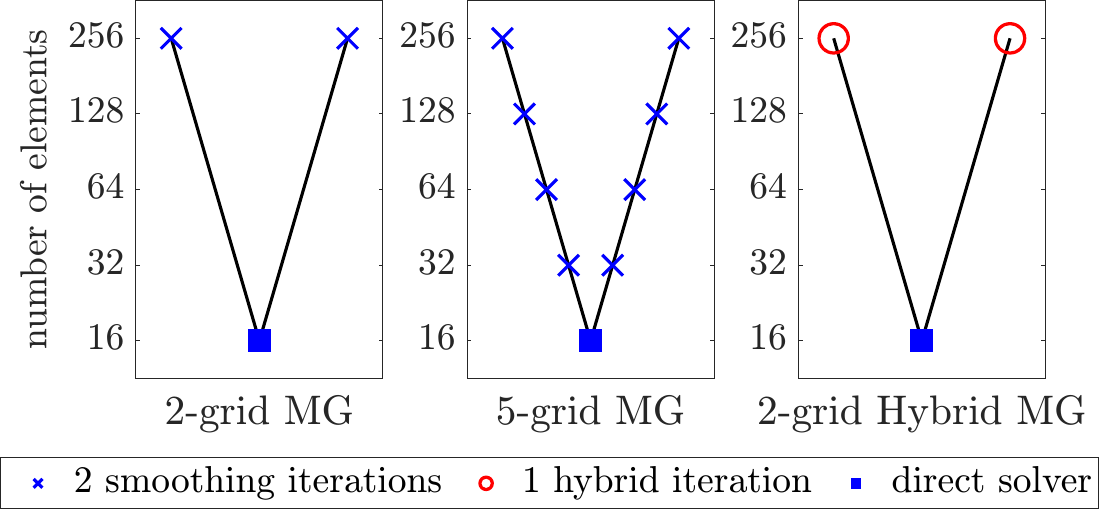}
			\caption{V-cycle structures for different approaches}
			\label{fig-1D-Helmholtz-MG-V-Cycle}
		\end{subfigure}
		\hspace{0.25cm}
		\begin{subfigure}{0.45\textwidth}
			\centering
			\includegraphics[width=\textwidth]{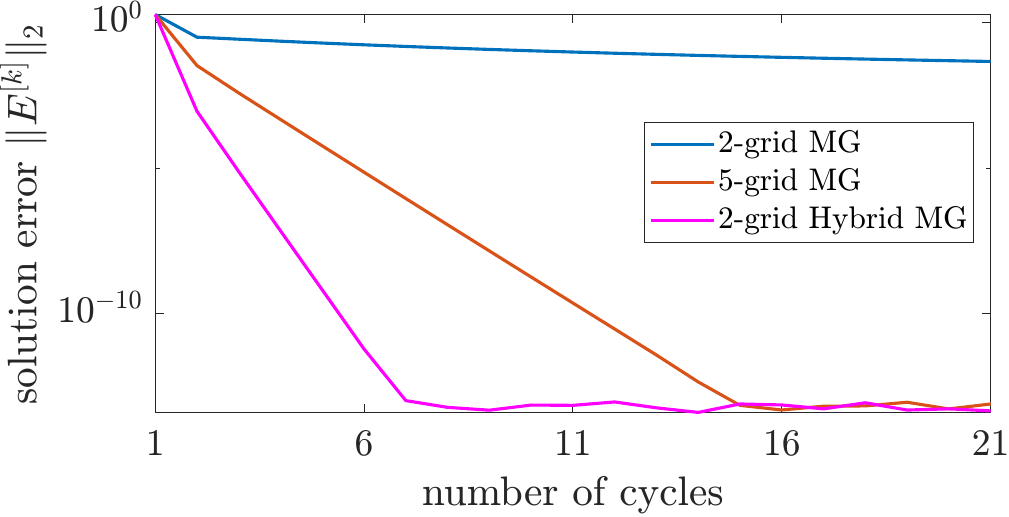}
			\caption{solution errors over iterative cycles}
			\label{fig-1D-Helmholtz-MG-Solu-Err}
		\end{subfigure}
		
		\caption{Comparison of multigrid methods with and without hybrid iterations for solving the discrete system of the Helmholtz equation \eqref{BVP-1D-Helmholtz}.}
		\label{Experiment-Helmholtz1D-HybridMG}
		\vspace{-0.3cm}
	\end{figure}

	Then, we incorporate the hybrid iterative approach \eqref{1D-Poisson-Hybrid-IteM} within the multigrid framework, employing the same strategies as applied to the Poisson's equations in section \ref{Section-1D-ToyModel-Hybrid-IteM} (see also \autoref{fig-1D-Helmholtz-MG-V-Cycle} for further details). Notably, applying our hybrid iterative scheme to both pre- and post-smoothing stages enables convergence to satisfactory accuracy in 7 cycles, even with the use of a two-grid structure as depicted in \autoref{fig-1D-Helmholtz-MG-Solu-Err}.
	
	%---------------------------%
	%---------------------------%
	
%%%%%%%%%%%%%%%%%%%%%%%%%%%%%%%%%%%%%%%%%%%%%%%%%%%%%

	%--------------------------------------%
	% mannully displayed here, content belonging to the next section
	\begin{figure}[!b]
		\centering
		\begin{subfigure}{0.31\textwidth}
			\centering
			\includegraphics[width=\textwidth]{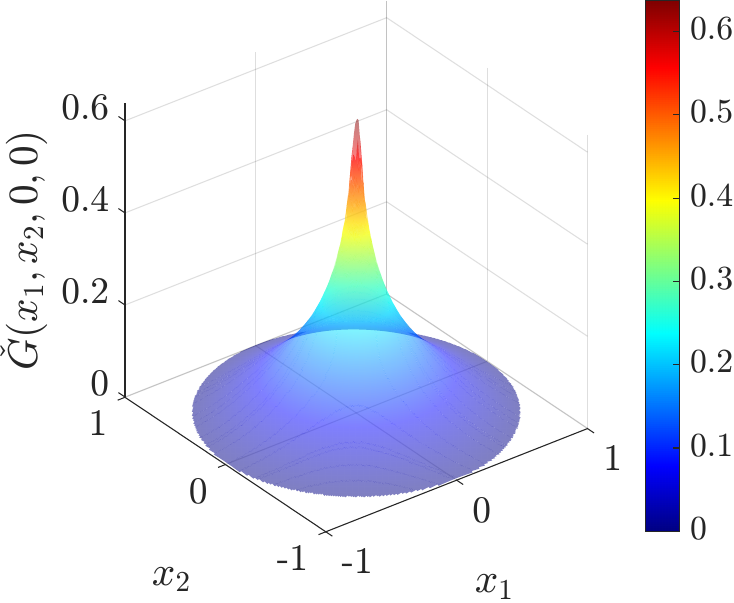}
			\caption{$\widecheck{G}(x_1,x_2,0,0)$}
			\label{fig-SgEncdGreen-2D-Poisson-surf-x0}
		\end{subfigure}
		\hspace{0.2cm}
		\begin{subfigure}{0.294\textwidth}
			\centering
			\includegraphics[width=\textwidth]{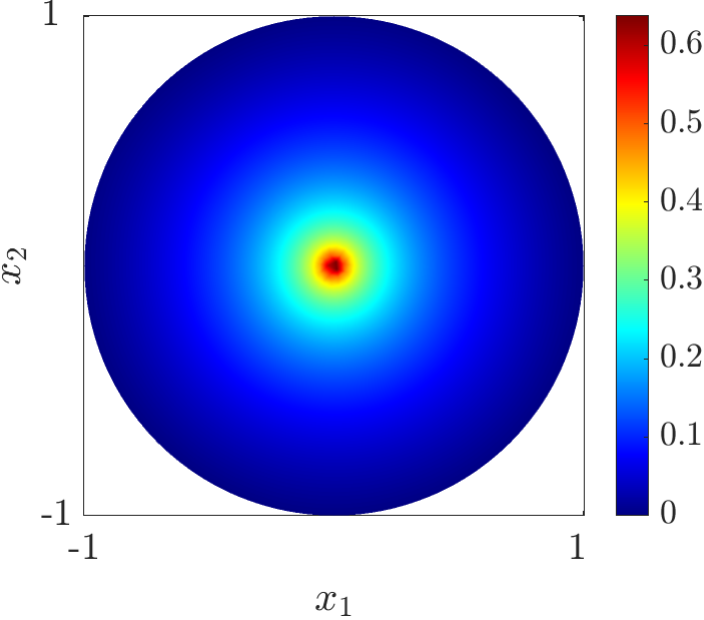}
			\caption{$\widecheck{G}(x_1,x_2,0,0)$}
			\label{fig-SgEncdGreen-2D-Poisson-x0}
		\end{subfigure}
		\hspace{0.2cm}
		\begin{subfigure}{0.31\textwidth}
			\centering
			\includegraphics[width=\textwidth]{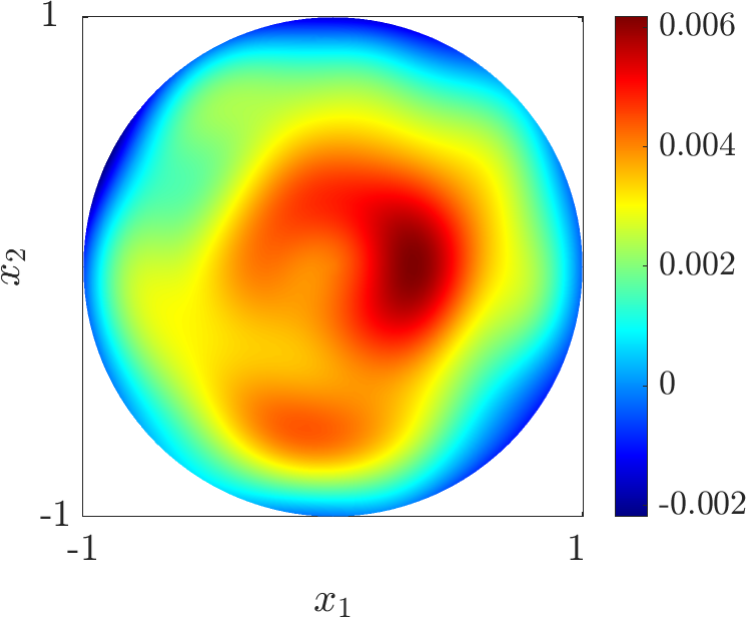}
			\caption{$G(x_1,x_2,0,0) - \widecheck{G}(x_1,x_2,0,0)$}
			\label{fig-Green-2D-Poisson-PtErr-x0}
		\end{subfigure}
	
		%\vspace{-0.1cm}
	
		\begin{subfigure}{0.31\textwidth}
			\centering
			\includegraphics[width=\textwidth]{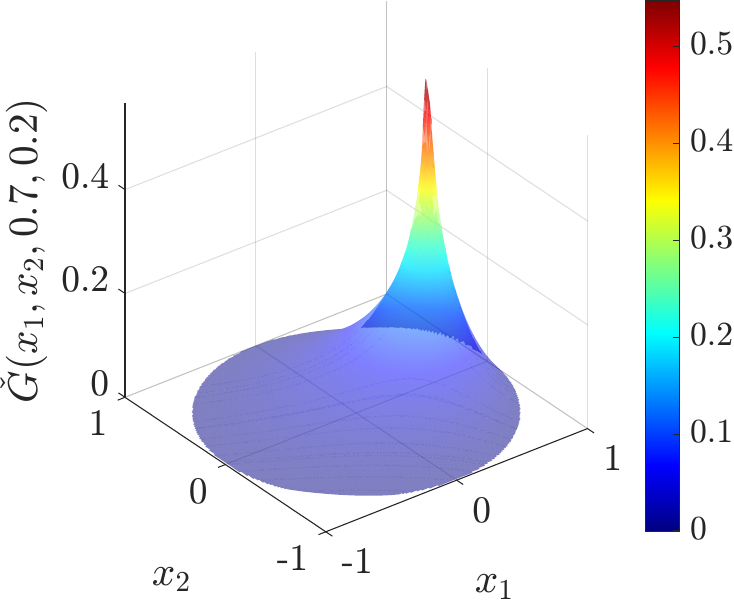}
			\caption{$\widecheck{G}(x_1,x_2,0.7,0.2)$}
			\label{fig-SgEncdGreen-2D-Poisson-surf-x1}
		\end{subfigure}
		\hspace{0.3cm}
		\begin{subfigure}{0.294\textwidth}
			\centering
			\includegraphics[width=\textwidth]{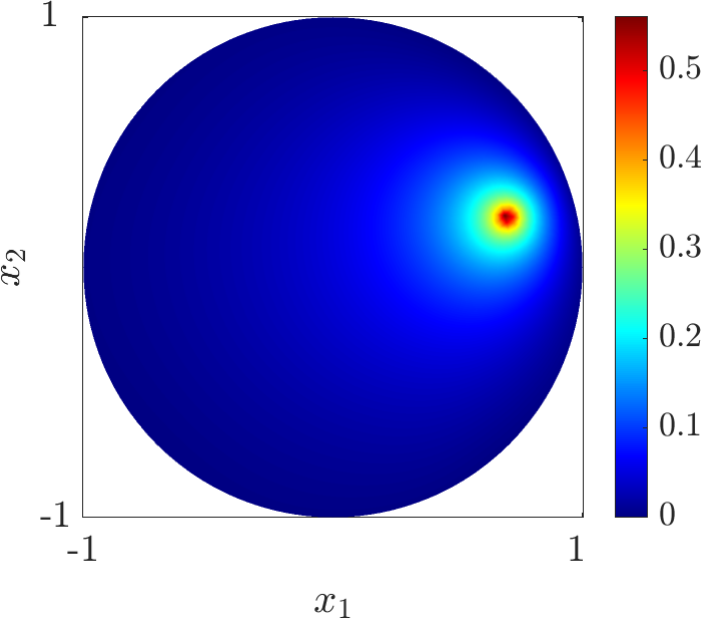}
			\caption{$\widecheck{G}(x_1,x_2,0.7,0.2)$}
			\label{fig-SgEncdGreen-2D-Poisson-x1}
		\end{subfigure}
		\hspace{0.2cm}
		\begin{subfigure}{0.31\textwidth}
			\centering
			\includegraphics[width=\textwidth]{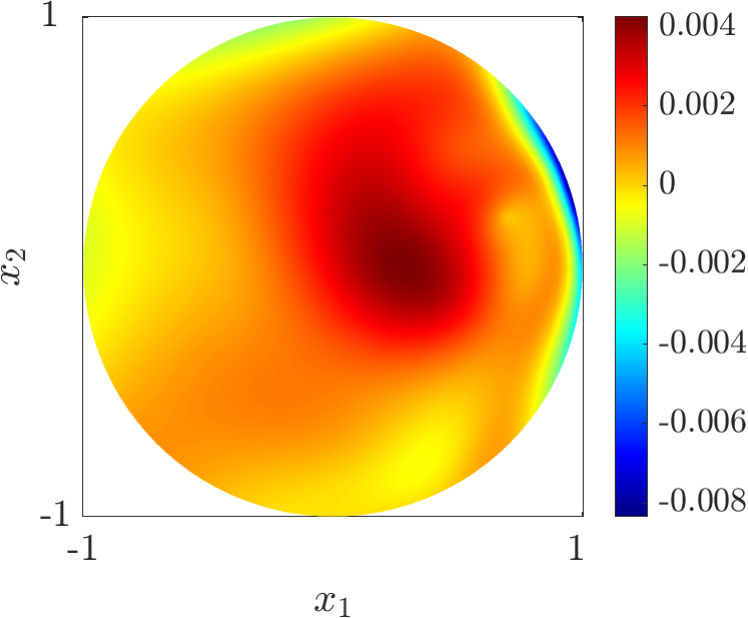}
			\caption{$G(x_1,x_2,0.7,0.2) - \widecheck{G}(x_1,x_2,0.7,0.2)$}
			\label{fig-Green-2D-Poisson-PtErr-x1}
		\end{subfigure}
		
		\vspace{-0.1cm}
	
		\caption{Numerical realizations of our singularity-encoded Green's function for the two-dimensional problem \eqref{BVP-2D-Laplacian}.}
		\label{Experiment-2D-Poisson-1}
		
		\vspace{-0.3cm}
	\end{figure}
	%--------------------------------------%
	
%%%%%%%%%%%%%%%%%%%%%%%%%%%%%%%%%%%%%%%%%%%%%%%%%%%%%
\subsection{Green's Function for Poisson's Equation in Two Dimension}
	
	We now extend our empirical studies to boundary value problems \eqref{BVP-ExactSolu} in higher dimensions, in which the Green's function exhibits fundamentally different singular behavior \eqref{Singularity-ExactGreen} compared to the one-dimensional cases. Here, we consider the benchmark Poisson's equation on a two-dimensional unit disc, namely,
	\begin{equation}
		-\Delta u(\bm{x}) = f(\bm{x}), \ \ \ \ \textnormal{for}\ \bm{x} \in \Omega = \{\bm{x} : \lVert \bm{x} \rVert < 1\}, 
		\label{BVP-2D-Laplacian}
	\end{equation}
	where $\bm{x} = (x_1,x_2)\in\mathbb{R}^2$ and $u(\bm{x})=0$ for $\bm{x}\in\partial\Omega$. Accordingly, the Green's function for problem \eqref{BVP-2D-Laplacian} satisfies 
	\begin{equation}
		\begin{array}{cl}
			-\Delta G(\bm{x},\bm{y}) = \delta(\bm{x} -\bm{y}), \ \ \ &\ \textnormal{for}\ \bm{x}\in\Omega,\\
			G(\bm{x},\bm{y}) = 0, \ \ \ &\ \textnormal{for}\ \bm{x}\in\partial\Omega,
		\end{array}
		\label{BVP-2D-Laplacian-ExactGreen}
	\end{equation}
	for any $\bm{y}=(y_1,y_2)\in\Omega$. By defining the point $\bm{y}^* = \frac{\bm{y}}{\lVert \bm{y} \rVert^2}$ dual to $\bm{y}$, the Green's function can be expressed as \cite{myint2007linear}  
	\begin{equation}
		G(\bm{x},\bm{y}) =  -\frac{1}{2\pi} \ln \left(  \frac{\lVert \bm{x} - \bm{y} \rVert}{\lVert \bm{x} - \bm{y}^* \rVert} \cdot \frac{1}{\lVert \bm{y} \rVert} \right) = - \frac{1}{4\pi}\ln \left( \dfrac{(x_1-y_1)^2+(x_2-y_2)^2}{(x_1y_2-x_2y_1)^2+(x_1y_1+x_2y_2-1)^2} \right),
		\label{ExactGreen-2D-Laplacian}
	\end{equation}
	which is highly expensive for traditional mesh-based methods to approximate due to the curse of dimensionality. Furthermore, the Green's function \eqref{ExactGreen-2D-Laplacian} exhibits a logarithmic blow-up, i.e., $G(\bm{x},\bm{y}) \sim \ln\lVert \bm{x}-\bm{y}\rVert$ as $\bm{x}\to\bm{y}$, which introduces additional challenges for classical mesh-based numerical methods.
	
	%--------------------------------------%
	\begin{figure}[!b]
		
		\begin{subfigure}{0.305\textwidth}
			\centering
			\includegraphics[width=\textwidth]{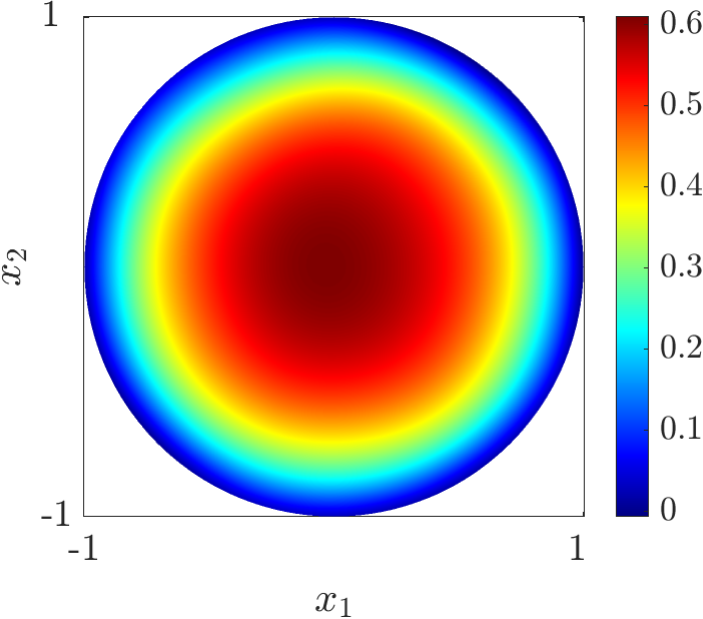}
			\caption{$\check{u}(x_1,x_2)$}
			\label{fig-SgEncdGreen-2D-solution-u}
		\end{subfigure}
		\hspace{0.2cm}
		\begin{subfigure}{0.30\textwidth}
			\centering
			\includegraphics[width=\textwidth]{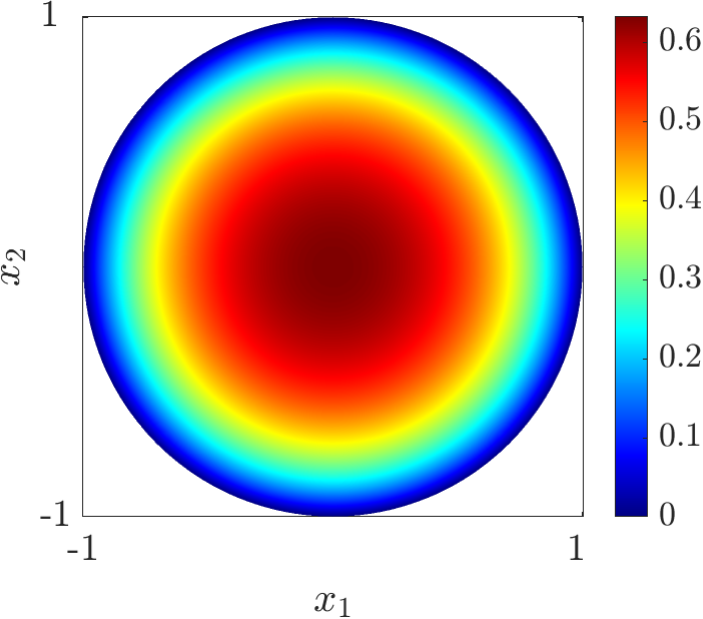}
			\caption{$u(x_1,x_2)$}
			\label{fig-2D-exact-solution-u}
		\end{subfigure}
		\hspace{0.2cm}
		\begin{subfigure}{0.308\textwidth}
			\centering
			\includegraphics[width=\textwidth]{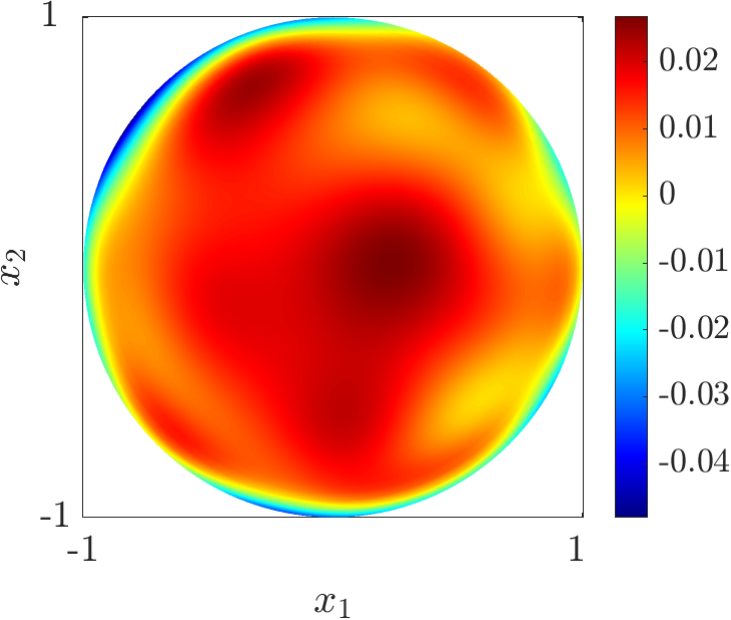}
			\caption{$u(x_1,x_2) - \check{u}(x_1,x_2)$}
			\label{fig-2D-PtErr-u}	
		\end{subfigure}
	
		%\vspace{-0.1cm}
	
		\begin{subfigure}{0.31\textwidth}
			\centering
			\includegraphics[width=\textwidth]{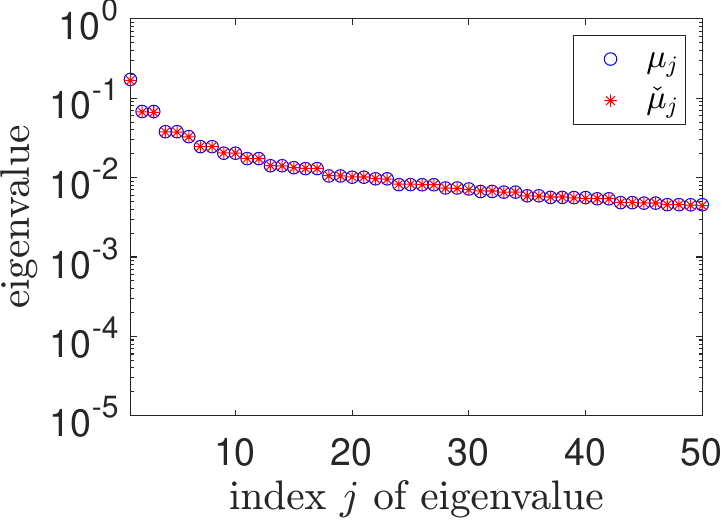}
			\caption{$\mu_j$, $\check{\mu}_j$ with $1\leq j\leq 50$}
			\label{fig-2D-eigenvalues-1to50}
		\end{subfigure}
		\hspace{0.1cm}
		\begin{subfigure}{0.31\textwidth}
			\centering
			\includegraphics[width=\textwidth]{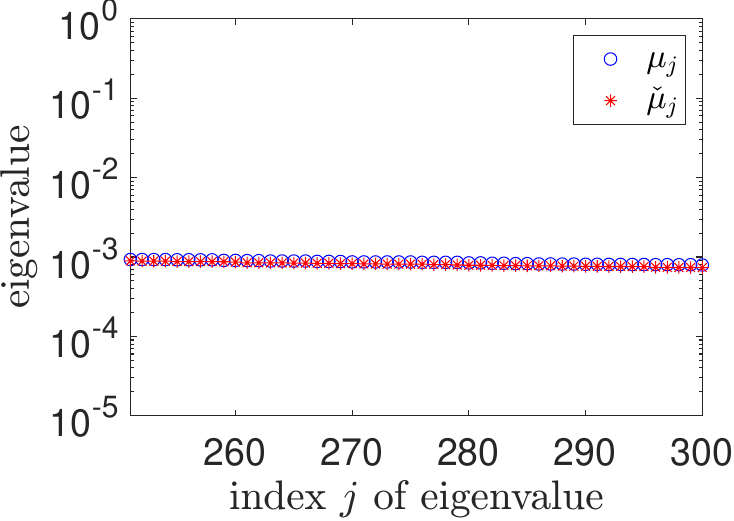}
			\caption{$\mu_j$, $\check{\mu}_j$ with $251\leq j\leq 300$}
			\label{fig-2D-eigenvalues-51to100}
		\end{subfigure}
		\hspace{0.1cm}
		\begin{subfigure}{0.31\textwidth}
			\centering
			\includegraphics[width=\textwidth]{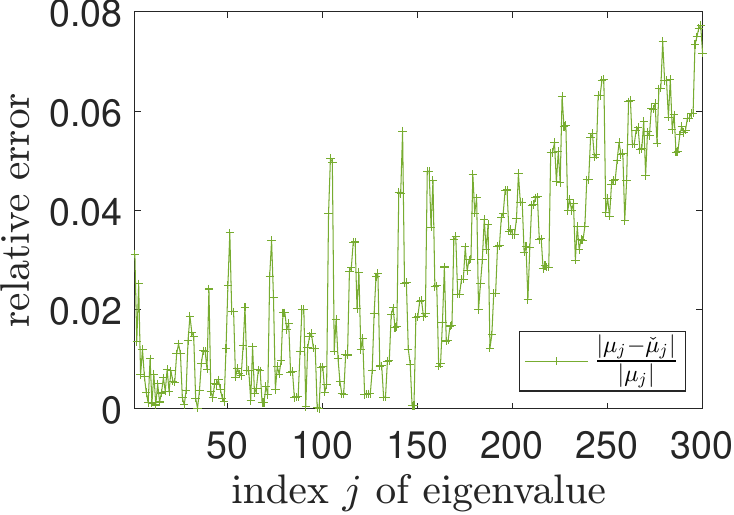}
			\caption{$\delta_{\mu_j}$ with $1\leq j\leq 300$}
			\label{fig-2D-eigenvalues-relerr}	
		\end{subfigure}	
		
		\vspace{-0.1cm}
	
		\caption{Numerical validation of our singularity-encoded Green's function for the two-dimensional problem \eqref{BVP-2D-Laplacian}.}
		\label{Experiment-2D-Poisson}
		
		\vspace{-0.3cm}
	\end{figure}
	%--------------------------------------%
	
	To address the logarithmic singularity, we augment the input variables of our deep surrogate model $\widehat{G}(\bm{x},\bm{y},\varphi(\bm{x},\bm{y}))$ with an extra feature $\varphi(\bm{x},\bm{y})=\ln\lVert \bm{x}-\bm{y}\rVert$ as discussed in section \ref{Sec-Green-Higher-Dimension}. Subsequently, by integrating the $\delta$-function over some ball $B_\epsilon(\bm{y})$, formula \eqref{BVP-ExactGreen-Encoded} implies that our singularity-encoded Green's function satisfies, for any $\bm{y}\in\Omega$,
	\begingroup
	\renewcommand*{\arraystretch}{2.2}
	\begin{equation}
		\begin{array}{cl}
			-\Delta_{\bm{x}}\widehat{G}(\bm{x},\bm{y},\varphi(\bm{x},\bm{y})) - 2 \nabla_{\bm{x}} \big( \partial_z\widehat{G}(\bm{x},\bm{y},\varphi(\bm{x},\bm{y})) \big) \cdot \dfrac{\bm{x}-\bm{y}}{\lVert \bm{x} - \bm{y} \rVert ^2}  - \dfrac{\partial_{zz}\widehat{G}(\bm{x},\bm{y},\varphi(\bm{x},\bm{y}))}{\lVert \bm{x} - \bm{y} \rVert ^2} = 0, \ \ \ &\ \textnormal{for}\ \bm{x}\in\Omega\setminus\Gamma, \\
			\displaystyle -\int_{\partial B_\epsilon (\bm{y})} \bigg(\nabla_{\bm{x}}G(\bm{x},\bm{y},\varphi(\bm{x},\bm{y})) \cdot \dfrac{\bm{x}-\bm{y}}{\lVert \bm{x} - \bm{y} \rVert} + \dfrac{\partial_z\widehat{G}(\bm{x},\bm{y},\varphi(\bm{x},\bm{y}))}{\lVert \bm{x} - \bm{y} \rVert}\bigg) d S(\bm{x}) = 1, \ \ \ &\ \textnormal{for} \ \bm{x}\in\Gamma,\\
			G(\bm{x},\bm{y},\varphi(\bm{x},\bm{y})) = 0, \ \ \ &\ \textnormal{for}\ \bm{x}\in\partial\Omega,
		\end{array}
		\label{BVP-2D-Laplacian-SgEncdGreen}
	\end{equation}
	\endgroup
	where $\Gamma \!=\! \{ \bm{x} \,|\, \bm{x}=\bm{y}\}$. Notably, as discussed in \eqref{SgEncd-Green-1st-Derivatives}, our singularity-encoded Green's function \eqref{BVP-2D-Laplacian-SgEncdGreen} now exhibits bounded partial derivatives $\nabla_{\bm{x}}\widehat{G}(\bm{x},\bm{y},\varphi(\bm{x},\bm{y}))$ and $\partial_z\widehat{G}(\bm{x},\bm{y},\varphi(\bm{x},\bm{y}))$, which differs from existing approaches (see section \ref{Sec-RelatedWork-GreenFunc}) but requires computation within a five-dimensional space. Fortunately, artificial neural networks \cite{goodfellow2016deep} have demonstrated remarkable effectiveness in handling high-dimensional smooth functions.			
		
	Consequently, we employ a fully-connected neural network with Tanh activation functions to parametrize the singularity-encoded Green's function, which is then optimized via Algorithm \ref{Algorithm-SgEncdGreen}. To capture the singular behavior governed by our normalization condition \eqref{BVP-2D-Laplacian-SgEncdGreen}, we adopt four concentric circles $B_\epsilon(\bm{y}_m)$ for each sampled $\bm{y}_m\in\Gamma$, with the radius chosen to be $\epsilon = 0.1$, $0.08$, $0.064$ and $0.005$. To visualize our numerical results, two slices of the trained neural network are depicted in \autoref{fig-SgEncdGreen-2D-Poisson-surf-x0} and \autoref{fig-SgEncdGreen-2D-Poisson-surf-x1} by fixing source points at $\bm{y}=(0,0)$ and $(0.7,0.2)$ respectively. In addition, \autoref{fig-Green-2D-Poisson-PtErr-x0} and \autoref{fig-Green-2D-Poisson-PtErr-x1} depict the discrepancies between the true Green's function \eqref{ExactGreen-2D-Laplacian} and our surrogate model for these two cases, demonstrating satisfactory resolution of the logarithmic singularity near $\bm{x}=\bm{y}$. As a direct application, our singularity-encoded Green's function is used to reconstruct the solution 
	\begin{equation*}
		u(\bm{x}) = - e^{x_1^2 + x_2^2 - 1} + 1
	\end{equation*}
	through the numerical integration \eqref{BVP-NumSolu-Representation-SgEncdGreen}, bypassing the need to numerically solve our Poisson's equation \eqref{BVP-2D-Laplacian}. The pointwise error depicted in \autoref{fig-2D-PtErr-u} implies close agreement between the exact and reconstructed solutions.

	%-------------------------%
	% mannully displayed here, content belonging to the next section
	\begin{table}[!b]
		\centering
		\renewcommand\arraystretch{2} 
		\begin{tabular}{|c|c|c|c|c|}
			\hline
			mesh size & eigenvalues of preconditioned matrix & condition number & $\#$ iterations & $\ell_2$-norm error \\
			\hline
			$h = 0.1$ &
			\begin{minipage}[b]{0.3\columnwidth}
				\centering
				\raisebox{-.4\height}{\includegraphics[width=\linewidth]{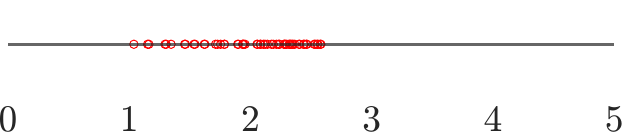}}			
			\end{minipage}
			& \makecell[c]{$\kappa (\widecheck{B}A) = 2.48$\\ $\kappa (A) = 20.55$} &\makecell[c]{$17$\\$(27)$} &\makecell[c]{$2.41 \times 10^{-14}$\\$(4.04 \times 10^{-14})$} \\
			\hline
			$h = 0.05$ &
			\begin{minipage}[b]{0.3\columnwidth}
				\centering
				\raisebox{-.4\height}{\includegraphics[width=\linewidth]{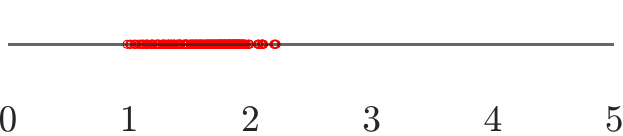}}
			\end{minipage}
			& \makecell[c]{$\kappa (\widecheck{B}A) = 2.25$ \\ $ \kappa (A) = 89.68$} &\makecell[c]{$17$\\$(59)$} &\makecell[c]{$3.01 \times 10^{-14}$\\$(3.34 \times 10^{-14})$}\\
			\hline
			$h = 0.025$ &
			\begin{minipage}[b]{0.3\columnwidth}
				\centering
				\raisebox{-.4\height}{\includegraphics[width=\linewidth]{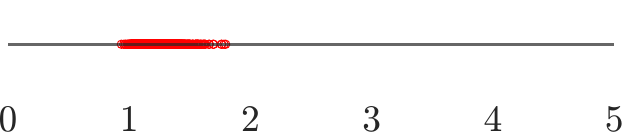}}
			\end{minipage}
			& \makecell[c]{$\kappa (\widecheck{B}A) = 1.92$ \\ $ \kappa (A) = 446.16$} &\makecell[c]{$16$\\$(129)$} &\makecell[c]{$2.51 \times 10^{-15}$\\$(6.06 \times 10^{-14})$}\\
			\hline
			$h = 0.0125$ &
			\begin{minipage}[b]{0.3\columnwidth}
				\centering
				\raisebox{-.4\height}{\includegraphics[width=\linewidth]{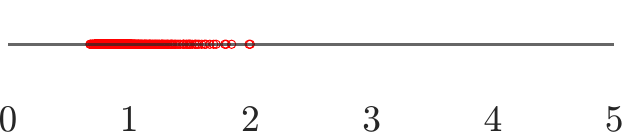}}
			\end{minipage}
			& \makecell[c]{$\kappa (\widecheck{B}A) = 2.96$ \\ $ \kappa (A) = 2041.57$} &\makecell[c]{$22$\\$(282)$} &\makecell[c]{$9.46 \times 10^{-15}$\\$(1.08 \times 10^{-14})$}\\
			\hline
		\end{tabular}
		\caption{Biconjugate gradient method with and without neural preconditioning for solving the discrete linear system of the two-dimensional problem \eqref{BVP-2D-Laplacian}.}
		\label{GreenPreconditioner-Poisson2DExp1}
	\end{table}
	%-------------------------%
		
	Particularly, eigenvalues of the exact Green's function \eqref{ExactGreen-2D-Laplacian} are available for this specific problem \cite{henrot2006extremum}, hence enabling a direct and rigorous comparison of those computed from our singularity-encoded Green’s function. By using linear elements with mesh width $h=0.045$ in formula \eqref{EigenProblem-SgEncdGreen}, we depict in \autoref{fig-2D-eigenvalues-1to50} and \autoref{fig-2D-eigenvalues-51to100} the exact eigenvalues and our approximated counterparts across a wide range of frequencies, along with their relative error in \autoref{fig-2D-eigenvalues-relerr}. Our numerical results reveal a clear manifestation of spectral bias \cite{xu2025understanding}, where low-frequency modes are reconstructed with high accuracy but high-frequency counterparts progressively deteriorates.
	
%%%%%%%%%%%%%%%%%%%%%%%%%%%%%%%%%%%%%%%%%%%%%%%%%%%%%

%%%%%%%%%%%%%%%%%%%%%%%%%%%%%%%%%%%%%%%%%%%%%%%%%%%%%
\subsubsection{Neural Preconditioning Matrix}
			
	Based on our pre-trained singularity-encoded Green's function, a neural preconditioner for the discrete system of problem \eqref{BVP-2D-Laplacian} can be constructed by applying Algorithm \ref{Alg-SgEncd-Green-Preconditioner}. The discretization is conducted using linear elements, while spectral properties of our preconditioned system are computationally evaluated and reported in \autoref{GreenPreconditioner-Poisson2DExp1}. As in the one-dimensional case (see section \ref{Section-1D-ToyModel-Preconditioner}), the discrete linear system is iteratively solved via the biconjugate gradient method \cite{golub2013matrix}, both with and without the use of our neural preconditioner. 
	
	As expected, the condition number of the original stiffness matrix grows with mesh refinement, leading to an increasing number of iterations to achieve convergence at machine precision. On the contrary, eigenvalues of our preconditioned matrix are notably clustered, allowing comparable accuracy to be achieved with fewer iterations.
	
%%%%%%%%%%%%%%%%%%%%%%%%%%%%%%%%%%%%%%%%%%%%%%%%%%%%%

%%%%%%%%%%%%%%%%%%%%%%%%%%%%%%%%%%%%%%%%%%%%%%%%%%%%%
\subsubsection{Hybrid Iterative Methods}	

	%--------------- 2DLaplacian-Example iterative methods ---------------%
	% mannully displayed here, originally in the next subsection, command from Qi
	\begin{figure}[!t]
		\centering
		\begin{subfigure}{0.47\textwidth}
			\centering
			\includegraphics[width=1\textwidth]{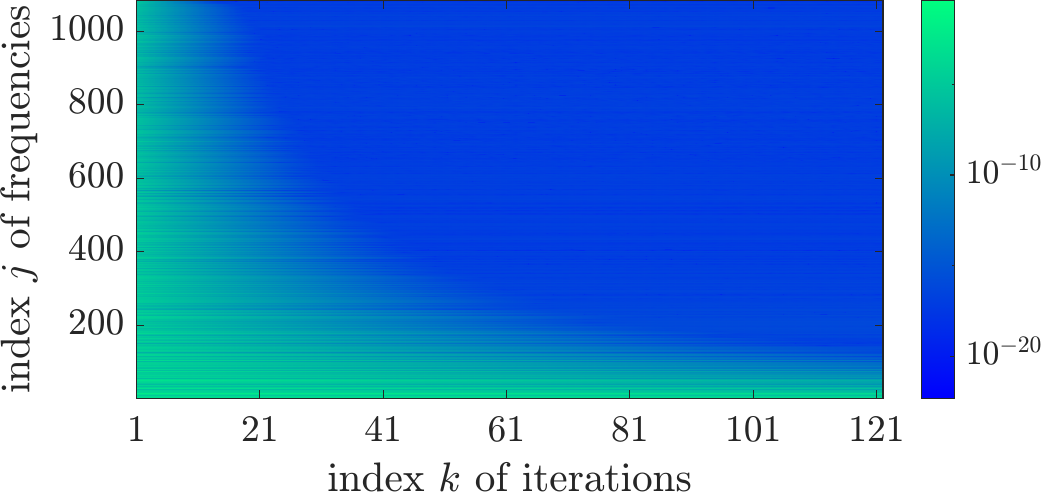}
			\caption{the spectral error of damped Jacobi method \eqref{1D-Poisson-Damped-Jacobi}}
			\label{fig-2D-Poisson-Jacobi-ModeErr}
		\end{subfigure}
		\hspace{0.2cm}		
		\begin{subfigure}{0.47\textwidth}
			\centering
			\includegraphics[width=1\textwidth]{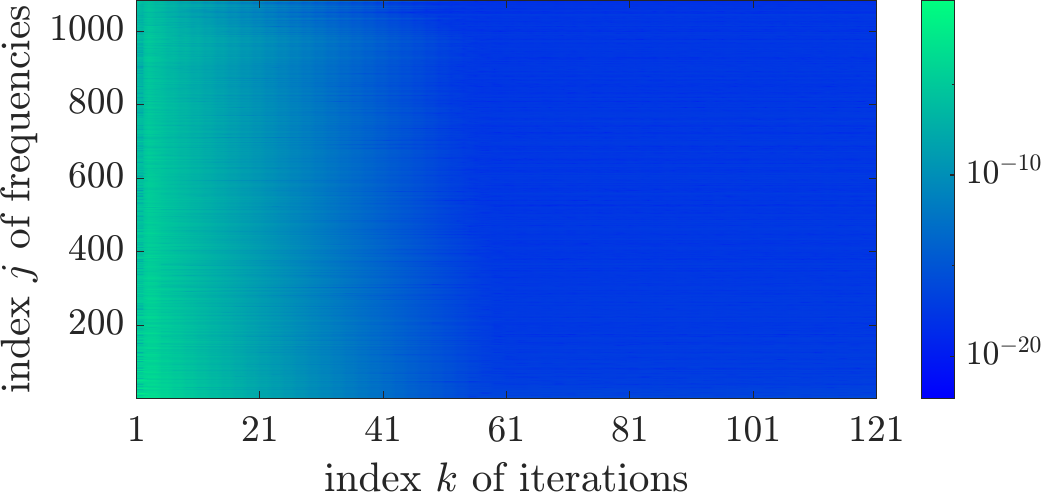}
			\caption{the spectral error of hybrid iterative method \eqref{1D-Poisson-Hybrid-IteM}}
			\label{fig-2D-Poisson-Hybrid-ModeErr}
		\end{subfigure}
	
		\centering
		\vspace{0.3cm}
		\begin{subfigure}{0.31\textwidth}
			\centering
			\includegraphics[width=1\textwidth]{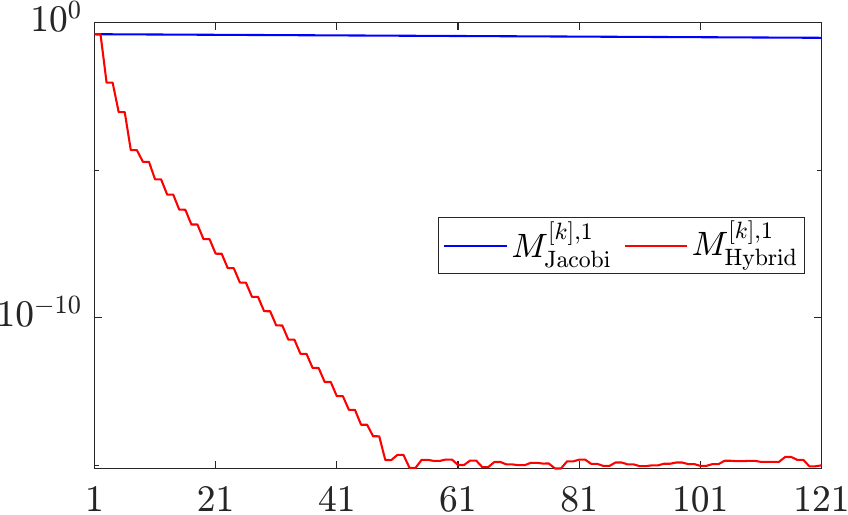}
			\caption{the lowest-frequency error component}
			\label{fig-2D-Poisson-Hybrid-ModeErr-j1} 
		\end{subfigure}
		\hspace{0.1cm}
		\begin{subfigure}{0.31\textwidth}
			\centering
			\includegraphics[width=1\textwidth]{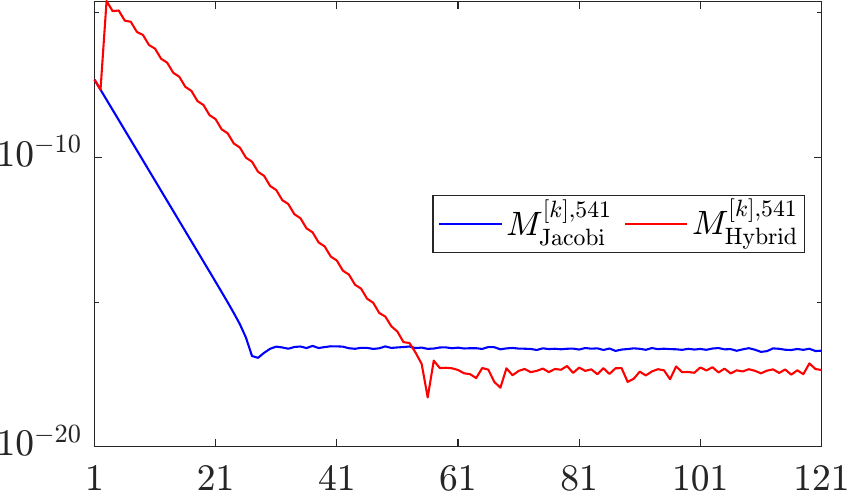}
			\caption{the mid-frequency error component}
			\label{fig-2D-Poisson-Hybrid-ModeErr-j127} 
		\end{subfigure}
		\hspace{0.1cm}
		\begin{subfigure}{0.31\textwidth}
			\centering
			\includegraphics[width=1\textwidth]{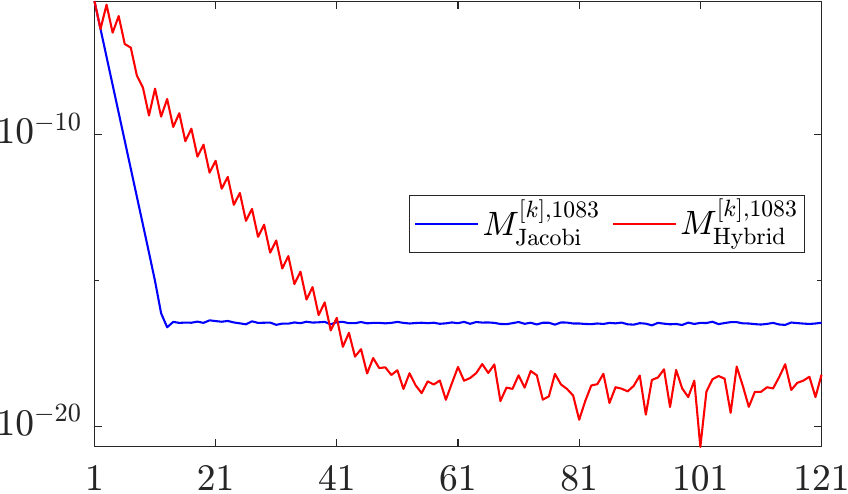}
			\caption{the highest-frequency error component}
			\label{fig-2D-Poisson-Hybrid-ModeErr-j255} 
		\end{subfigure}
		
		\vspace{-0.1cm}		
	
		\caption{Comparison of spectral errors for solving the discrete linear system of the two-dimensional problem \eqref{BVP-2D-Laplacian}.}
		\label{Experiments-Poisson2D-IteM-1}
		\vspace{-0.4cm}
	\end{figure}
	%--------------- 2DLaplacian-Example iterative methods ---------------%
	
\begin{figure}[!t]	
		\centering
		\begin{subfigure}{0.46\textwidth}
			\centering
			\includegraphics[width=1\textwidth]{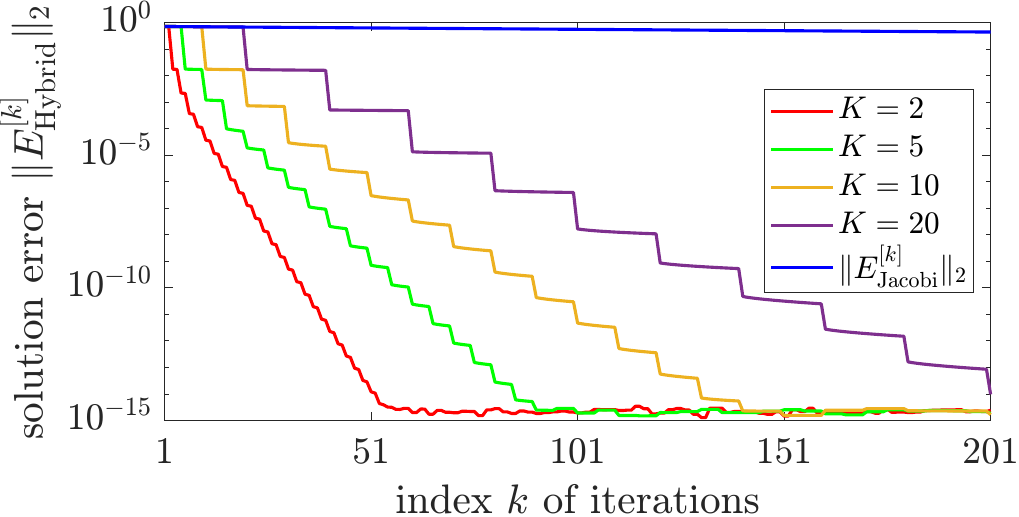}
			\caption{comparison of the solution error over iterations}
			\label{fig-2D-Poisson-Solu-Err} 
		\end{subfigure}
		\hspace{0.3cm}
		\begin{subfigure}{0.46\textwidth}
			\centering
			\includegraphics[width=1\textwidth]{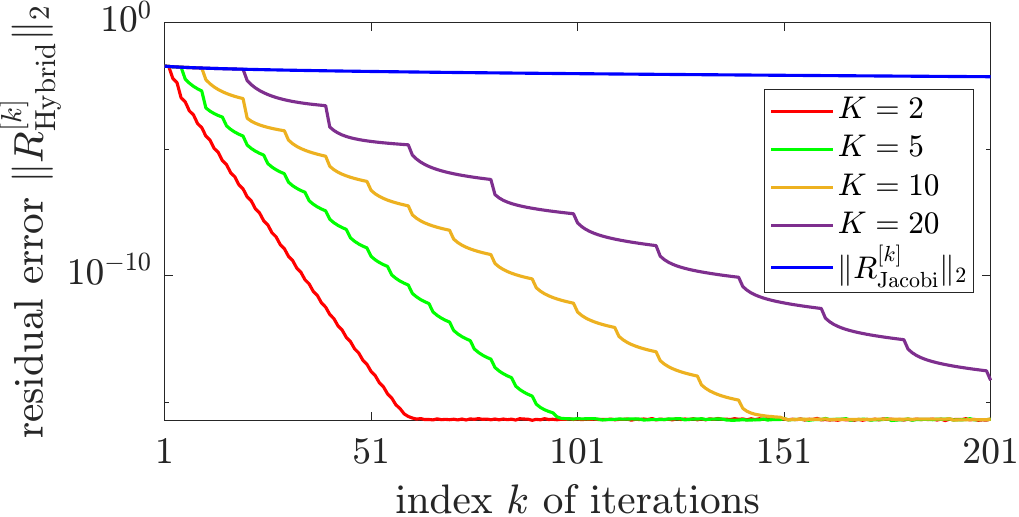}
			\caption{comparison of the residual error over iterations}
			\label{fig-2D-Poisson-Res-Err}
		\end{subfigure}
		
		\vspace{-0.1cm}		
	
		\caption{Hybrid iteration \eqref{1D-Poisson-Hybrid-IteM} with different switching periods for the discrete linear system of the two-dimensional problem \eqref{BVP-2D-Laplacian}.}
		\label{Experiments-Poisson2D-IteM-2}
		\vspace{-0.4cm}
	\end{figure}

	To exploit the spectral bias inherent in our singularity-encoded Green's function, the hybrid iterative scheme \eqref{1D-Poisson-Hybrid-IteM} is applied to solve the discrete system of \eqref{BVP-2D-Laplacian} with mesh size $h=0.025$, followed by the comparison with that of the damped Jacobi method \eqref{1D-Poisson-Damped-Jacobi}. As evident in \autoref{fig-2D-Poisson-Jacobi-ModeErr} and \autoref{fig-2D-Poisson-Hybrid-ModeErr}, the damped Jacobi method is effective in attenuating high-frequency errors but exhibits limited efficacy in eliminating low-frequency errors, as reported in section \ref{Section-1D-ToyModel-Hybrid-IteM}. On the contrary, our hybrid method provides rapid error reduction across all frequencies, achieving machine accuracy within approximately 50 number of iterations. Moreover, \autoref{fig-2D-Poisson-Solu-Err} and \autoref{fig-2D-Poisson-Res-Err} suggest that $K=2$ remains the most efficient choice of switching period, consistent with our findings in the one dimension.

	%--------------- 2DLaplacian-Example iterative methods ---------------%
	%--------------- 2DLaplacian-Example iterative methods ---------------%
	
\begin{figure}[H]
		\centering
		\begin{subfigure}{0.5\textwidth}
			\centering
			\includegraphics[width=\textwidth]{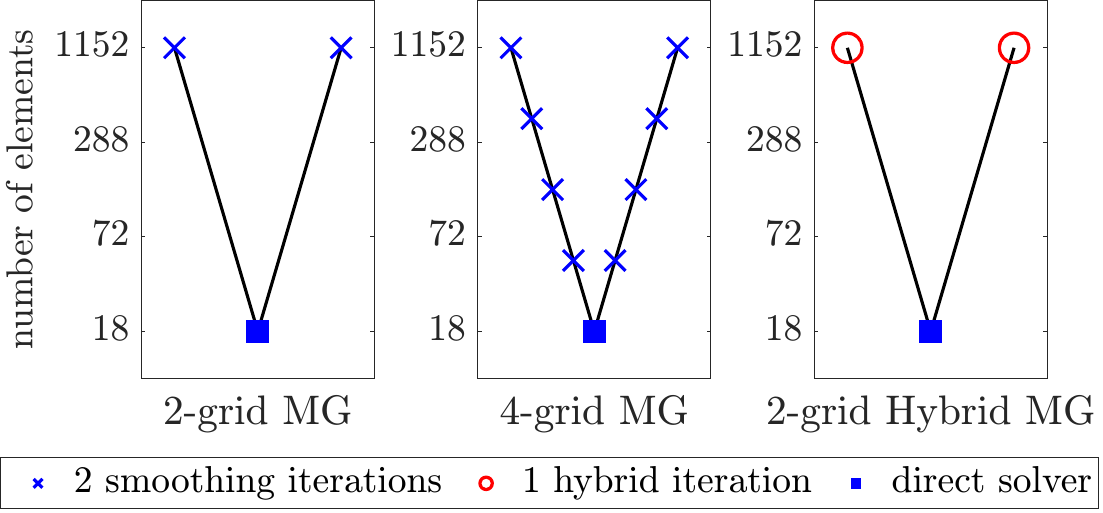}
			\caption{V-cycle structures for different approaches}
			\label{fig-2D-Poisson-MG-V-Cycle}
		\end{subfigure}
		\hspace{0.25cm}
		\begin{subfigure}{0.45\textwidth}
			\centering
			\includegraphics[width=\textwidth]{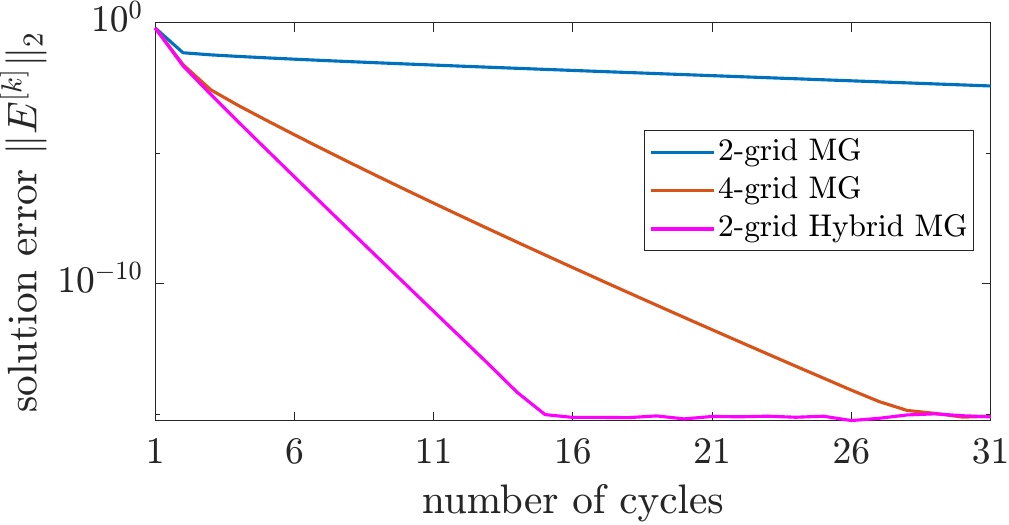}
			\caption{solution errors over iterative cycles}
			\label{fig-2D-Poisson-MG-Solu-Err}
		\end{subfigure}
		
		\vspace{-0.1cm}		
	
		\caption{Comparison of multigrid methods with and without hybrid iterations for solving the discrete system of the two-dimensional problem \eqref{BVP-2D-Laplacian}.}
		\label{Experiment-Poisson2D-HybridMG}
		\vspace{-0.4cm}
	\end{figure}

	Next, a comparative study between classical multigrid methods and our hybrid-smoothing variants are shown in \autoref{fig-2D-Poisson-MG-V-Cycle} and \autoref{fig-2D-Poisson-MG-Solu-Err}. Remarkably, the two-grid variant achieves a faster convergence than the traditional four-grid scheme, reducing the reliance on deeply nested grids without compromising solution accuracy.
		
	%---------------------------%
	%---------------------------%

%%%%%%%%%%%%%%%%%%%%%%%%%%%%%%%%%%%%%%%%%%%%%%%%%%%%%
%%%%%%%%%%%%%%%%%%%%%%%%%%%%%%%%%%%%%%%%%%%%%%%%%%%%%

%%%%%%%%%%%%%%%%%%%%%%%%%%%%%%%%%%%%%%%%%%%%%%%%%%%%%
%%%%%%%%%%%%%%%%%%%%%%%%%%%%%%%%%%%%%%%%%%%%%%%%%%%%%
\section{Conclusions}\label{Sec-Conclusion}
	
	In this paper, we propose a singularity-encoded learning framework to reconstruct the Green’s function of elliptic boundary value problems without supervised training data. Based on the estimate of its singularities, our Green's function is embedded into a one-order higher-dimensional space by incorporating an augmented variable, leading to a numerically tractable surrogate with improved regularity. This enables us to reconstruct the Green's function using fully-connected neural networks with smooth activation functions, while the sampling of collocation points remains unchanged regardless of the increased dimensionality. After projecting the trained neural network model back onto its physical domain, the spectral bias inherent to our singularity-encoded Green's function is applied to accelerate classical iterative solvers, either as an effective preconditioner or as a correction step to eliminate low- frequency errors. Numerical studies across benchmark problems are reported to demonstrate the effectiveness and efficiency of our proposed methods. 
	
	Our proposed methods offer a promising framework for integrating deep learning techniques with traditional numerical schemes and the rigorous theoretical analysis will be the focus of future research. It is also noteworthy that the proposed singularity-encoded technique admits a natural extension to a broad class of singular problems, a direction that will be explored in our subsequent work.
	
%%%%%%%%%%%%%%%%%%%%%%%%%%%%%%%%%%%%%%%%%%%%%%%%%%%%%
%%%%%%%%%%%%%%%%%%%%%%%%%%%%%%%%%%%%%%%%%%%%%%%%%%%%%	

\appendix
\section{Green's Function for Poisson's Equation in Two Dimension}\label{Appendix-1}
	
	Unlike the logarithmic augmented variable \eqref{Singularity-SgEncdGreen}, we now consider an alternative approach of the form $\varphi(\bm{x},\bm{y}) = \lVert \bm{x} - \bm{y} \rVert ^{-0.2}$ for \eqref{BVP-ExactGreen-Encoded}. Numerical results for the four-dimensional problem \eqref{BVP-2D-Laplacian-ExactGreen} are displayed in \autoref{Appendix-2D-Poisson}, which imply that our singularity-encoded Green's function achieves satisfactory accuracy. Specifically, the realizations of our singularity-encoded Green's function at fixed points $\bm{y}=(0,0)$ and $\bm{y}=(0.7,0.2)$ are depicted in \autoref{fig-Appendix-SgEncdGreen-2D-Poisson-surf-x0} and \autoref{fig-Appendix-SgEncdGreen-2D-Poisson-surf-x1}, along with their error profiles in \autoref{fig-Appendix-Green-2D-Poisson-PtErr-x0} and \autoref{fig-Appendix-Green-2D-Poisson-PtErr-x1}. Moreover, the comparison between exact and numerical eigenvalues are presented in \autoref{fig-Appendix-2D-eigenvalues-1to50}, \autoref{fig-Appendix-2D-eigenvalues-51to100}, and \autoref{fig-Appendix-2D-eigenvalues-relerr}, demonstrating the spectral bias of our trained model.
	
% The derivatives $\nabla_{\bm{x}}\widehat{G}(\bm{x},\bm{y},\varphi(\bm{x},\bm{y}))$ and  $\nabla_{\bm{x}}\widehat{G}(\bm{x},\bm{y},\varphi(\bm{x},\bm{y}))$ remain bounded despite the stronger singularity of $\varphi(\bm{x},\bm{y})$ than $\ln\lVert\bm{x}-\bm{y}\rVert$ in \eqref{Singularity-ExactGreen} as $\bm{x}\to\bm{y}$. 	
%	Here, our singularity-encoded Green's function satisfies, for any $\bm{x}\in\Omega$, 	
%	\begin{equation}
%		\begin{array}{cl}
%			-\Delta_{\bm{x}}\widehat{G}(\bm{x},\bm{y},\varphi(\bm{x},\bm{y})) + 0.4 \nabla_{\bm{x}} \big( \partial_z\widehat{G}(\bm{x},\bm{y},\varphi(\bm{x},\bm{y})) \big) \cdot \dfrac{\bm{x}-\bm{y}}{\lVert \bm{x} - \bm{y} \rVert ^{2.2}}  - 0.04 \dfrac{\partial_{zz}\widehat{G}(\bm{x},\bm{y},\varphi(\bm{x},\bm{y}))}{\lVert \bm{x} - \bm{y} \rVert ^{2.4}} = 0, \ \ \ &\ \textnormal{for}\ \bm{x}\in\Omega\setminus\Gamma, \\
%			\displaystyle -\int_{\partial B_\epsilon (\bm{y})} \bigg(\nabla_{\bm{x}}G(\bm{x},\bm{y},\varphi(\bm{x},\bm{y})) \cdot \dfrac{\bm{x}-\bm{y}}{\lVert \bm{x} - \bm{y} \rVert} - 0.2 \dfrac{\partial_z\widehat{G}(\bm{x},\bm{y},\varphi(\bm{x},\bm{y}))}{\lVert \bm{x} - \bm{y} \rVert ^ {1.2}}\bigg) d S(\bm{x}) = 1, \ \ \ &\ \textnormal{for} \ \bm{x}\in\Gamma,\\
%			G(\bm{x},\bm{y},\varphi(\bm{x},\bm{y})) = 0, \ \ \ &\ \textnormal{for}\ \bm{x}\in\partial\Omega,
%		\end{array}
%	\end{equation}

	%--------------------------------------%
	\begin{figure}[!ht]
		\centering
		\begin{subfigure}{0.31\textwidth}
			\centering
			\includegraphics[width=\textwidth]{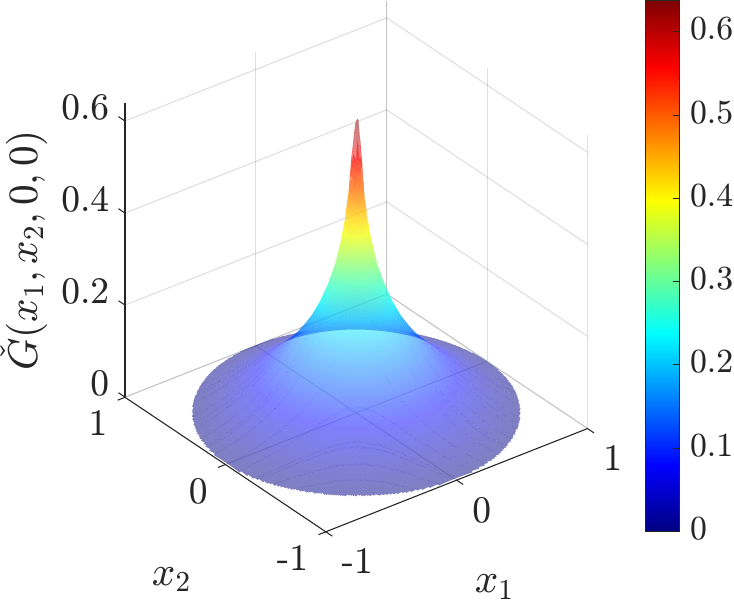}
			\caption{$\widecheck{G}(x_1,x_2,0,0)$}
			\label{fig-Appendix-SgEncdGreen-2D-Poisson-surf-x0}
		\end{subfigure}
		\hspace{0.2cm}
		\begin{subfigure}{0.294\textwidth}
			\centering
			\includegraphics[width=\textwidth]{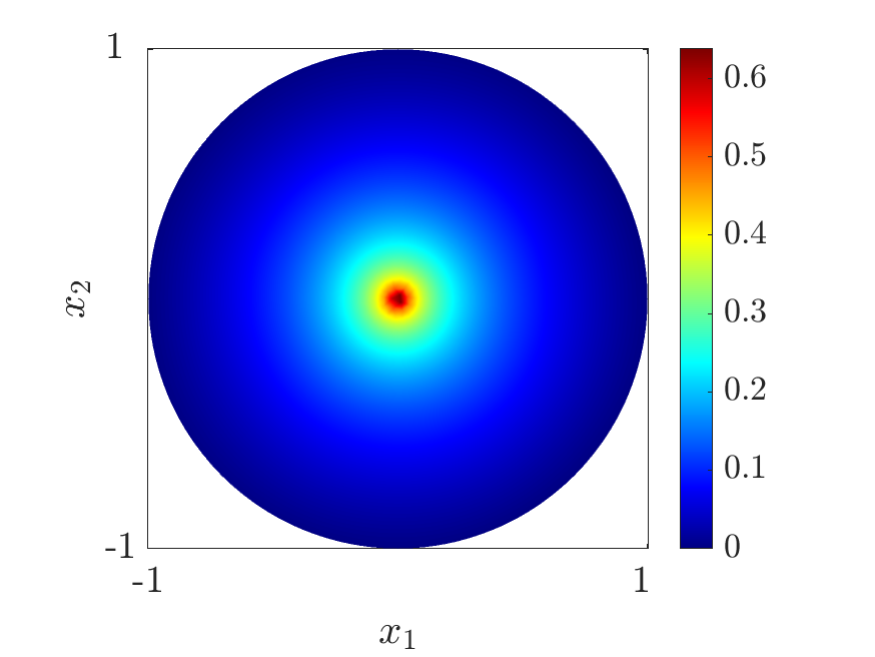}
			\caption{$\widecheck{G}(x_1,x_2,0,0)$}
			\label{fig-Appendix-SgEncdGreen-2D-Poisson-x0}
		\end{subfigure}
		\hspace{0.2cm}
		\begin{subfigure}{0.31\textwidth}
			\centering
			\includegraphics[width=\textwidth]{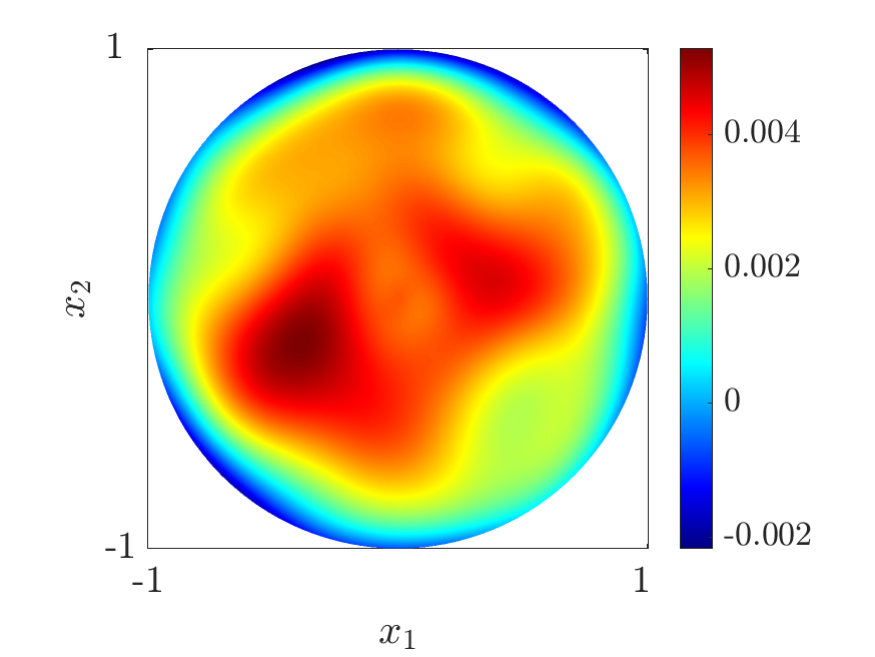}
			\caption{$G(x_1,x_2,0,0) - \widecheck{G}(x_1,x_2,0,0)$}
			\label{fig-Appendix-Green-2D-Poisson-PtErr-x0}
		\end{subfigure}
		
		\vspace{0.1cm}
		
		\begin{subfigure}{0.31\textwidth}
			\centering
			\includegraphics[width=\textwidth]{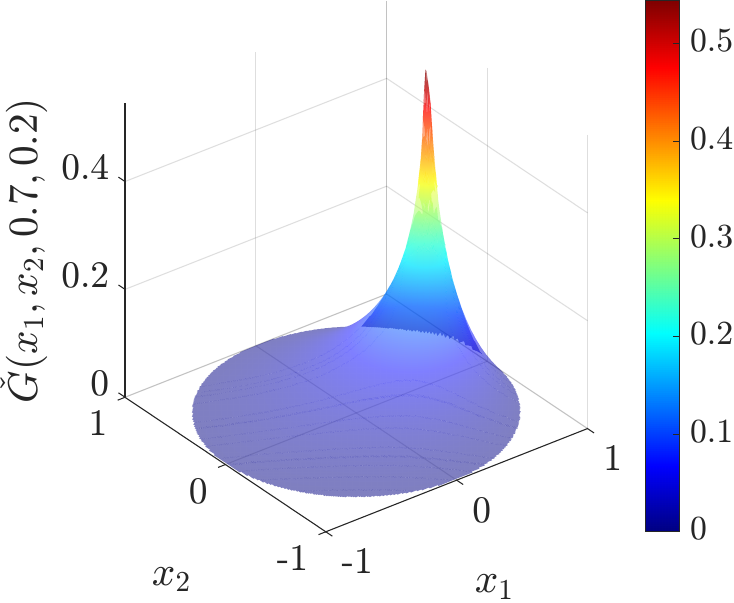}
			\caption{$\widecheck{G}(x_1,x_2,0.7,0.2)$}
			\label{fig-Appendix-SgEncdGreen-2D-Poisson-surf-x1}
		\end{subfigure}
		\hspace{0.3cm}
		\begin{subfigure}{0.294\textwidth}
			\centering
			\includegraphics[width=\textwidth]{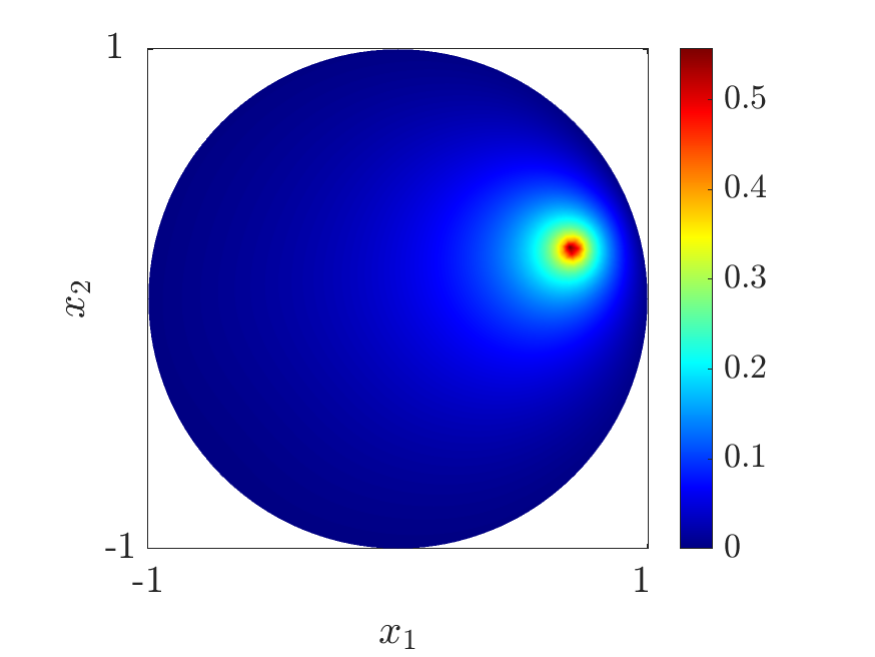}
			\caption{$\widecheck{G}(x_1,x_2,0.7,0.2)$}
			\label{fig-Appendix-SgEncdGreen-2D-Poisson-x1}
		\end{subfigure}
		\hspace{0.2cm}
		\begin{subfigure}{0.31\textwidth}
			\centering
			\includegraphics[width=\textwidth]{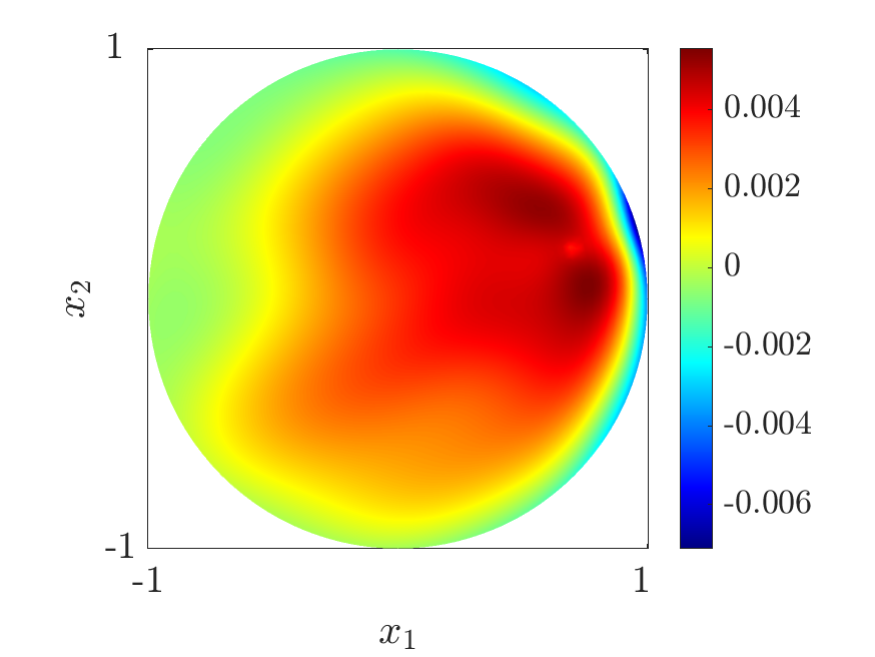}
			\caption{$G(x_1,x_2,0.7,0.2) - \widecheck{G}(x_1,x_2,0.7,0.2)$}
			\label{fig-Appendix-Green-2D-Poisson-PtErr-x1}
		\end{subfigure}
		
		\vspace{0.1cm}
		
		\begin{subfigure}{0.31\textwidth}
			\centering
			\includegraphics[width=\textwidth]{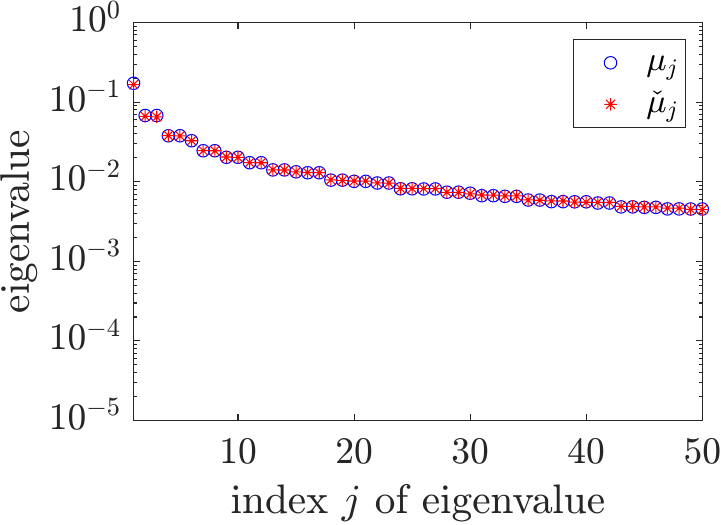}
			\caption{$\mu_j$, $\check{\mu}_j$ with $1\leq j\leq 50$}
			\label{fig-Appendix-2D-eigenvalues-1to50}
		\end{subfigure}
		\hspace{0.1cm}
		\begin{subfigure}{0.31\textwidth}
			\centering
			\includegraphics[width=\textwidth]{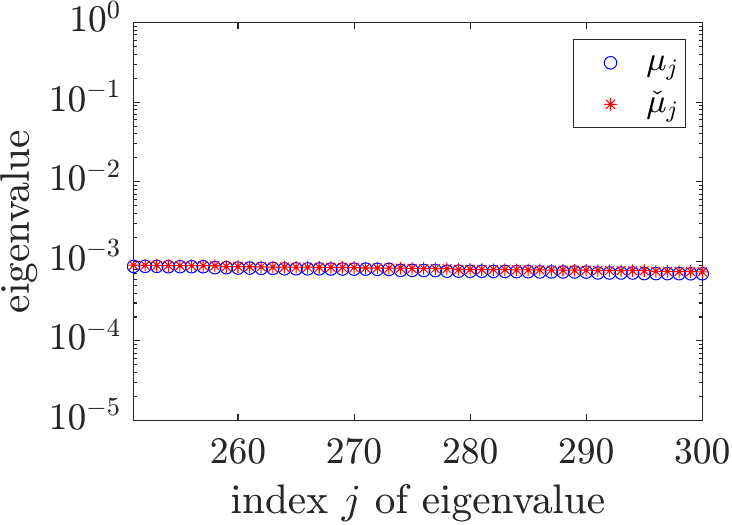}
			\caption{$\mu_j$, $\check{\mu}_j$ with $251\leq j\leq 300$}
			\label{fig-Appendix-2D-eigenvalues-51to100}
		\end{subfigure}
		\hspace{0.1cm}
		\begin{subfigure}{0.31\textwidth}
			\centering
			\includegraphics[width=\textwidth]{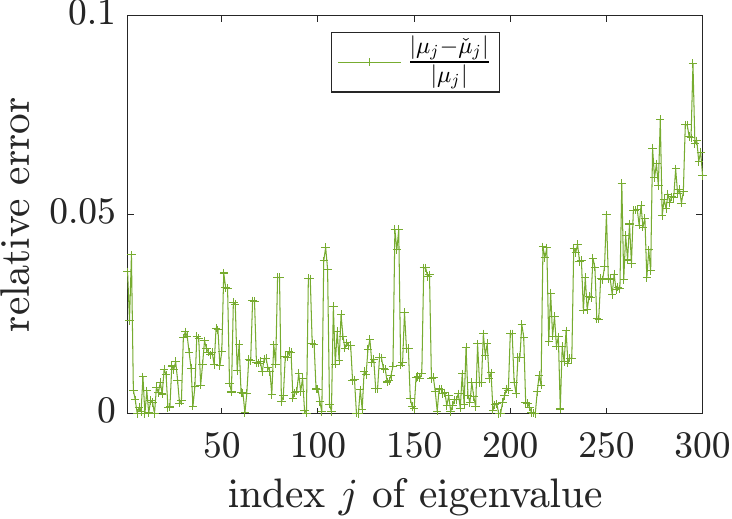}
			\caption{$\delta_{\mu_j}$ with $1\leq j\leq 300$}
			\label{fig-Appendix-2D-eigenvalues-relerr}	
		\end{subfigure}	
		
		\caption{Numerical validation of our singularity-encoded Green's function for the two-dimensional problem \eqref{BVP-2D-Laplacian}.}
		\label{Appendix-2D-Poisson}
	\end{figure}
	%--------------------------------------%
	
%%%%%%%%%%%%%%%%%%%%%%%%%%%%%%%%%%%%%%%%%%%%%%%%%%%%%
%%%%%%%%%%%%%%%%%%%%%%%%%%%%%%%%%%%%%%%%%%%%%%%%%%%%% 	




\bibliographystyle{elsarticle-num}
\bibliography{refs}

\end{document}